\documentclass[11pt]{article}
\usepackage{float}
\usepackage[margin=1in]{geometry}
\usepackage{amsmath,amssymb,amsthm}
\usepackage{bbm,enumitem,color,comment}
\usepackage{graphicx}
\usepackage{subcaption}
\usepackage{amsfonts} 
\usepackage{pdfpages}
\usepackage{dsfont}
\usepackage{caption}
\usepackage{mathrsfs}
\numberwithin{equation}{section}
\usepackage{dsfont}

\usepackage{ragged2e}
\usepackage{algorithm}
\usepackage{algpseudocode}
\usepackage{lastpage}
\usepackage{fancyhdr}
\usepackage[utf8]{inputenc}
\usepackage[numbers,sort&compress]{natbib}
\usepackage[title]{appendix}
\usepackage{indentfirst}
\usepackage{float}
\usepackage[section]{placeins}
\usepackage{hyperref}
\hypersetup{
    colorlinks=false,
    pdfborder={0 0 1},
    citebordercolor={0 1 0}
}
\begin{document}
\newtheorem{assumption}{Assumption}
\newtheorem{theorem}{Theorem}[section]
\newtheorem{corollary}[theorem]{Corollary}
\newtheorem{proposition}[theorem]{Proposition}
\newtheorem{lemma}[theorem]{Lemma}
\newtheorem{definition}[theorem]{Definition} 
\newtheorem{remark}[theorem]{Remark}

\newtheorem{stabilityrestatement}{Theorem}
\renewcommand{\thestabilityrestatement}{%
  \ref{thm:quantitative-stability}}
\newtheorem{convergencerestatement}{Theorem}
\renewcommand{\theconvergencerestatement}{\ref{thm:convergence-to-mot}}
\newtheorem{example}{Example}

\title{Entropic Martingale Optimal Transport with Wasserstein Relaxed Marginal Constraints}
\author{
Haowen Chen\thanks{
Department of Finance and Risk Engineering,
Tandon School of Engineering, New York University.
Email: \texttt{hc4698@nyu.edu}.
}
\and
Renyuan Xu\thanks{
Department of Management Science and Engineering,
Stanford University.
Email: \texttt{renyuanxu@stanford.edu}.
}
\and
Xin Zhang\thanks{
Department of Finance and Risk Engineering,
Tandon School of Engineering, New York University.
Email: \texttt{xz1662@nyu.edu}.
}
}
\maketitle

\begin{abstract}
Martingale optimal transport (MOT) is a fundamental framework in robust finance, yet its stability under marginal perturbation no longer holds in high dimensions $d\geq 2$. In this paper, we introduce Wasserstein-relaxed entropic martingale optimal transport (WEMOT), which combines Wasserstein relaxation of the second marginal with  entropy regularization. This formulation ensures stability under perturbations of the marginals, thereby facilitating numerical implementation. We establish strong duality and dual attainment for WEMOT, which leads to efficient Sinkhorn-type algorithms. 
We further prove a   convergence rate of WEMOT to MOT as regularization parameters vanish, and also the quantitative stability of WEMOT with respect to marginal perturbations. 
\end{abstract}

\tableofcontents

\section{Introduction}
\label{sec:introduction}

Martingale optimal transport (MOT) provides a powerful framework for robust pricing and hedging \cite{beiglbock2016problem, beiglbock2013model, beiglbock2017complete}. For $\mu,\nu\in\mathcal P_2(\mathbb R^d)$ in convex order and a cost function $c$, MOT is defined via 
\[
M(\mu,\nu):=\inf_{\pi\in\Pi_M(\mu,\nu)}\int_{\mathbb R^d\times\mathbb R^d}c(x,y)\,\pi(dx,dy),
\]
where $\Pi_M(\mu,\nu)$ denotes the set of martingale couplings $\pi$ from $\mu$ to $\nu$. Here a coupling $\pi$ is said to be a martingale coupling if there exist a martingale $(X,Y)$, i.e. $\!\mathbb E[Y| X]\!\!=\!\!X$, whose joint distribution is $\pi$.   In mathematical finance, this martingale constraint encodes the absence of arbitrage thanks to the fundamental theorem of asset pricing \citep{delbaen1994general}, and MOT yields model-independent price bounds with a dual representing semi-static hedging  \cite{beiglbock2013model}. To facilitate computation, entropic martingale optimal transport (EMOT) adds relative entropy regularization to MOT cost functional \cite{nutz2023martingale, NutzWiesel2024,backhoff2026bridging}, allowing for Sinkhorn-type algorithms \cite{demarch2018entropic, KazeykinaRenWang2026}.

However, the martingale constraint makes both MOT and EMOT unstable to marginal perturbations. By Strassen's theorem, $\Pi_M(\mu,\nu)\neq\varnothing$ if and only if $\mu\preceq_{cx}\nu$ \cite{strassen1965existence}, i.e. $\int f d\mu \leq \int f d\nu$ for any convex function $f$. Discretizations $\mu^N\to\mu$ and $\nu^N\to\nu$ may violate the convex order condition and thus make both problems infeasible. Even when convex order between $\mu^N, \nu^N$ is preserved, continuity of $M(\mu^N, \nu^N)$ to $M(\mu, \nu)$ does not hold in general. Stability on the real line was established by \citep{backhoff2022stability} and \citep{wiesel2023continuity}. In dimensions $d\geq2$, \citet{bruckerhoff2022instability} construct convergent marginal sequences for which both the value and the optimizers of MOT fail to converge. The same construction also yields the instability of  EMOT.

We address the instability of MOT and EMOT by relaxing the second marginal constraint. For $\lambda,\epsilon>0$, we introduce the Wasserstein-relaxed entropic martingale optimal transport  (WEMOT):
\begin{equation}
M^{\lambda,\epsilon}(\mu,\nu):=\inf_{\substack{\eta\in\mathcal P_2(\mathbb R^d)\\ \pi\in\Pi_M(\mu,\eta)}}\biggl\{\int c\,d\pi+\epsilon H(\pi\mid\mu\otimes\eta)+\lambda W_2^2(\eta,\nu)\biggr\}.
\label{eq:intro-WEMOT}
\end{equation}
 The Wasserstein relaxation allows the second marginal to vary, while entropic regularization yields a Gibbs-type density of the optimal martingale coupling relative to $\mu\otimes\eta$.
 The optimization problem \eqref{eq:intro-WEMOT} is well-posed for arbitrary $\mu,\nu\in\mathcal P_2(\mathbb R^d)$ without requiring convex order. At $\epsilon=0$, it reduces to Wasserstein-penalized MOT, which fits the general approach of penalizing marginal constraints developed by \citet{doldi2023entropy}. Our Wasserstein-Entropy regularization is conceptually related to the Gamma-divergence of \citep{dupuis2022divergence}, which is defined as an inf-convolution of relative entropy and a transport-type discrepancy.

Our first main result is a quantitative stability estimate with respect to both marginals in arbitrary dimension. Under the Lipschitz assumption on $c$, Theorem~\ref{thm:quantitative-stability} gives an explicit rate
\[
\bigl|M^{\lambda,\epsilon}(\widetilde\mu,\widetilde\nu)-M^{\lambda,\epsilon}(\mu,\nu)\bigr|=O\bigl(W_2(\mu,\widetilde\mu)^{2/3}+W_2(\nu,\widetilde\nu)\bigr).
\]
The asymmetry in the rates for the two marginals arises as only the second marginal constraint is relaxed. And the $2/3$ exponent in $W_2(\mu,\widetilde\mu)^{2/3}$ is due to the Gaussian smoothing construction in the proof, which we do not claim that it is optimal, while the Lipschitz rate in the target marginal is sharp. For fixed $\lambda$ and $\epsilon$, the estimate ensures the continuity of WEMOT under marginal approximations even when the approximating marginals are not in convex order.

We next show that WEMOT approximates classical MOT as the regularization parameters vanish. For $\mu\preceq_{cx}\nu$,  Theorem~\ref{thm:convergence-to-mot} establishes the convergence
\[
\lim_{\lambda\to\infty}\lim_{\epsilon\to0^+}M^{\lambda,\epsilon}(\mu,\nu)=M(\mu,\nu)\]
Under the additional geometric and compact-support assumptions of $\mu, \nu$, we obtain quantitative estimate
\[
\bigl|M^{\lambda,\epsilon}(\mu,\nu)-M(\mu,\nu)\bigr|\leq C\biggl\{\frac1\lambda+\epsilon\log(1+\lambda)\biggr\},
\]
where $C$ is a positive constant  independent of $\lambda$ and $\epsilon$. The proof combines Gaussian smoothing, which produces finite-entropy martingale couplings, with a uniform bound on  dual potentials of MOT established in \citet{carlier2026weak}.

We also establish strong duality and dual attainment. For bounded continuous costs and arbitrary $\mu,\nu\in\mathcal P_2(\mathbb R^d)$, we prove
\[
M^{\lambda,\epsilon}(\mu,\nu)=D^{\lambda,\epsilon}(\mu,\nu):=\sup_{(f,h)\in\mathcal D}\Bigl\{\mu(f)+\nu\bigl(Q_\lambda\Gamma^{\epsilon}_{f,h}\bigr)\Bigr\},
\]
where $\mathcal D$ is the class of admissible dual potentials. In the dual formulation, we use the logarithmic transform $\Gamma^{\epsilon}_{f,h}$ defined in \eqref{eq:gamma-transform}  associated with entropy and the quadratic infimal convolution $Q_\lambda$ defined in \eqref{eq:Q transform} associated with the Wasserstein penalty. Dual attainment holds under Assumption~\ref{ass:dual-attainment}, which requires compact support of the first marginal $\mu$ but not of the second marginal $\nu$.

Our proof of WEMOT dual attainment relies on the dual attainment for  EMOT problem with marginals $(\mu,\eta^*)$, where $\eta^*$ is the optimal relaxed second marginal in \eqref{eq:intro-WEMOT}. Assumption~\ref{ass:dual-attainment} ensures the existence of a reference martingale coupling whose density relative to $\mu\otimes\eta^*$ is bounded above and away from zero, and Proposition~\ref{prop:emot-dual-attainment} yields an EMOT dual optimizer under this bounded density condition. In the existing literature, \citet{NutzWiesel2024} establish dual attainment for EMOT with $c=0$ in one dimension. In arbitrary dimension, \citet{KazeykinaRenWang2026, carlier2026weak} establish EMOT dual attainment requiring both marginals $\mu,\nu$ are compactly supported.  Our result allows the second marginal to have noncompact support in arbitrary dimension. This is particularly relevant to WEMOT, where the relaxed marginal $\eta$ may not be compactly supported.

For the dual attainment of $\operatorname{EMOT}(\mu,\eta^*)$ in \eqref{eq:emot}, we develop a novel Galerkin method that approximates the infinite-dimensional martingale constraint by finite dimensional constraints. The existence of bounded-density reference coupling ensures dual attainment for the approximating problems. Then, by a compactness argument, we prove the convergence of the dual optimizers, obtaining an admissible dual optimizer for the original problem. Extending dual attainment from EMOT to WEMOT requires an additional argument. The key step is to show that the dual potential of  $\operatorname{EMOT}(\mu,\eta^*)$ corresponding to $\eta^*$ coincides with the dual potential  of $W_2^2(\eta^*,\nu)$ corresponding to $\eta^*$ after a normalization.

The stability of WEMOT facilitates more reliable numerical implementation, which typically relies on discrete approximations of the marginals. We implement numerical experiments to illustrate the stability of WEMOT under perturbations of  marginals and its convergence to classical MOT. We compare the stability of WEMOT and MOT with sampled marginals. Figure~\ref{fig:mot-WEMOT-stability} shows lower standard deviation in the WEMOT comparing to MOT, indicating reduced sensitivity to marginal perturbations caused by sampling. Moreover, Figure~\ref{fig:stability-rate} provides numerical evidence that the quantitative stability rate in Theorem~\ref{thm:quantitative-stability} is close to sharp. For the  counterexample for MOT stability in \citet{bruckerhoff2022instability}, Figure~\ref{fig:Bruckerhoff_joint_stability_WEMOT} shows that the WEMOT values approach the limiting MOT value, while the MOT values  diverge from it. To assess the computational effect of entropy regularization, we propose a Sinkhorn-type algorithm for WEMOT and compare its convergence speed with that of Wasserstein-relaxed MOT. Figure~\ref{fig:computational-convergence} shows that entropy regularization accelerates numerical convergence in both iteration count and runtime in  examples considered.

The paper is organized as follows. Section~\ref{sec:problem-formulation} introduces the notation and the WEMOT formulation, and Section~\ref{sec:main-results} states the main results. Section~\ref{sec:proof-strong-duality} proves strong duality. Section~\ref{sec:primal and dual attainment} establishes  existence and uniqueness of the primal optimizer and proves dual attainment. Section~\ref{sec:quantitative_stability} develops the stability and convergence estimates, and Section~\ref{sec:numerics} presents numerical experiments. The proof of the EMOT dual-attainment result is given in Appendix~\ref{app:emot-dual-attainment}.

\section{Problem Formulation}
\label{sec:problem-formulation}
\subsection{Notation and preliminaries}\label{subsec:notations}

Fix $d\in\mathbb N$, and let $|\cdot|$ denote the Euclidean norm on
$\mathbb R^d$. For a matrix $A\in\mathbb R^{d\times d}$, write $\|A\|_{\mathrm{op}}:=\sup_{|v|=1}|Av|$ for the induced operator norm. When $A$ is symmetric, let $\theta_{\min}(A)$ denote its smallest eigenvalue. For a Polish space $E$, let $\mathcal P(E)$ denote the
set of Borel probability measures on $E$, and set $
\mathcal P_2(\mathbb R^d)
:=\{
\mu\in\mathcal P(\mathbb R^d):
\int_{\mathbb R^d}|x|^2\,\mu(dx)<\infty\}.$ For $\mu\in\mathcal P_2(\mathbb R^d)$, denote its mean by $m_\mu:=\int_{\mathbb R^d}x\,\mu(dx)$ and its second moment by $m_2(\mu):=\int_{\mathbb R^d}|x|^2\,\mu(dx)$.
For $k\in\mathbb N$, let $C_b(E;\mathbb R^k)$ denote the space of bounded continuous functions from $E$ to $\mathbb R^k$, and write $C_b(E):=C_b(E;\mathbb R)$. For $\varphi\in C_b(E)$, set $\|\varphi\|_\infty:=\sup_{x\in E}|\varphi(x)|$. For a measure $\mu$ and a
$\mu$-integrable function $\varphi$, we use the notation
$\mu(\varphi):=\int\varphi\,d\mu$.

Given $\mu,\nu\in\mathcal P_2(\mathbb R^d)$, a coupling between $\mu$ and $\nu$ is a probability measure on
$\mathbb R^d\times\mathbb R^d$ with first marginal $\mu$ and second
marginal $\nu$. We denote the set of such couplings by
$\Pi(\mu,\nu)$. A coupling $\pi\in\Pi(\mu,\nu)$ is a martingale
coupling if, for a disintegration
$\pi(dx,dy)=\mu(dx)\pi_x(dy)$, it satisfies $
\int_{\mathbb R^d}y\,\pi_x(dy)=x$ $\mu\text{-a.e.}$ We denote the set of martingale couplings by $
\Pi_M(\mu,\nu).$ Equivalently, $\pi\in\Pi_M(\mu,\nu)$ if and only if
$\pi\in\Pi(\mu,\nu)$ and
$\int q(x)\cdot(y-x)\,\pi(dx,dy)=0$ for every
$q\in C_b(\mathbb R^d;\mathbb R^d)$. If
$K:=\operatorname{supp}(\mu)$ is compact, it suffices to test this
identity against $q\in C(K;\mathbb R^d)$.

For probability measures $\pi$ and $\rho$ on the same measurable
space, their relative entropy is defined by
\begin{align}\label{eq:relative-entropy}
H(\pi\mid\rho)
:=
\begin{cases}
\displaystyle
\int\log\biggl(\frac{d\pi}{d\rho}\biggr)\,d\pi,
& \pi\ll\rho,\\[1ex]
+\infty,
& \text{otherwise}.
\end{cases}
\end{align}
We also use relative entropy for finite nonnegative measures
$\alpha$ and $\beta$ of the same total mass $m$.
For $m>0$, we set $
H(\alpha\mid\beta)
:=m\,H(\frac{\alpha}{m}\,|\,
            \frac{\beta}{m}),$
and adopt the convention $H(0\mid0):=0$.
For $\eta,\nu\in\mathcal P_2(\mathbb R^d)$, the squared
$2$-Wasserstein distance is defined by
\[
W_2^2(\eta,\nu)
:=
\inf_{\gamma\in\Pi(\eta,\nu)}
\int_{\mathbb R^d\times\mathbb R^d}
|z-y|^2\,\gamma(dz,dy).
\]
Similarly, the $1$-Wasserstein distance is defined by
\[
W_1(\eta,\nu):=\inf_{\gamma\in\Pi(\eta,\nu)}\int_{\mathbb R^d\times\mathbb R^d}|z-y|\,\gamma(dz,dy).
\]
Throughout the paper, the infimum over the empty set is understood to
be $+\infty$.

\subsection{Wasserstein-relaxed entropic martingale optimal transport }\label{subsec: WEMOT}
Let $\mu,\nu\in\mathcal P_2(\mathbb R^d)$ and
$c\in C_b(\mathbb R^d\times\mathbb R^d)$. For
$\lambda,\epsilon>0$, define the Wasserstein-relaxed entropic martingale optimal transport  (WEMOT)  by
\begin{equation}
M^{\lambda,\epsilon}(\mu,\nu)
:=
\inf_{\eta\in\mathcal P_2(\mathbb R^d)}
\biggl\{
\lambda W_2^2(\eta,\nu)
+
\inf_{\pi\in\Pi_M(\mu,\eta)}
\Bigl\{
\int_{\mathbb R^{d}\times \mathbb{R}^d}c(x,z)\,\pi(dx,dz)
+
\epsilon H(\pi\mid\mu\otimes\eta)
\Bigr\}
\biggr\}.
\label{eq:WEMOT}
\end{equation}
The WEMOT problem optimizes jointly over an intermediate marginal
$\eta$ and a martingale coupling $\pi\in\Pi_M(\mu,\eta)$. 
The parameter $\lambda$ controls  the marginal relaxation by penalizing the squared $2$-Wasserstein distance between $\eta$ and the target marginal $\nu$, whereas $\epsilon$ controls the entropy regularization of $\pi$ relative to the product measure $\mu\otimes\eta$. As $\lambda \to \infty$, $\epsilon \to 0$, we expect that $M^{\lambda,\epsilon}(\mu,\nu)$ converges to the martingale optimal transport problem \cite{beiglbock2013model}, defined by
\begin{equation}
M(\mu,\nu)
:=
\inf_{\pi\in\Pi_M(\mu,\nu)}
\int_{\mathbb R^d\times\mathbb R^d}
c(x,y)\,\pi(dx,dy).
\label{eq:mot}
\end{equation}

Let us discuss  two other variations of MOT. 
For a fixed $\epsilon>0$, the entropic martingale optimal transport problem (EMOT) is defined by
\begin{equation}
\operatorname{EMOT}(\mu,\nu):=
\inf_{\pi\in\Pi_M(\mu,\nu)}
\biggl\{
\int_{\mathbb R^d\times\mathbb R^d}c(x,y)\,\pi(dx,dy)
+\epsilon H(\pi\mid\mu\otimes\nu)
\biggr\}.
\label{eq:emot}
\end{equation}
This formulation and related martingale Schr\"odinger bridge problems have been studied from both computational and theoretical perspectives; see \cite{demarch2018entropic,NutzWiesel2024,nutz2023martingale,ChenConfortiRenWang2026,KazeykinaRenWang2026}. 
Another variation preserves the martingale constraint while replacing hard marginal constraints with penalties \cite{doldi2023entropy}.  Here we only consider the relaxation of the second marginal
\begin{equation}
M^\lambda(\mu,\nu):=
\inf_{\substack{\eta\in\mathcal P_2(\mathbb R^d)\\
                \pi\in\Pi_M(\mu,\eta)}}
\biggl\{
\int_{\mathbb R^d\times\mathbb R^d}c(x,y)\,\pi(dx,dy)
+\lambda W_2^2(\eta,\nu)
\biggr\}.
\label{eq:wmot}
\end{equation}

Our WEMOT formulation combines the regularization by relative entropy and the marginal relaxation by $2$-Wasserstein distance. Equivalently, it can be written as
\[
M^{\lambda,\epsilon}(\mu,\nu)
=\inf_{\eta\in\mathcal P_2(\mathbb R^d)}
\bigl\{\operatorname{EMOT}(\mu,\eta)
+\lambda W_2^2(\eta,\nu)\bigr\}.
\]
Marginal relaxation ensures a finite WEMOT value for
arbitrary given $\mu,\nu\in\mathcal P_2(\mathbb R^d)$ and
provides stability under perturbations of the input
marginals. Entropy regularization provides Gibbs representations of the optimal coupling
and a dual formulation suitable for Sinkhorn-type algorithms,
which improve computational efficiency in our numerical
experiments. Let us discuss these properties and compare the formulations
in more detail below.

\medskip

\noindent \emph{Finiteness.} First, WEMOT has a finite value for every
$\mu,\nu\in\mathcal P_2(\mathbb R^d)$.
For the classical MOT problem, Strassen's theorem
\cite{strassen1965existence} states that
$\Pi_M(\mu,\nu)\neq\varnothing$ if and only if
$\mu\preceq_{cx}\nu$.
Here $\mu\preceq_{cx}\nu$ means that
$\int\varphi\,d\mu\leq\int\varphi\,d\nu$ for every
convex $\varphi:\mathbb R^d\to\mathbb R$ for which
both integrals are finite.
WEMOT removes this requirement on the given marginals
by optimizing over an intermediate  marginal $\eta$. To see finiteness of WEMOT, let $X\sim\mu$ and let
$G\sim N(0,I_d)$ be independent of $X$.
Consider $\eta:=\operatorname{Law}(X+G)$ and
$\pi:=\operatorname{Law}(X,X+G)$.
Then $\eta\in\mathcal P_2(\mathbb R^d)$ and
$\pi\in\Pi_M(\mu,\eta)$.
Lemma~\ref{lem:gaussian_information} gives
$H(\pi\mid\mu\otimes\eta)<\infty$ and $W_2(\eta,\nu)<\infty$.
This admissible pair has a finite objective value, and
thus yields finiteness of $M^{\lambda,\epsilon}(\mu,\nu)$. 

This  flexibility of marginal distributions is useful for numerical implementation.
Empirical approximations $\mu^N$ and $\nu^N$ may fail
to satisfy $\mu^N\preceq_{cx}\nu^N$ even when
$\mu\preceq_{cx}\nu$, in which case
$M(\mu^N,\nu^N)=+\infty$.
In particular, convex order requires equal means,
so any difference between the empirical means makes
the classical MOT problem ill-posed.
WEMOT remains finite for these empirical marginals while
preserving the martingale constraint.

\medskip

\noindent \emph{Convexity.}
The objective function in \eqref{eq:WEMOT} is
\begin{equation}\label{eq:objective-function}
J_{\mu, \nu}(\pi,\eta):=
\int_{\mathbb R^d\times\mathbb R^d}c(x,y)\,\pi(dx,dy)
+\epsilon H(\pi\mid\mu\otimes\eta)
+\lambda W_2^2(\eta,\nu).
\end{equation}
The functional $J$ is jointly convex in $(\pi,\eta)$,
since the cost term is linear in $\pi$, relative entropy
is jointly convex, and $\eta\mapsto W_2^2(\eta,\nu)$
is convex. If $\nu$ is absolutely continuous with respect to
Lebesgue measure, $J$ is strictly convex on the set
of admissible pairs with finite relative entropy.

Since the marginal and martingale constraints are linear, the set of admissible pairs is convex. Together with the strict convexity discussed above, this implies the uniqueness of optimizer when $\nu$ is absolutely continuous; existence and uniqueness are established in Theorem~\ref{thm:primal-attainment}. Classical MOT, however, may admit multiple optimizers when $\mu,\nu$ are both absolutely continuous; see \citep{ghoussoub2019structure} for an explicit example.

\medskip
\noindent \emph{Convergence to MOT.}
Entropy regularization alone does not guarantee convergence
to MOT as $\epsilon\downarrow0$.
Indeed, the marginal constraints in EMOT
\eqref{eq:emot} may force every admissible coupling
to have infinite relative entropy, even when the
classical MOT value is finite.

For example, suppose $\mu=\nu$ is nonatomic.
For any $(X,Y)$ with law in $\Pi_M(\mu,\mu)$,
the martingale property gives
$\mathbb E|Y-X|^2=\mathbb E|Y|^2-\mathbb E|X|^2=0$.
Hence $Y=X$ almost surely.
The resulting diagonal coupling is singular with respect
to $\mu\otimes\mu$, so
$\operatorname{EMOT}(\mu,\mu)=+\infty$ for every
$\epsilon>0$, whereas
$M(\mu,\mu)=\int_{\mathbb R^d}c(x,x)\,\mu(dx)<\infty$.
Thus EMOT does not converge to MOT in this example. Marginal relaxation in WEMOT  removes this obstacle by allowing
 martingale couplings with finite relative entropy whose second marginals
approximate the target distribution; see the construction in \emph{Finiteness}.

This observation motivates us to combine marginal relaxation
with entropy regularization.
For $\mu\preceq_{cx}\nu$ and under the regularity
assumptions stated in
Theorem~\ref{thm:convergence-to-mot}, we prove
\[
\lim_{\lambda\to\infty}\lim_{\epsilon\downarrow0}
M^{\lambda,\epsilon}(\mu,\nu)=M(\mu,\nu).
\]
For fixed $\lambda$, letting $\epsilon\downarrow0$
recovers the Wasserstein-relaxed problem $M^\lambda$;
letting $\lambda\to\infty$ then recovers MOT.
Thus the combined formulation permits convergence
to the classical problem even in cases where
fixed-marginal EMOT fails to converge.

\medskip

\noindent \emph{Stability with respect to given marginals.}
Marginal relaxation also addresses instability of the
classical MOT problem.
In dimensions $d\geq2$, empirical approximations
$\mu^N$ and $\nu^N$ may satisfy
$\mu^N\preceq_{cx}\nu^N$ and converge to $\mu$ and $\nu$
in $W_2$, while $M(\mu^N,\nu^N)$ fails to converge
to $M(\mu,\nu)$ \cite{bruckerhoff2022instability}.
The same discrete construction also shows that
fixed-marginal EMOT can be unstable, so entropy
regularization alone does not resolve this problem.

By allowing the second marginal $\eta$ to adjust,
WEMOT remains stable under perturbations of the given
marginals while preserving the martingale constraint.
For fixed $\lambda,\epsilon>0$, we establish qualitative
stability when $c$ is bounded and continuous.
Under additional Lipschitz assumptions on $c$,
Theorem~\ref{thm:quantitative-stability} gives quantitative
bounds on the change in $M^{\lambda,\epsilon}$
in terms of the $W_2$ distances between the given
marginals. These bounds yield a local H\"older rate
of order $2/3$ with respect to the first marginal
and a local Lipschitz rate with respect to the
second marginal.

\medskip
\noindent \emph{Computational benefits.}
The computational role of entropy regularization can
be seen by comparing WEMOT with the Wasserstein-relaxed
MOT problem $M^\lambda$ defined in \eqref{eq:wmot}.
Both formulations allow the second marginal to adjust
while preserving the martingale constraint.
The entropy term in WEMOT additionally yields Gibbs
representations of the coupling density and log-sum-exp
expressions in the dual problem, which will be discussed in Theorem \ref{thm:dual-attainment}.
These expressions support explicit derivatives and
Sinkhorn-type iterative updates; see
\cite{demarch2018entropic} for related computational
methods in EMOT.
The numerical experiments in Section~\ref{sec:numerics}
compare WEMOT with $M^\lambda$ and show faster
convergence of the Sinkhorn type algorithm.

\section{Main Results}\label{sec:main-results}

This section collects the main theoretical results of the paper. Throughout this section, let $\mu,\nu\in\mathcal P_2(\mathbb R^d)$, $c\in C_b(\mathbb R^d\times\mathbb R^d;\mathbb R)$, and regularization parameters $\lambda,\epsilon>0$. Specifically, we present: (i) strong duality for WEMOT; (ii) dual attainment under additional geometric conditions on the marginals $\mu$ and $\nu$; (iii) stability with respect to perturbations of the marginals; and (iv) convergence to the classical MOT problem. Quantitative estimates are provided for the stability and convergence results, while the proofs are deferred to the subsequent sections.

\subsection{Strong duality}

For measurable functions $f:\mathbb R^d\to\mathbb R$ and $h:\mathbb R^d\to\mathbb R^d$, define
\begin{equation}\label{eq:gamma-transform}
\Gamma^\epsilon_{f,h}(z)
:=
-\epsilon\log\biggl(
\int_{\mathbb R^d}
\exp\biggl(
\frac{f(x)+h(x)\cdot(z-x)-c(x,z)}{\epsilon}
\biggr)\mu(dx)
\biggr),
\qquad z\in\mathbb R^d.
\end{equation}
For a measurable function $g:\mathbb R^d\to[-\infty,\infty]$, define its infimal convolution with $\lambda|\cdot|^2$ by
\begin{equation}\label{eq:Q transform}
Q_\lambda g(y)
:=
\inf_{z\in\mathbb R^d}
\bigl\{\lambda|y-z|^2+g(z)\bigr\},
\qquad y\in\mathbb R^d.
\end{equation}
We define the dual value by
\[
D^{\lambda,\epsilon}(\mu,\nu)
:=
\sup_{f\in C_b(\mathbb R^d),\,h\in C_c(\mathbb R^d;\mathbb R^d)}
\bigl\{\mu(f)+\nu(Q_\lambda\Gamma^{\epsilon}_{f,h})\bigr\}.
\]

\begin{theorem}\label{thm:strong-duality}
For any $\mu,\nu \in \mathcal{P}_2(\mathbb R^d)$ and $c \in C_b(\mathbb R^d \times \mathbb R^d;\mathbb R)$, we have  $M^{\lambda,\epsilon}(\mu,\nu)=D^{\lambda,\epsilon}(\mu,\nu)$.
\end{theorem}

The admissible class of dual potentials can be enlarged
beyond $C_b(\mathbb R^d)\times C_c(\mathbb R^d;\mathbb R^d)$
as follows.

\begin{definition}
\label{def:dual-admissible-class}
Let $\mathcal D$ denote the set of measurable pairs
$(f,h)$, with $f:\mathbb R^d\to\mathbb R$ and
$h:\mathbb R^d\to\mathbb R^d$, satisfying
$f\in L^1(\mu)$
and $Q_\lambda\Gamma^{\epsilon}_{f,h}\in L^1(\nu).
$
In addition, for every $\eta\in\mathcal P_2(\mathbb R^d)$
and every $\pi\in\Pi_M(\mu,\eta)$ with
$\pi\ll\mu\otimes\eta$, we require
$\Gamma^{\epsilon}_{f,h}\in L^1(\eta)$ and 
$h(x)\cdot(z-x)\in L^1(\pi)$.
\end{definition}

The following corollary shows that the dual value is
unchanged when the supremum is taken over $\mathcal D$.
This class will also be used in the dual attainment
analysis in Section~\ref{subsec:dual-attainment}.

\begin{corollary}\label{cor:extended-dual-domain}
The dual value is unchanged when the supremum is taken over $\mathcal D$:
\[
D^{\lambda,\epsilon}(\mu,\nu)
=
\sup_{(f,h)\in\mathcal D}
\bigl\{\mu(f)+\nu(Q_\lambda\Gamma^{\epsilon}_{f,h})\bigr\}.
\]
\end{corollary}

\subsection{Dual attainment}\label{subsec:dual-attainment}

We next give conditions under which the dual problem admits an optimizer. Denote the means of $\mu$ and $\nu$ by $m_\mu:=\int_{\mathbb R^d}x\,\mu(dx)$ and $m_\nu:=\int_{\mathbb R^d}y\,\nu(dy)$. Define $R_\mu:=\sup_{x\in\operatorname{supp}(\mu)}|x-m_\mu|$. Define the mean-corrected translate of $\nu$ by $\widehat\nu:=(y\mapsto y+m_\mu-m_\nu)_\#\nu$, so that $\widehat\nu$ has mean $m_\mu$.

For $\xi\in\mathcal P_1(\mathbb R^d)$ with mean $m$, define 
\begin{equation}\label{eq:spread-function}
A(\xi)
:=
\int_{\mathbb R^d}
\frac{(z-m)(z-m)^\top}{|z-m|}\,\xi(dz),\qquad S(\xi)
:=
\theta_{\min}\bigl(A(\xi)\bigr)
\end{equation}
where the integrand is set equal to zero at $z=m$, and $\theta_{\min}(A)$ denotes the smallest eigenvalue of a symmetric matrix $A$. The matrix $A(\xi)$ measures the directional spread of $\xi$ around its mean, while $S(\xi)$ captures the spread in the least dispersed direction. In particular, $A(\xi)$ is positive semidefinite, and $S(\xi)>0$ if and only if $\operatorname{aff}(\operatorname{supp}(\xi))=\mathbb R^d$. Define 
\begin{equation}\label{eq:hdefinition}
\mathfrak h(q):=\frac12[(1+q)\log(1+q)+(1-q)\log(1-q)], \quad \forall \, q \in [0,1).
\end{equation}

When $\mu$ is non-Dirac, we impose the following conditions.

\begin{assumption}\label{ass:dual-attainment}
The following conditions hold:
\begin{itemize}
\item[(i)] $\mu$ is compactly supported and $\operatorname{aff}(\operatorname{supp}(\mu))=\mathbb R^d$.
\item[(ii)] $S(\widehat\nu)>2R_\mu$ and, with $q:=2R_\mu/S(\widehat\nu)\in[0,1)$,
\begin{equation}
\lambda\bigl(S(\widehat\nu)-2R_\mu\bigr)^2
>
9\bigl(2\|c\|_\infty+\epsilon\mathfrak h(q)\bigr).
\end{equation}
\end{itemize}
\end{assumption}

\begin{theorem}\label{thm:dual-attainment}
If $\mu$ is a Dirac measure, the dual problem admits an optimizer in $\mathcal D$. If $\mu$ is non-Dirac and $(\mu,\nu)$ satisfies Assumption~\ref{ass:dual-attainment}, the same conclusion holds. In either case, there exists $(f^*,h^*)\in\mathcal D$ such that
\[
D^{\lambda,\epsilon}(\mu,\nu)
=
\mu(f^*)+\nu\bigl(Q_\lambda\Gamma^\epsilon_{f^*,h^*}\bigr).
\]
\end{theorem}

The proof of Theorem~\ref{thm:dual-attainment} relies on
the following dual-attainment result for fixed-marginal
EMOT. This result is also of independent interest.
Its main feasibility assumption is the existence of a
martingale coupling whose density with respect to the
product measure is bounded above and away from zero.
We first define the corresponding admissible class.

\begin{definition}
\label{def:emot-dual-admissible-class}
For $\mu,\eta\in\mathcal P_2(\mathbb R^d)$, let
$\mathcal D_{\operatorname{EMOT}}(\mu,\eta)$ denote
the set of measurable pairs $(f,h)$, with
$f:\mathbb R^d\to\mathbb R$ and
$h:\mathbb R^d\to\mathbb R^d$, satisfying
$f\in L^1(\mu)$ and $\Gamma^{\epsilon}_{f,h}\in L^1(\eta)$.
In addition, for every $\pi\in\Pi_M(\mu,\eta)$ with
$H(\pi\mid\mu\otimes\eta)<\infty$, we require
$h(x)\cdot(z-x)\in L^1(\pi)$.
\end{definition}

\begin{proposition}
\label{prop:emot-dual-attainment}
Let $\mu,\eta\in\mathcal P_2(\mathbb R^d)$. Suppose that
$\operatorname{supp}(\mu)$ is compact and
$\operatorname{aff}(\operatorname{supp}(\mu))=\mathbb R^d$. Assume further that there exist $\bar\pi\in\Pi_M(\mu,\eta)$ and
$\kappa\in(0,1)$ such that $\bar\pi\ll\mu\otimes\eta$ and
\[
\kappa
\leq
\frac{d\bar\pi}{d(\mu\otimes\eta)}
\leq
\kappa^{-1}
\qquad
(\mu\otimes\eta)\text{-a.e.}
\]
Then strong duality holds and the dual supremum is attained. There exists
$(f_0,h_0)\in\mathcal D_{\operatorname{EMOT}}(\mu,\eta)$ such that
\[
\operatorname{EMOT}(\mu,\eta)
=
\sup_{(f,h)\in\mathcal D_{\operatorname{EMOT}}(\mu,\eta)}
\bigl\{
\mu(f)+\eta(\Gamma^{\epsilon}_{f,h})
\bigr\}
=
\mu(f_0)+\eta(\Gamma^\epsilon_{f_0,h_0}).
\]
\end{proposition}

Proposition~\ref{prop:emot-dual-attainment} complements
existing dual-attainment results for EMOT by allowing the second marginal to have noncompact support.
In particular, \cite{carlier2026weak} establish
multidimensional dual attainment for compactly supported
marginals, obtaining
bounded dual potentials.
Here we also assume compactness of the first marginal but
allow the second marginal to have unbounded support,
under the existence of bounded-density martingale coupling stated
in the proposition.
This is needed for our application to WEMOT, where
the second marginal is optimized and need not be
compactly supported.
Without the compactness of $\operatorname{supp}(\eta)$, their $L^\infty$ estimates
do not apply directly, and we instead establish
attainment in the admissible class
$\mathcal D_{\operatorname{EMOT}}(\mu,\eta)$,
which requires integrability of the marginal potentials.

Our proof uses a Galerkin approximation of the martingale
constraint. We first impose the constraint against
finite-dimensional spaces of test functions and obtain
dual potentials for the finite-dimensional problems.
We then pass to the limit using strong convergence of
the Gibbs densities and convergence of the potentials
after affine normalization.
To the best of our knowledge, this Galerkin approach
has not previously been used to establish dual
attainment for EMOT.

\subsection{Stability under marginal perturbations}

We next analyze the stability of the regularized problem with respect to perturbations of the marginals.

\begin{assumption}\label{ass:c-lipschitz}
The cost function $c$ is  Lipschitz continuous. There exists $L_c>0$ such that
\[|c(x,y)-c(x',y')|\leq L_c\bigl(|x-x'|+|y-y'|\bigr)\]
for all $(x,y),(x',y')\in\mathbb R^d\times\mathbb R^d$.
\end{assumption}

To state the quantitative estimate, define
\begin{align}\label{eq:K estimate of stability}
K(\mu,\nu)
:=
\Big(
\frac{2\|c\|_\infty}{\lambda}
+
\frac{\epsilon}{2\lambda}m_2(\mu)
+
\big(\sqrt d+W_2(\mu,\nu)\big)^2
\Big)^{1/2},
\end{align}
where $m_2(\mu):=\int_{\mathbb R^d}|x|^2\,\mu(dx)$.

\begin{theorem}\label{thm:quantitative-stability}
Suppose Assumption~\ref{ass:c-lipschitz} holds. Let
$\mu,\widetilde\mu,\nu,\widetilde\nu\in\mathcal P_2(\mathbb R^d)$. Set
$\delta:=W_2(\mu,\widetilde\mu)$ and
$\rho:=W_2(\nu,\widetilde\nu)$. Define
$\overline K:=\max\{K(\mu,\nu),K(\widetilde\mu,\widetilde\nu)\}$.
Then there exists a constant \mbox{$C=C(d,L_c,\lambda,\epsilon,\overline K)>0$} such that
\[
\bigl|
M^{\lambda,\epsilon}(\widetilde\mu,\widetilde\nu)
-
M^{\lambda,\epsilon}(\mu,\nu)
\bigr|
\leq
C\bigl\{s(\delta,\rho)+s(\delta,\rho)^2\bigr\},
\]
where $s(\delta,\rho):=\delta^{2/3}+\delta+\rho$. In particular, for fixed
$(\mu,\nu)$, as $\delta,\rho\to0$,
\[
\bigl|
M^{\lambda,\epsilon}(\widetilde\mu,\widetilde\nu)
-
M^{\lambda,\epsilon}(\mu,\nu)
\bigr|
=
O\bigl(\delta^{2/3}+\rho\bigr).
\]
\end{theorem} 

\begin{remark}
\begin{enumerate}
    \item Theorem~\ref{thm:quantitative-stability} gives a local H\"older estimate of order $2/3$ in the first marginal and a local Lipschitz estimate in the second marginal. The different rates reflect the asymmetric treatment of the two marginals in our formulation. The first marginal is not relaxed in WEMOT, so perturbing $\mu$ therefore requires constructing a new coupling that preserves the martingale constraint, whereas perturbing $\nu$ only changes the Wasserstein penalty, which yields a linear rate.
    \item The linear dependence on the second-marginal perturbation is sharp. Indeed, in dimension one, take $c\equiv0$, $\mu=\delta_0$, $\nu=\delta_1$, and $\widetilde\nu=\delta_{1+t}$ with $t>0$. For any $a\in\mathbb R$, the martingale constraint requires the intermediate marginal to have mean zero, so $M^{\lambda,\epsilon}(\delta_0,\delta_a)=\lambda a^2$, attained by $\eta=\delta_0$. Consequently,
\[
M^{\lambda,\epsilon}(\mu,\widetilde\nu)-M^{\lambda,\epsilon}(\mu,\nu)=2\lambda t+\lambda t^2,\qquad W_2(\nu,\widetilde\nu)=t,
\]
showing that the linear order cannot generally be improved.
\end{enumerate}
\end{remark}

The significance of this estimate is that stability holds in every dimension and does not require the given marginals to remain in convex order. This addresses the instability of multidimensional MOT exhibited by \citet{bruckerhoff2022instability}, whose discrete construction also yields instability for fixed-marginal EMOT.  By relaxing the second marginal, our formulation provides an explicit bound on the effect of marginal approximation errors while preserving the martingale constraint. This is particularly useful for empirical and discretized inputs. The numerical experiments in Sections~\ref{subsec:stability-comparison} and~\ref{subsec:bruckerhoff-instability} illustrate this stability through lower variance across empirical samples and  convergence along the Br\"uckerhoff--Juillet perturbation sequence, respectively.

\subsection{Convergence to martingale optimal transport}

We finally study the approximation to the classical martingale optimal transport problem as $\lambda\to\infty$ and $\epsilon\to0$. The following theorem first establishes convergence and then gives a quantitative estimate that allows the parameters to vary simultaneously. The bound separates the error from marginal relaxation, of order $1/\lambda$, from the contribution of entropy regularization, controlled by $\epsilon\log(1+\lambda)$. It therefore provides a sufficient condition for choosing the parameters jointly. For the quantitative estimate, we additionally impose the following compactness assumption on $\nu$.

\begin{assumption}\label{ass:nu-compact}
The measure $\nu$ is compactly supported.
\end{assumption}

\begin{theorem}\label{thm:convergence-to-mot}
Suppose Assumption~\ref{ass:c-lipschitz} holds. Let $\mu, \nu \in \mathcal{P}_2(\mathbb{R}^d)$ and $\mu\preceq_{cx}\nu$.
\begin{itemize}
\item[(i)] $
\lim_{\lambda\to\infty}\lim_{\epsilon\to0}
M^{\lambda,\epsilon}(\mu,\nu)
=
M(\mu,\nu).$

\item[(ii)] Suppose in addition that Assumption \ref{ass:dual-attainment} and \ref{ass:nu-compact} hold. Then there exists a positive constant
$C=C(d,\mu,\nu,L_c,\|c\|_\infty)$, independent of $\lambda$ and $\epsilon$, such that
\[
\bigl|
M^{\lambda,\epsilon}(\mu,\nu)-M(\mu,\nu)
\bigr|
\leq
C\biggl\{
\frac{1}{\lambda}
+
\epsilon\log(1+\lambda)
\biggr\}.
\]
Consequently, $M^{\lambda,\epsilon}(\mu,\nu)\to M(\mu,\nu)$ for any $\mu,\nu$ satisfying the above assumptions such that $\lambda\to\infty$ and $\epsilon\log(1+\lambda)\to0$. 
\end{itemize}
\end{theorem}

The proof uses Gaussian smoothing to control the entropy term and a bounded MOT martingale dual potential to derive the $1/\lambda$ estimate. Assumptions~\ref{ass:dual-attainment} and~\ref{ass:nu-compact} allow us to apply the uniform bounds on dual potentials established by \citet{carlier2026weak} in arbitrary dimension. Combining these estimates yields the stated convergence rate.


\section{Strong Duality}
\label{sec:proof-strong-duality}

In this section, we prove Theorem~\ref{thm:strong-duality} by applying the Fenchel--Moreau theorem. We first introduce the topological vector spaces used in the argument.

For $p\in\{0,1,2\}$, let $\mathcal M_p(\mathbb R^n)$ denote the space of finite signed Borel measures $\kappa$ on $\mathbb R^n$ satisfying
$\int_{\mathbb R^n}(1+|u|^p)\,|\kappa|(du)<\infty$, and let
$C_p(\mathbb R^n)$ be the space of continuous functions $\varphi:\mathbb R^n\to\mathbb R$ such that
$\sup_{u\in\mathbb R^n}|\varphi(u)|/(1+|u|^p)<\infty$.
We write $\mathcal M(\mathbb R^n):=\mathcal M_0(\mathbb R^n)$ and
$C_0(\mathbb R^n):=C_b(\mathbb R^n)$. For $0\leq p\leq k\leq2$, let
$\sigma(\mathcal M_k(\mathbb R^n),C_p(\mathbb R^n))$ denote the coarsest topology on $\mathcal M_k(\mathbb R^n)$ for which the map
$\kappa\mapsto\int_{\mathbb R^n}\varphi(u)\,\kappa(du)$
is continuous for every $\varphi\in C_p(\mathbb R^n)$.

Set
\[
\mathcal Y
:=
\mathcal M(\mathbb R^d)
\times
\mathcal M_1(\mathbb R^d)
\times
\mathcal M_2(\mathbb R^d)
\times
\mathcal M(\mathbb R^d;\mathbb R^d),
\]
where $\mathcal M(\mathbb R^d;\mathbb R^d)$ denotes the space of finite $\mathbb R^d$-valued signed Borel measures. We equip $\mathcal Y$ with the product topology
\[
\tau
:=
\sigma\bigl(\mathcal M(\mathbb R^d),C_b(\mathbb R^d)\bigr)
\times
\sigma\bigl(\mathcal M_1(\mathbb R^d),C_1(\mathbb R^d)\bigr)
\times
\sigma\bigl(\mathcal M_2(\mathbb R^d),C_2(\mathbb R^d)\bigr)
\times
\sigma\bigl(
\mathcal M(\mathbb R^d;\mathbb R^d),
C_c(\mathbb R^d;\mathbb R^d)
\bigr).
\]
Each of the test-function spaces appearing above is a separating vector space of linear functionals on the corresponding measure space. Hence, by \cite[Theorem~3.10]{rudin1991functional}, $(\mathcal Y,\tau)$ is a Hausdorff locally convex topological vector space, and its continuous dual can be identified with
\[
\mathcal Y^*
=
C_b(\mathbb R^d)
\times
C_1(\mathbb R^d)
\times
C_2(\mathbb R^d)
\times
C_c(\mathbb R^d;\mathbb R^d).
\]
We also set $
\mathcal X
:=
\mathcal M_1(\mathbb R^{2d})
\times
\mathcal M_2(\mathbb R^{2d}),$
and equip it with the product topology
\[
\omega
:=
\sigma\bigl(\mathcal M_1(\mathbb R^{2d}),C_b(\mathbb R^{2d})\bigr)
\times
\sigma\bigl(\mathcal M_2(\mathbb R^{2d}),C_b(\mathbb R^{2d})\bigr).
\]
Thus, convergence with respect to $\omega$ is the ordinary weak convergence against bounded continuous functions.

For $\pi\in\mathcal M_1(\mathbb R^{2d})$, define the finite
$\mathbb R^d$-valued signed measure $m(\pi)$ by
\[
\langle h,m(\pi)\rangle
:=
\int_{\mathbb R^{2d}}
h(x)\cdot(z-x)\,\pi(dx,dz),
\qquad
h\in C_c(\mathbb R^d;\mathbb R^d).
\]
If $\pi\geq0$, $\pi_1=\mu$, and
$\pi(dx,dz)=\mu(dx)\pi_x(dz)$, then $
m(\pi)(dx)
=(
\int_{\mathbb R^d}(z-x)\,\pi_x(dz))\mu(dx).$
Consequently, $m(\pi)=\mathbf 0$ if and only if $\pi$ satisfies the martingale constraint.
Define the linear map $A:\mathcal X\to\mathcal Y$ by 
\[
A(\pi,\gamma)
:=
\bigl(
\pi_1,\,
\pi_2-\gamma_1,\,
\gamma_2,\,
m(\pi)
\bigr).\]

We next define the extended-real-valued functional $\mathcal F^{\lambda,\epsilon}:\mathcal X\to(-\infty,\infty]$ by
\[
\mathcal F^{\lambda,\epsilon}(\pi,\gamma):=\int_{\mathbb R^d\times\mathbb R^d}c(x,z)\,\pi(dx,dz)+\epsilon H(\pi\mid\mu\otimes\pi_2)+\lambda\int_{\mathbb R^d\times\mathbb R^d}|z-y|^2\,\gamma(dz,dy)
\]
when $\pi,\gamma\geq0$, and set $\mathcal F^{\lambda,\epsilon}(\pi,\gamma):=+\infty$ otherwise.
For $a\in\mathcal Y$, define the perturbation value function $\Phi: \mathcal{Y}\to \mathbb{R}$
\[
\Phi(a)
:=
\inf\bigl\{
\mathcal F^{\lambda,\epsilon}(\pi,\gamma):
(\pi,\gamma)\in\mathcal X,\ 
A(\pi,\gamma)=a
\bigr\},
\]
with the convention $\inf\varnothing=+\infty$. Finally, set $a^*:=(\mu,0,\nu,\mathbf 0)\in\mathcal Y$. By construction, whenever $\mathcal F^{\lambda,\epsilon}(\pi,\gamma)<\infty$, the constraint $A(\pi,\gamma)=a^*$ is equivalent to the existence of some $\eta\in\mathcal P_2(\mathbb R^d)$ such that $\pi\in\Pi_M(\mu,\eta)$ and $\gamma\in\Pi(\eta,\nu)$. It follows directly $\Phi(a^*)=M^{\lambda,\epsilon}(\mu,\nu)$. 

It remains to show that $\Phi$ is proper, convex, and lower semicontinuous on $(\mathcal Y,\tau)$. The Fenchel--Moreau theorem then gives $\Phi=\Phi^{**}$, and the desired dual formulation follows by computing the convex conjugate of $\Phi$.

\begin{proposition}\label{prop:Phi-properties}
The function $\Phi:\mathcal Y\to(-\infty,\infty]$ is proper, convex, and lower semicontinuous with respect to $\tau$.
\end{proposition}

\begin{proof}
\noindent\emph{Properness.}
Since $\Phi(a^*)=M^{\lambda,\epsilon}(\mu,\nu)<\infty$, the function $\Phi$ is not identically equal to $+\infty$. To see that $\Phi$ never
takes the value $-\infty$, fix $a=(a_1,a_2,a_3,a_4)\in\mathcal Y$. If
$\{(\pi,\gamma):A(\pi,\gamma)=a\}$ is empty, then
$\Phi(a)=+\infty$. Otherwise, every feasible pair in this fiber satisfies
$\pi_1=a_1$, and hence
$\pi(\mathbb R^{2d})=a_1(\mathbb R^d)$. Since $c$ is bounded and the
entropy and quadratic transport terms are nonnegative,
\[
\mathcal F^{\lambda,\epsilon}(\pi,\gamma)
\geq
-\|c\|_\infty a_1(\mathbb R^d).
\]
It follows that $\Phi(a)>-\infty$. Hence $\Phi$ is proper.

\noindent\emph{Convexity.}
The functional $\mathcal F^{\lambda,\epsilon}$ is convex. Indeed, the cost and quadratic transport terms are linear in $\pi$ and $\gamma$, respectively, while the map $
\pi\longmapsto H(\pi\mid\mu\otimes\pi_2)$
is convex by the joint convexity of relative entropy and the linearity of the map
$\pi\mapsto(\pi,\mu\otimes\pi_2)$; see, for example,
\cite[Appendix~D.3]{dembozeitouni2010large}. Since $A$ is linear, taking the infimum of $\mathcal F^{\lambda,\epsilon}$ over the affine fibers of $A$ preserves convexity. Thus $\Phi$ is convex.

\noindent\emph{Lower semicontinuity.}
Let $(a^\alpha)_\alpha$ be a net in $\mathcal Y$ converging to
$a=(a_1,a_2,a_3,a_4)$ with respect to $\tau$, where
$a^\alpha=(a_1^\alpha,a_2^\alpha,a_3^\alpha,a_4^\alpha)$. We claim that $L:=\liminf_\alpha\Phi(a^\alpha)\geq \phi (a)$. The assertion is immediate if
$L=+\infty$, so we suppose that $L<\infty$. Passing to a subnet, which we do not relabel, we may assume that $\Phi(a^\alpha)\to L$.

Fix $\delta>0$. For each $\alpha$, choose nonnegative measures
$(\pi^\alpha,\gamma^\alpha)$ such that
$A(\pi^\alpha,\gamma^\alpha)=a^\alpha$ and
\[
\mathcal F^{\lambda,\epsilon}(\pi^\alpha,\gamma^\alpha)
\leq
\Phi(a^\alpha)+\delta.
\]
After restricting to a tail of the net, we may assume that
$B:=\sup_\alpha\mathcal F^{\lambda,\epsilon}(\pi^\alpha,\gamma^\alpha)<\infty$. By the definition of $\mathcal F^{\lambda,\epsilon}$, this implies that
$\pi^\alpha,\gamma^\alpha\geq0$ for every $\alpha$.

 Since
$\pi_1^\alpha=a_1^\alpha$ and
$a_1^\alpha\to a_1$ in
$\sigma(\mathcal M(\mathbb R^d),C_b(\mathbb R^d))$, testing against the constant function $1$ gives
$\pi^\alpha(\mathbb R^{2d})\to a_1(\mathbb R^d)$. In particular,
$M_0:=\sup_\alpha \pi^\alpha(\mathbb R^{2d})<\infty$. Since
$\int c\,d\pi^\alpha\geq-\|c\|_\infty \pi^\alpha(\mathbb R^{2d})$, we obtain
\[
\sup_\alpha
H(\pi^\alpha\mid\mu\otimes\pi_2^\alpha)
\leq
\frac{B+\|c\|_\infty M_0}{\epsilon}
=:C_H
\]
and
\[
\sup_\alpha
\int_{\mathbb R^{2d}}
|z-y|^2\,\gamma^\alpha(dz,dy)
\leq
\frac{B+\|c\|_\infty M_0}{\lambda}
=:C_Q.
\]
We first prove uniform tightness of $(\gamma^\alpha)_\alpha$. Since
$\gamma_2^\alpha=a_3^\alpha$ and
$a_3^\alpha\to a_3$ in
$\sigma(\mathcal M_2(\mathbb R^d),C_2(\mathbb R^d))$, testing against
$1$ and $y\mapsto|y|^2$ yields
\[
\sup_\alpha\gamma^\alpha(\mathbb R^{2d})=\sup_\alpha\gamma_2^\alpha(\mathbb R^{d})<\infty,
\qquad
\sup_\alpha
\int_{\mathbb R^d}|y|^2\,\gamma_2^\alpha(dy)<\infty.
\]
Using $|z|^2\leq2|z-y|^2+2|y|^2$, we further obtain
\begin{align}
\sup_\alpha
\int_{\mathbb R^d}|z|^2\,\gamma_1^\alpha(dz)
\leq
2C_Q
+
2\sup_\alpha
\int_{\mathbb R^d}|y|^2\,\gamma_2^\alpha(dy)
<\infty.
\label{eq:gamma-first-second-moment}
\end{align}
Thus both marginals of $\gamma^\alpha$ are uniformly tight, and hence so is the family $(\gamma^\alpha)_\alpha$.

We next consider the second marginals of $\pi^\alpha$. Since
$\pi_2^\alpha-\gamma_1^\alpha=a_2^\alpha$, we have
\[
\int_{\mathbb R^d}|z|\,\pi_2^\alpha(dz)
\leq
\int_{\mathbb R^d}|z|\,\gamma_1^\alpha(dz)
+
\bigl|
\langle |z|,a_2^\alpha\rangle
\bigr|.
\]
The first term is uniformly bounded by
\eqref{eq:gamma-first-second-moment} and the Cauchy--Schwarz inequality. The second term is uniformly bounded because
$a_2^\alpha\to a_2$ in
$\sigma(\mathcal M_1(\mathbb R^d),C_1(\mathbb R^d))$ and
$z\mapsto|z|$ belongs to $C_1(\mathbb R^d)$. Therefore $
\sup_\alpha
\int_{\mathbb R^d}|z|\,\pi_2^\alpha(dz)<\infty,$
and $(\pi_2^\alpha)_\alpha$ is uniformly tight.

It remains to establish uniform tightness of $(\pi_1^\alpha)_\alpha$. This does not follow merely from the weak convergence of the net
$(a_1^\alpha)_\alpha$. Instead, we use the entropy bound. The first marginal of $\mu\otimes\pi_2^\alpha$ is $\pi^\alpha(\mathbb{R}^{2d})\mu$, and hence the data-processing inequality for relative entropy gives
\[
H(\pi_1^\alpha\mid \pi^\alpha(\mathbb{R}^{2d})\mu)
\leq
H(\pi^\alpha\mid\mu\otimes\pi_2^\alpha)
\leq C_H.
\]
By the variational representation of relative entropy, for every Borel set
$E\subset\mathbb R^d$ and every $t>0$,
\[
t\pi_1^\alpha(E)
\leq
C_H
+
\pi^\alpha(\mathbb{R}^{2d})\mu(E)(e^t-1).
\]
Let $K\subset\mathbb R^d$ be compact and set $r:=\mu(K^c)$. If $r\in(0,1)$, choosing $t=\log(1/r)$ yields $
\sup_\alpha\pi_1^\alpha(K^c)
\leq
\frac{C_H+M_0}{\log(1/r)}.$
Since $\mu$ is tight, $K$ can be chosen so that $r$ is arbitrarily small, and the right-hand side then converges to zero. If $r=0$, the same conclusion follows from
$\pi_1^\alpha\ll \pi^{\alpha}(\mathbb{R}^{2d})\mu$. Hence
$(\pi_1^\alpha)_\alpha$ is uniformly tight. Together with the uniform tightness of $(\pi_2^\alpha)_\alpha$, this implies that
$(\pi^\alpha)_\alpha$ is uniformly tight.

The families $(\pi^\alpha)_\alpha$ and $(\gamma^\alpha)_\alpha$ are therefore uniformly tight, consist of nonnegative measures, and have uniformly bounded total masses. By Prokhorov's theorem, their weak closures are compact. Passing to a  subnet, which we again do not relabel, there exist nonnegative finite Borel measures $\pi$ and $\gamma$ such that
\begin{align}
(\pi^\alpha,\gamma^\alpha)
\to
(\pi,\gamma)
\qquad\text{in }\omega.
\label{eq:omega-convergence}
\end{align}

We now verify that $A(\pi,\gamma)=a$. From
\eqref{eq:omega-convergence},
$\pi_1^\alpha\Rightarrow\pi_1$ and
$\gamma_2^\alpha\Rightarrow\gamma_2$. On the other hand,
$\pi_1^\alpha=a_1^\alpha\to a_1$ and
$\gamma_2^\alpha=a_3^\alpha\to a_3$ weakly. The uniqueness of weak limits gives $
\pi_1=a_1
$ and $
\gamma_2=a_3.$
Similarly, $\gamma_1^\alpha\Rightarrow\gamma_1$. Moreover,
\eqref{eq:gamma-first-second-moment} implies that the first moments of
$(\gamma_1^\alpha)_\alpha$ are uniformly integrable, since
\[
\sup_\alpha
\int_{\{|z|>R\}}
|z|\,\gamma_1^\alpha(dz)
\leq
\frac{1}{R}
\sup_\alpha
\int_{\mathbb R^d}|z|^2\,\gamma_1^\alpha(dz)
\to0
\]
as $R\to\infty$. Therefore, the weak convergence
$\gamma_1^\alpha\Rightarrow\gamma_1$ together with the uniform
integrability of the first moments yields
\[
\gamma_1^\alpha
\to
\gamma_1
\qquad
\text{in }
\sigma(\mathcal M_1(\mathbb R^d),C_1(\mathbb R^d)).
\]
Since
$\pi_2^\alpha=\gamma_1^\alpha+a_2^\alpha$ and
$a_2^\alpha\to a_2$ in $\sigma(\mathcal M_1(\mathbb R^d),C_1(\mathbb R^d))$, it follows that 
\[
\pi_2^\alpha
\to
\gamma_1+a_2 \qquad
\text{in }
\sigma(\mathcal M_1(\mathbb R^d),C_1(\mathbb R^d)).\]
At the same time, \eqref{eq:omega-convergence} implies
$\pi_2^\alpha\Rightarrow\pi_2$. Hence
$\pi_2=\gamma_1+a_2$, and consequently 
\begin{align}
\pi_2-\gamma_1=a_2,\qquad  \pi_2^\alpha \to
\pi_2
\qquad
\text{in }
\sigma(\mathcal M_1(\mathbb R^d),C_1(\mathbb R^d)).
\label{eq:second-constraint-limit}
\end{align}

It remains to identify the martingale component. Fix
$h\in C_c(\mathbb R^d;\mathbb R^d)$ and set
$G_h(x,z):=h(x)\cdot(z-x)$. Since $h$ has compact support,
$G_h$ is continuous and satisfies
$|G_h(x,z)|\leq C_h(1+|z|)$ for some $C_h>0$. Let  $B_R$ be a ball of radius $R>0$ centered at $0$. Let $\chi_R\in C_c(\mathbb R^d)$ satisfy
$0\leq\chi_R\leq1$, $\chi_R=1$ on $B_R$, and
$\chi_R=0$ outside $B_{2R}$.
By \eqref{eq:second-constraint-limit},
\[
\lim_\alpha
\int (1+|z|)(1-\chi_R(z))\,\pi_2^\alpha(dz)
=
\int (1+|z|)(1-\chi_R(z))\,\pi_2(dz).
\]
Moreover, $(x,z)\mapsto G_h(x,z)\chi_R(z)$ is bounded
and continuous, so weak convergence gives
\[
\int G_h(x,z)\chi_R(z)\,\pi^\alpha(dx,dz)
\to
\int G_h(x,z)\chi_R(z)\,\pi(dx,dz).
\]
Using $|G_h(x,z)|\leq C_h(1+|z|)$, we obtain $
\limsup_\alpha
\left|\int G_h\,d\pi^\alpha-\int G_h\,d\pi\right|
\leq
2C_h\int (1+|z|)(1-\chi_R(z))\,\pi_2(dz).$
Since $\pi_2\in\mathcal M_1(\mathbb R^d)$, the right-hand
side tends to zero as $R\to\infty$ by dominated convergence.
Consequently, 
\begin{align}
\int_{\mathbb R^{2d}}G_h(x,z)\,\pi^\alpha(dx,dz)
\to
\int_{\mathbb R^{2d}}G_h(x,z)\,\pi(dx,dz). \label{eq:martingale-integral-limit}
\end{align}

Since $m(\pi^\alpha)=a_4^\alpha$ and
$a_4^\alpha\to a_4$ against
$C_c(\mathbb R^d;\mathbb R^d)$, the left-hand side of
\eqref{eq:martingale-integral-limit} also converges to
$\langle h,a_4\rangle$. It follows that
\begin{equation}\label{eq:a4-equality}
\int_{\mathbb R^{2d}}
h(x)\cdot(z-x)\,\pi(dx,dz)
=
\int h(x)\, a_4(dx)
\qquad
\text{for every }
h\in C_c(\mathbb R^d;\mathbb R^d).
\end{equation}
Define the finite vector measure $
\zeta(E):=\int_{E\times\mathbb R^d}z\,\pi(dx,dz).$
\eqref{eq:a4-equality} gives $x\pi_1=\zeta-a_4$, and hence 
\[
\int_{\mathbb R^d}|x|\,\pi_1(dx)
\leq
\int_{\mathbb R^d}|z|\,\pi_2(dz)
+|a_4|(\mathbb R^d)
<\infty.
\]
Thus $\pi\in\mathcal M_1(\mathbb R^{2d})$.
Since both marginals of $\gamma^\alpha$ have uniformly
bounded second moments, the Portmanteau theorem gives
\[
\int_{\mathbb R^{2d}}(|z|^2+|y|^2)\,\gamma(dz,dy)
\leq
\liminf_\alpha
\int_{\mathbb R^{2d}}(|z|^2+|y|^2)\,
\gamma^\alpha(dz,dy)
<\infty.
\]
Hence $\gamma\in\mathcal M_2(\mathbb R^{2d})$.
Consequently, $(\pi,\gamma)\in\mathcal X$ and $m(\pi)=a_4$.
Together with \eqref{eq:a4-equality}, this proves
$A(\pi,\gamma)=a$.

Finally, \eqref{eq:omega-convergence} and the bounded continuity of $c$ yield
\[
\int_{\mathbb R^{2d}}c(x,z)\,\pi^\alpha(dx,dz)
\to
\int_{\mathbb R^{2d}}c(x,z)\,\pi(dx,dz).
\]
Moreover, $\pi_2^\alpha\Rightarrow\pi_2$, and hence
$\mu\otimes\pi_2^\alpha\Rightarrow\mu\otimes\pi_2$. The joint lower semicontinuity of relative entropy
\cite[Appendix~D.3]{dembozeitouni2010large} and the Portmanteau theorem for the nonnegative lower semicontinuous function
$(z,y)\mapsto|z-y|^2$ imply $
\mathcal F^{\lambda,\epsilon}(\pi,\gamma)
\leq
\liminf_\alpha
\mathcal F^{\lambda,\epsilon}(\pi^\alpha,\gamma^\alpha).$
Since $A(\pi,\gamma)=a$, we conclude that
\[
\Phi(a)
\leq
\mathcal F^{\lambda,\epsilon}(\pi,\gamma)
\leq
\liminf_\alpha
\mathcal F^{\lambda,\epsilon}(\pi^\alpha,\gamma^\alpha)
\leq
L+\delta.
\]
Letting $\delta\downarrow0$ gives
$\Phi(a)\leq L=\liminf_\alpha\Phi(a^\alpha)$. Hence $\Phi$ is
$\tau$-lower semicontinuous.
\end{proof}

To prove Theorem~\ref{thm:strong-duality}, we use the following version of the Fenchel--Moreau theorem; see, for example, \cite[Theorem~2.22]{barbuprecupanu2012convexity}.

\begin{theorem}[Fenchel--Moreau]\label{thm:fenchel-moreau}
Let $E$ be a Hausdorff locally convex topological vector space with continuous dual $E^*$. For a proper function $f:E\to(-\infty,+\infty]$, define its convex conjugate by
\[
f^*(\ell)
:=
\sup_{x\in E}
\bigl\{
\langle \ell,x\rangle-f(x)
\bigr\},
\qquad
\ell\in E^*,
\]
and its biconjugate by
\[
f^{**}(x)
:=
\sup_{\ell\in E^*}
\bigl\{
\langle \ell,x\rangle-f^*(\ell)
\bigr\},
\qquad
x\in E.
\]
If $f$ is convex and lower semicontinuous, then $f=f^{**}$ on $E$.
\end{theorem}

\begin{proof}[Proof of Theorem~\ref{thm:strong-duality}]
By Proposition~\ref{prop:Phi-properties}, $\Phi$ is proper, convex, and lower semicontinuous on $(\mathcal Y,\tau)$. Hence, the Fenchel--Moreau theorem yields $\Phi=\Phi^{**}$ on $\mathcal Y$. Since
$\Phi(a^*)=M^{\lambda,\epsilon}(\mu,\nu)$, we have
\begin{align*}
M^{\lambda,\epsilon}(\mu,\nu)
&=
\Phi^{**}(a^*)=
\sup_{\ell\in\mathcal Y^*}
\bigl\{
\langle\ell,a^*\rangle-\Phi^*(\ell)
\bigr\}=
\sup_{\ell\in\mathcal Y^*}
\inf_{a\in\mathcal Y}
\bigl\{
\Phi(a)+\langle\ell,a^*-a\rangle
\bigr\}\\
&=
\sup_{\ell\in\mathcal Y^*}
\inf_{\substack{(\pi,\gamma)\in\mathcal X\\ \pi,\gamma\geq0}}
\bigl\{
\mathcal F^{\lambda,\epsilon}(\pi,\gamma)
+
\langle\ell,a^*-A(\pi,\gamma)\rangle
\bigr\}.
\end{align*}
The last equality follows from the definition of $\Phi$ and the fact that $\mathcal F^{\lambda,\epsilon}$ is equal to $+\infty$ outside the cone of nonnegative measures.

Write $
\ell=(f,q,u,h)
\in
C_b(\mathbb R^d)
\times
C_1(\mathbb R^d)
\times
C_2(\mathbb R^d)
\times
C_c(\mathbb R^d;\mathbb R^d).$
Using $a^*=(\mu,0,\nu,\mathbf 0)$ and the definition of $A$, set the Lagrangian 
\begin{align*}
&\mathcal L(\pi,\gamma;f,q,u,h)
:=\mathcal F^{\lambda,\epsilon}(\pi,\gamma)+\langle\ell,a^*-A(\pi,\gamma)\rangle\\
=&\mu(f)+\nu(u)+\epsilon H(\pi\mid\mu\otimes\pi_2)
+\int_{\mathbb R^{2d}}\bigl[c(x,z)-f(x)-q(z)-h(x)\cdot(z-x)\bigr]\,\pi(dx,dz)\\
&\quad+\int_{\mathbb R^{2d}}\bigl[\lambda|z-y|^2+q(z)-u(y)\bigr]\,\gamma(dz,dy).
\end{align*}
It follows that
\[
M^{\lambda,\epsilon}(\mu,\nu)
=
\sup_{(f,q,u,h)\in\mathcal Y^*}
\inf_{\substack{(\pi,\gamma)\in\mathcal X\\ \pi,\gamma\geq0}}
\mathcal L(\pi,\gamma;f,q,u,h).
\]

Fix
$(f,q,u,h)\in
\mathcal{Y}^*$.
We first minimize the $\pi$-dependent part of the Lagrangian. For $\pi\geq0$, write
$\eta:=\pi_2$ and disintegrate
$\pi(dx,dz)=\pi_z(dx)\eta(dz)$. The chain rule for relative entropy \cite[Theorem~D.13]{dembozeitouni2010large} gives $
H(\pi\mid\mu\otimes\eta)
=
\int_{\mathbb R^d}
H(\pi_z\mid\mu)\,\eta(dz).$
For each $z\in\mathbb R^d$, the Gibbs variational formula \cite[Appendix~D]{dembozeitouni2010large} yields
\[
\inf_{\rho\in\mathcal{P}(\mathbb{R}^d)}
\biggl\{
\int_{\mathbb R^d}
\bigl[
c(x,z)-f(x)-h(x)\cdot(z-x)
\bigr]\,\rho(dx)
+
\epsilon H(\rho\mid\mu)
\biggr\}
=
\Gamma^{\epsilon}_{f,h}(z).
\]
Therefore the $\pi$-dependent part of the Lagrangian is bounded from below by $
\int_{\mathbb R^d}
\bigl[
\Gamma^{\epsilon}_{f,h}(z)-q(z)
\bigr]\,\eta(dz).$
It follows that
\[
\inf_{\substack{\pi\in\mathcal M_1(\mathbb R^{2d})\\ \pi\geq0}}
\biggl\{
\int_{\mathbb R^{2d}}\bigl[c(x,z)-f(x)-q(z)-h(x)\cdot(z-x)\bigr]\,\pi(dx,dz)
+\epsilon H(\pi\mid\mu\otimes\pi_2)
\biggr\}
=0
\]
if $q\leq\Gamma^{\epsilon}_{f,h}$, and the infimum is $-\infty$ otherwise.
Indeed, when $q\leq\Gamma^{\epsilon}_{f,h}$, the value zero is attained at $\pi=0$. If $q(z_0)>\Gamma^{\epsilon}_{f,h}(z_0)$ for some $z_0$, one may concentrate $\pi_2$ at $z_0$ and scale its total mass to obtain the value $-\infty$.

Similarly,
\[
\inf_{\substack{\gamma\in\mathcal M_2(\mathbb R^{2d})\\ \gamma\geq0}}
\int_{\mathbb R^{2d}}\bigl[\lambda|z-y|^2+q(z)-u(y)\bigr]\,\gamma(dz,dy)
=0
\]
if $u(y)\leq\lambda|z-y|^2+q(z)$ for all $y,z\in\mathbb R^d$, and the infimum is $-\infty$ otherwise.
Therefore,
\begin{equation}\label{eq:intermediate-dual}
M^{\lambda,\epsilon}(\mu,\nu)
=
\sup\bigl\{\mu(f)+\nu(u):(f,q,u,h)\in\mathcal Y^*,\ q\leq\Gamma^{\epsilon}_{f,h},\ u\leq Q_\lambda q\bigr\},
\end{equation}
where $Q_\lambda q(y):=\inf_{z\in\mathbb R^d}\{\lambda|z-y|^2+q(z)\}$.

We now eliminate $u$ and $q$. For $q\in C_1(\mathbb R^d)$, the function $
Q_\lambda q(y)
:=
\inf_{z\in\mathbb R^d}\{
\lambda|y-z|^2+q(z)\}$
is finite and continuous. Moreover, if
$|q(z)|\leq K(1+|z|)$, then $
-K(1+|y|)-\frac{K^2}{4\lambda}
\leq
Q_\lambda q(y)
\leq
K(1+|y|),$
and hence $Q_\lambda q\in C_1(\mathbb R^d)\subset C_2(\mathbb R^d)$. The constraint on $u$ in \eqref{eq:intermediate-dual} is equivalent to
$u\leq Q_\lambda q$. Since $\nu$ is nonnegative, the largest admissible choice is therefore
$u=Q_\lambda q$.

For $f\in C_b(\mathbb R^d)$ and
$h\in C_c(\mathbb R^d;\mathbb R^d)$, the function
$\Gamma^{\epsilon}_{f,h}$ is finite and continuous. Moreover, since $f$ and $c$ are bounded and both $h$ and $x\mapsto h(x)\cdot x$ are bounded, there exists a constant $C>0$ such that
$|\Gamma^{\epsilon}_{f,h}(z)|\leq C(1+|z|)$. Thus
$\Gamma^{\epsilon}_{f,h}\in C_1(\mathbb R^d)$.

Finally, $Q_\lambda$ is monotone: if $q_1\leq q_2$, then
$Q_\lambda q_1\leq Q_\lambda q_2$. Hence, under the constraint
$q\leq\Gamma^{\epsilon}_{f,h}$, the largest dual value is obtained by choosing
$q=\Gamma^{\epsilon}_{f,h}$. We conclude that
\[
M^{\lambda,\epsilon}(\mu,\nu)
=
\sup_{\substack{
f\in C_b(\mathbb R^d)\\
h\in C_c(\mathbb R^d;\mathbb R^d)
}}
\bigl\{
\mu(f)+\nu(Q_\lambda\Gamma^{\epsilon}_{f,h})
\bigr\}
=
D^{\lambda,\epsilon}(\mu,\nu).
\]
\end{proof}

We next extend the dual formulation to the larger admissible class introduced in Corollary~\ref{cor:extended-dual-domain}.
\begin{proof}[Proof of Corollary~\ref{cor:extended-dual-domain}]
Firstly, we show that
$C_b(\mathbb R^d)\times C_c(\mathbb R^d;\mathbb R^d)\subset\mathcal D$.
Let $f\in C_b(\mathbb R^d)$ and
$h\in C_c(\mathbb R^d;\mathbb R^d)$.  As shown in the proof of Theorem~\ref{thm:strong-duality},
$\Gamma^{\epsilon}_{f,h}\in C_1(\mathbb R^d)$ and
$Q_\lambda\Gamma^{\epsilon}_{f,h}\in C_1(\mathbb R^d)$, and hence
$Q_\lambda\Gamma^{\epsilon}_{f,h}\in L^1(\nu)$. 

Now let $\eta\in\mathcal P_2(\mathbb R^d)$ and
$\pi\in\Pi_M(\mu,\eta)$ with $\pi\ll\mu\otimes\eta$. The linear growth
of $\Gamma^{\epsilon}_{f,h}$ implies $\Gamma^{\epsilon}_{f,h}\in L^1(\eta)$. Since
$h$ has compact support, there exists $C_h>0$ such that
$|h(x)\cdot(z-x)|\leq C_h(1+|z|)$, and therefore
$h(x)\cdot(z-x)\in L^1(\pi)$. Thus $(f,h)\in\mathcal D$. Theorem~\ref{thm:strong-duality} therefore gives
\[
M^{\lambda,\epsilon}(\mu,\nu)
\leq
\sup_{(f,h)\in\mathcal D}
\bigl\{
\mu(f)+\nu(Q_\lambda\Gamma^{\epsilon}_{f,h})
\bigr\}.
\]

It remains to prove the reverse inequality. Fix $(f,h)\in\mathcal D$
and let $\eta\in\mathcal P_2(\mathbb R^d)$,
$\pi\in\Pi_M(\mu,\eta)$, and
$\gamma\in\Pi(\eta,\nu)$. If
$\mathcal F^{\lambda,\epsilon}(\pi,\gamma)=+\infty$, there is nothing
to prove, so assume that it is finite. In particular,
$\pi\ll\mu\otimes\eta$. Disintegrate
$\pi(dx,dz)=\pi_z(dx)\eta(dz)$. By the Gibbs variational inequality,
for $\eta$-a.e. $z$,
\[
\epsilon H(\pi_z\mid\mu)
\geq
\int_{\mathbb R^d}
\bigl[
f(x)+h(x)\cdot(z-x)-c(x,z)+\Gamma^{\epsilon}_{f,h}(z)
\bigr]
\,\pi_z(dx).
\]
Integrating with respect to $\eta$ and using the chain rule for relative
entropy yields
\[
\int_{\mathbb R^{2d}}c(x,z)\,\pi(dx,dz)
+
\epsilon H(\pi\mid\mu\otimes\eta)
\geq
\mu(f)
+
\int_{\mathbb R^{2d}}
h(x)\cdot(z-x)\,\pi(dx,dz)
+
\eta(\Gamma^{\epsilon}_{f,h}).
\]
All terms are well defined by the definition of $\mathcal D$. Since
$\pi$ is a martingale coupling, the middle term vanishes, and hence
\[
\int_{\mathbb R^{2d}}c(x,z)\,\pi(dx,dz)
+
\epsilon H(\pi\mid\mu\otimes\eta)
\geq
\mu(f)+\eta(\Gamma^{\epsilon}_{f,h}).
\]
Adding the quadratic transport cost gives
\begin{align*}
\mathcal F^{\lambda,\epsilon}(\pi,\gamma)
&\geq
\mu(f)
+
\int_{\mathbb R^{2d}}
\bigl[
\lambda|z-y|^2+\Gamma^{\epsilon}_{f,h}(z)
\bigr]
\,\gamma(dz,dy)\geq
\mu(f)
+
\int_{\mathbb R^d}
Q_\lambda\Gamma^{\epsilon}_{f,h}(y)\,\nu(dy).
\end{align*}
Taking the infimum over all admissible
$\eta$, $\pi\in\Pi_M(\mu,\eta)$, and $\gamma\in\Pi(\eta,\nu)$ yields
\[
M^{\lambda,\epsilon}(\mu,\nu)
\geq
\mu(f)+\nu(Q_\lambda\Gamma^{\epsilon}_{f,h}).
\]
Since $(f,h)\in\mathcal D$ was arbitrary, we obtain the reversed inequality
\[
M^{\lambda,\epsilon}(\mu,\nu)
\geq
\sup_{(f,h)\in\mathcal D}
\bigl\{
\mu(f)+\nu(Q_\lambda\Gamma^{\epsilon}_{f,h})
\bigr\}.
\]
Combining the two inequalities proves
\[
D^{\lambda,\epsilon}(\mu,\nu)
=
\sup_{(f,h)\in\mathcal D}
\bigl\{
\mu(f)+\nu(Q_\lambda\Gamma^{\epsilon}_{f,h})
\bigr\}.
\]
\end{proof}

\section{Primal and Dual Attainment of WEMOT}\label{sec:primal and dual attainment}

We next study attainment in the primal and dual WEMOT problems. We first establish the existence of a primal optimizer and its uniqueness when $\nu$ is absolutely continuous. We then prove dual attainment under additional assumptions.

\subsection{Existence and uniqueness of primal optimizers}\label{subsec:primal-welldefinedness}

We first establish primal attainment by a compactness argument. When $\nu$ is absolutely continuous, uniqueness follows from the strict convexity of the Wasserstein penalty and the entropy term.

\begin{theorem}\label{thm:primal-attainment}
Let $\mu,\nu\in\mathcal P_2(\mathbb R^d)$, $c\in C_b(\mathbb R^d\times\mathbb R^d)$, and $\lambda,\epsilon>0$. Then there exists an optimal pair $(\pi^*,\eta^*)$, with $\eta^*\in\mathcal P_2(\mathbb R^d)$ and $\pi^*\in\Pi_M(\mu,\eta^*)$, such that
\[
M^{\lambda,\epsilon}(\mu,\nu)=\int c\,d\pi^*+\epsilon H(\pi^*\mid\mu\otimes\eta^*)+\lambda W_2^2(\eta^*,\nu).
\]
Moreover, if $\nu$ is absolutely continuous with respect to Lebesgue measure, then the optimal pair is unique.
\end{theorem}

\begin{proof}
For an admissible pair $(\pi,\eta)$, write
\[
J_{\mu,\nu}(\pi,\eta):=\int c\,d\pi+\epsilon H(\pi\mid\mu\otimes\eta)+\lambda W_2^2(\eta,\nu).
\]
Let $X\sim\mu$ and let $G\sim\mathcal N(0,I_d)$ be independent of $X$. The pair $\pi^G:=\mathcal L(X,X+G)$ and $\eta^G:=\mathcal L(X+G)$ is admissible. By Lemma~\ref{lem:gaussian_information} and the bound $W_2(\eta^G,\nu)\leq\sqrt d+W_2(\mu,\nu)$, we obtain
\begin{equation}\label{eq:primal-value-bound}
M^{\lambda,\epsilon}(\mu,\nu)\leq J_{\mu,\nu}(\pi^G,\eta^G)\leq\|c\|_\infty+\frac{\epsilon}{2}m_2(\mu)+\lambda\bigl(\sqrt d+W_2(\mu,\nu)\bigr)^2<\infty.
\end{equation} Since $J_{\mu,\nu}\geq-\|c\|_\infty$, the value is also bounded below. Choose a minimizing sequence $(\pi_n,\eta_n)$ satisfying $\pi_n\in\Pi_M(\mu,\eta_n)$ and $J_{\mu,\nu}(\pi_n,\eta_n)\leq M^{\lambda,\epsilon}(\mu,\nu)+\frac1n.$
The nonnegativity of relative entropy and the boundedness of $c$ imply $
\lambda W_2^2(\eta_n,\nu)\leq M^{\lambda,\epsilon}(\mu,\nu)+\|c\|_\infty+1.$
Using the triangular inequality $m_2(\eta_n)^{1/2}\leq m_2(\nu)^{1/2}+W_2(\eta_n,\nu)$, we obtain $\sup_n m_2(\eta_n)<\infty$. Consequently,
\[
\sup_n\int_{\mathbb R^{2d}}\bigl(|x|^2+|y|^2\bigr)\,\pi_n(dx,dy)=m_2(\mu)+\sup_n m_2(\eta_n)<\infty.
\]
Thus $(\pi_n)_n$ is tight. By Prokhorov's theorem, after passing to a subsequence, $\pi_n$ converges weakly to a probability measure $\pi^*$. Its first marginal is $\mu$, and its second marginal, denoted by $\eta^*$, is the weak limit of $(\eta_n)_n$. Lower semicontinuity of the second moment gives $\eta^*\in\mathcal P_2(\mathbb R^d)$.

We next verify the martingale constraint. For every $q\in C_b(\mathbb R^d;\mathbb R^d)$, the uniform second-moment bound implies uniform integrability of $q(x)\cdot(y-x)$ with respect to $(\pi_n)_n$. Since this function is continuous, weak convergence yields
\[
\int_{\mathbb R^{2d}}q(x)\cdot(y-x)\,\pi^*(dx,dy)=\lim_{n\to\infty}\int_{\mathbb R^{2d}}q(x)\cdot(y-x)\,\pi_n(dx,dy)=0.
\]
Hence $\pi^*\in\Pi_M(\mu,\eta^*)$.

Moreover, $\mu\otimes\eta_n$ converges weakly to $\mu\otimes\eta^*$. Joint lower semicontinuity of relative entropy, weak lower semicontinuity of the quadratic transport cost, and bounded continuity of $c$ therefore imply
\[
J_{\mu,\nu}(\pi^*,\eta^*)\leq\liminf_{n\to\infty}J_{\mu,\nu}(\pi_n,\eta_n)=M^{\lambda,\epsilon}(\mu,\nu).
\]
Since $(\pi^*,\eta^*)$ is admissible, it is an optimal pair.

Suppose now that $\nu$ is absolutely continuous. We first recall why $\eta\mapsto W_2^2(\eta,\nu)$ is strictly convex under linear interpolation. Let $\eta_0\neq\eta_1$, let $\gamma_i\in\Pi(\nu,\eta_i)$ be optimal transport plans, and fix $t\in(0,1)$. Set $\eta_t:=(1-t)\eta_0+t\eta_1$. The coupling $(1-t)\gamma_0+t\gamma_1$ gives
\[
W_2^2(\nu,\eta_t)\leq(1-t)W_2^2(\nu,\eta_0)+tW_2^2(\nu,\eta_1).
\]
If equality held, this mixture would be optimal. By Brenier's theorem \cite[Theorem~2.12]{villani2003topics}, each optimal plan with first marginal $\nu$ is induced by a transport map. Writing $\gamma_i=(\mathrm{id},T_i)_\#\nu$, the conditional law $(1-t)\delta_{T_0(y)}+t\delta_{T_1(y)}$ could therefore be a Dirac measure only if $T_0=T_1$ $\nu$-almost everywhere. This would imply $\eta_0=\eta_1$, a contradiction.

Finally, the admissible set is convex, the cost term is linear, and $(\pi,\eta)\mapsto H(\pi\mid\mu\otimes\eta)$ is jointly convex. Two optimal pairs with different second marginals would therefore contradict the strict convexity of the Wasserstein penalty. If their second marginals coincide, strict convexity of relative entropy in its first argument forces their couplings to coincide as well. Thus the optimal pair is unique.
\end{proof}

\subsection{Dual attainment}\label{subsec:proof-dual-attainment}
In this section, we prove Theorem~\ref{thm:dual-attainment}. We first consider the case where $\mu$ is non-Dirac and Assumption~\ref{ass:dual-attainment} holds. Recall that $\eta\in\mathcal P_2(\mathbb R^d)$ denotes the second marginal of the martingale coupling, which is allowed to differ from $\nu$ through the Wasserstein penalty. By Theorem~\ref{thm:primal-attainment}, choose an optimal pair $(\pi^*,\eta^*)$ for the primal WEMOT problem. Then $\pi^*$ also minimizes the fixed-marginal problem $\operatorname{EMOT}(\mu,\eta^*)$, and
\[
M^{\lambda,\epsilon}(\mu,\nu)=\operatorname{EMOT}(\mu,\eta^*)+\lambda W_2^2(\eta^*,\nu).
\]
The proof proceeds in three steps. First, we show that the pair $(\mu,\eta^*)$ satisfies the hypotheses of Proposition~\ref{prop:emot-dual-attainment}, which yields a dual optimizer $(f_0,h_0)$ for the fixed-marginal problem $\operatorname{EMOT}(\mu,\eta^*)$. Second, using the strong duality result in Theorem~\ref{thm:strong-duality} and Corollary~\ref{cor:extended-dual-domain}, we choose a maximizing sequence for the WEMOT dual problem. Decomposing the WEMOT duality gap shows that this sequence approaches the optimal dual value of the fixed-marginal problem $\operatorname{EMOT}(\mu,\eta^*)$. At the same time, The same sequence also induces an asymptotically optimal sequence for the dual problem associated with $\lambda W_2^2(\eta^*,\nu)$. Finally, we normalize the potentials by affine gauge
transformations, establish compactness of the normalized sequence, and
pass to the limit to obtain an WEMOT dual optimizer. The case where $\mu$ is a Dirac measure is treated separately at the end of the section.

We begin with the first step. The compactness and full-dimensionality of $\operatorname{supp}(\mu)$ are already guaranteed by Assumption~\ref{ass:dual-attainment}. It remains to construct a martingale coupling between $\mu$ and $\eta^*$ whose density relative to $\mu\otimes\eta^*$ is bounded above and away from zero. To this end, we first establish a stability estimate for the map $\xi\mapsto S(\xi)$ defined in \eqref{eq:spread-function}. This estimate allows us to transfer the nondegeneracy condition on $\widehat\nu$ in Assumption~\ref{ass:dual-attainment}(ii) to the optimal relaxed marginal $\eta^*$.

\begin{lemma}
\label{lem:score-spread-stability}
Let $\xi_1,\xi_2\in\mathcal P_1(\mathbb R^d)$ have the same mean $m$. Then
\[
\bigl|
S(\xi_1)-S(\xi_2)
\bigr|
\leq
3W_1(\xi_1,\xi_2).
\]
\end{lemma}

\begin{proof}
For each $v\in\mathbb R^d$ with $|v|=1$, define
\[
F_v(z)
:=
\begin{cases}
\dfrac{(v\cdot(z-m))^2}{|z-m|}, & z\neq m,\\[2mm]
0, & z=m.
\end{cases}
\]
The function $F_v$ is continuous at $m$. For $z\neq m$,
\[
\nabla F_v(z)
=
2\frac{v\cdot(z-m)}{|z-m|}v
-
\frac{(v\cdot(z-m))^2}{|z-m|^3}(z-m).
\]
Since $|v|=1$, we have $|\nabla F_v(z)|\leq3$. Hence  $F_v$ is globally $3$-Lipschitz. By Kantorovich--Rubinstein duality \cite{villani2009optimal}, $W_1(\xi_1,\xi_2)=\sup_{\operatorname{Lip}(\varphi)\leq1}|\int\varphi\,d\xi_1-\int\varphi\,d\xi_2|$, where the supremum is taken over all $1$-Lipschitz functions $\varphi:\mathbb R^d\to\mathbb R$. Applying this representation to $F_v/3$ yields
\[
3W_1(\xi_1,\xi_2)\geq\Big|\int F_v\,d\xi_1-\int F_v\,d\xi_2\Big|=\Bigl|v^\top\bigl(A(\xi_1)-A(\xi_2)\bigr)v\Bigr|,
\]
where $A(\cdot)$ is defined in \eqref{eq:spread-function}. 
Since $A(\xi_1)-A(\xi_2)$ is symmetric, taking the supremum over $|v|=1$ yields $
\|A(\xi_1)-A(\xi_2)\|_{\mathrm{op}}
\leq
3W_1(\xi_1,\xi_2).$
It implies that $
\bigl|S(\xi_1)-S(\xi_2)\bigr|\leq3W_1(\xi_1,\xi_2),$
which proves the claim.
\end{proof}

We now combine Lemma~\ref{lem:score-spread-stability} with the optimality condition of the relaxed marginal $\eta^*$ to verify the hypotheses of Proposition~\ref{prop:emot-dual-attainment} for $\operatorname{EMOT}(\mu, \eta^*)$.

\begin{proposition}
\label{prop:input-condition}
Suppose Assumption~\ref{ass:dual-attainment} holds. Then
$\operatorname{aff}(\operatorname{supp}(\eta^*))=\mathbb R^d$.
Moreover, there exist $\bar\pi^*\in\Pi_M(\mu,\eta^*)$ and
$\kappa^*\in(0,1)$ such that
$\bar\pi^*\ll\mu\otimes\eta^*$ and, writing
$\bar r^*:=d\bar\pi^*/d(\mu\otimes\eta^*)$,
\[
\kappa^*
\leq
\bar r^*
\leq
(\kappa^*)^{-1}
\qquad
(\mu\otimes\eta^*)\text{-a.e.}
\]
Consequently, the pair $(\mu,\eta^*)$ satisfies the hypotheses of
Proposition~\ref{prop:emot-dual-attainment}.
\end{proposition}

\begin{proof}

Since the primal optimizer of WEMOT $\pi^*\in\Pi_M(\mu,\eta^*)$, we have $\mu\preceq_{cx}\eta^*$. Lemma~\ref{lem:support-convex-order} and $\operatorname{aff}(\operatorname{supp}\mu)=\mathbb R^d$ therefore imply $\operatorname{aff}(\operatorname{supp}\eta^*)=\mathbb R^d$.

For a probability measure $\xi\in \mathcal{P}_2 (\mathbb{R}^d)$ with mean $m$, define 
\[
s_\xi(z)
:=
\begin{cases}
\dfrac{z-m}{|z-m|}, & z\neq m,\\[2mm]
0, & z=m,
\end{cases}
\qquad
\psi_\xi(z):=s_\xi(z)-a_\xi, \qquad \forall z\in \mathbb{R}^d
\]
where $a_\xi:=\int s_\xi\,d\xi$.
Then $\int\psi_\xi\,d\xi=0$, $|\psi_\xi|\leq2$, and  $
\int (z-m)\psi_\xi(z)^\top\,\xi(dz)=A(\xi).$

We first construct a  martingale coupling between $\mu$ and $\widehat\nu$, where $\widehat\nu:=(y\mapsto y+m_\mu-m_\nu)_\#\nu$. And $\hat{\nu}$ has mean $m_{\mu}$.  By Assumption \ref{ass:dual-attainment} (ii),
$S(\widehat\nu)>2R_\mu$, and hence
$A(\widehat\nu)$ is invertible. Define $
\bar r(x,z)
:=
1+
\psi_{\widehat\nu}(z)^\top
A(\widehat\nu)^{-1}(x-m_\mu).$
Since $\int\psi_{\widehat\nu}\,d\widehat\nu=0$ and
$\int(x-m_\mu)\,\mu(dx)=0$, the measure
$\bar\pi:=\bar r(\mu\otimes\widehat\nu)$ has marginals
$\mu$ and $\widehat\nu$. Moreover,
\[
\int(z-m_\mu)\bar r(x,z)\,\widehat\nu(dz)
=
A(\widehat\nu)
A(\widehat\nu)^{-1}(x-m_\mu)
=
x-m_\mu,
\]
so $\bar\pi\in\Pi_M(\mu,\widehat\nu)$.

Since $|\psi_{\widehat\nu}|\leq2$, $
|\bar r(x,z)-1|
\leq
\frac{2R_\mu}{S(\widehat\nu)}
=:
q<1$ by Assumption \ref{ass:dual-attainment}.
Hence $1-q\leq\bar r\leq1+q$. Moreover,  convexity of
$u\mapsto u\log u$ gives $
 H(\bar\pi\,|\,\mu\otimes\widehat\nu)
\leq \mathfrak{h}(q)$, where $\mathfrak{h}$ is defined in \eqref{eq:hdefinition}.
Using $(\widehat\nu,\bar\pi)$ as an admissible pair for WEMOT and observing that
$W_2^2(\widehat\nu,\nu)=|m_\mu-m_\nu|^2$, we obtain
\[
M^{\lambda,\epsilon}(\mu,\nu)
\leq
\|c\|_{\infty}+\epsilon \mathfrak{h}(q)+\lambda|m_\mu-m_\nu|^2.
\]

On the other hand, since $\mu \preceq_{cx} \eta^*$,  $\eta^*$ has
mean $m_\mu$. Since $\widehat\nu:=(y\mapsto y+m_\mu-m_\nu)_{\#}\nu$, one has
\[
W_2^2(\eta^*,\nu)
=
|m_\mu-m_\nu|^2
+
W_2^2(\eta^*,\widehat\nu).
\]
Since relative
entropy is nonnegative,
\[
M^{\lambda,\epsilon}(\mu,\nu)
\geq
-\|c\|_{\infty}
+
\lambda|m_\mu-m_\nu|^2
+
\lambda W_2^2(\eta^*,\widehat\nu).
\]
Combining the preceding two estimates yields $
W_2(\eta^*,\widehat\nu)
\leq
\sqrt{
\frac{
2\|c\|_{\infty}+\epsilon \mathfrak{h}(q)
}{\lambda}
}.$

By Lemma~\ref{lem:score-spread-stability} and $W_1(\eta^*, \widehat{\nu})\leq W_2(\eta^*, \widehat{\nu})$, $
S(\eta^*)
\geq
S(\widehat\nu)
-
3
\sqrt{
\frac{
2\|c\|_{\infty}+\epsilon \mathfrak{h}(q)
}{\lambda}
}.$
Assumption \ref{ass:dual-attainment}~(ii) gives $
3
\sqrt{
\frac{
2\|c\|_{\infty}+\epsilon \mathfrak{h}(q)
}{\lambda}
}
<
S(\widehat\nu)-2R_\mu.$
Therefore $S(\eta^*)>2R_\mu$. We then apply the same construction to $\eta^*$. Define $
\bar r^*(x,z)
:=
1+
\psi_{\eta^*}(z)^\top
A(\eta^*)^{-1}(x-m_\mu).$
As above, $\bar r^*(\mu\otimes\eta^*)$ defines a martingale coupling
$\bar\pi^*\in\Pi_M(\mu,\eta^*)$. Writing
$q^*:=2R_\mu/S(\eta^*)<1$, we have
$1-q^*\leq\bar r^*\leq1+q^*$. Therefore, with
$\kappa^*:=1-q^*\in(0,1)$, $
\kappa^*
\leq
\bar r^*
\leq
1+q^*
\leq
\frac{1}{1-q^*}
=
(\kappa^*)^{-1},$
which proves the required density bound.
\end{proof}

Consequently, $(\mu,\eta^*)$ satisfies the assumptions of Proposition~\ref{prop:emot-dual-attainment}.  Therefore the  dual problem of $\operatorname{EMOT}(\mu, \eta^*)$ admits an optimizer in $\mathcal D_{\operatorname{EMOT}}(\mu,\eta^*)$. Since $(\eta^*, \pi^*)$ is an optimizer for $M^{\lambda, \epsilon}(\mu, \nu)$, $\pi^*$ is also the primal optimizer for $\operatorname{EMOT}(\mu,\eta^*)$. The following lemma shows that,  for any $(f,h)$ in $ \mathcal{D}_{\operatorname{EMOT}}(\mu, \eta^*)$, the primal--dual gap $\operatorname{EMOT}(\mu, \eta^*)-(\mu(f)+\eta^*(\Gamma_{f,h}^\epsilon))$  equals $\epsilon H(\pi^*\mid r_{f,h}(\mu\otimes \eta^*))$, where the density $r_{f,h}$ is defined in the lemma.

\begin{lemma}
\label{lem:fixed-marginal-entropy-gap}
Let $(f,h)\in\mathcal D_{\operatorname{EMOT}}(\mu,\eta^*)$ and set
$g:=\Gamma^{\epsilon}_{f,h}$. Define
\[
r_{f,h}(x,z)
:=
\exp\biggl(
\frac{
f(x)+h(x)\cdot(z-x)-c(x,z)+g(z)
}{\epsilon}
\biggr).
\]
Then $r_{f,h}(\mu\otimes\eta^*)$ is a probability measure with second marginal $\eta^*$, and
\[
\operatorname{EMOT}(\mu,\eta^*)
-
\mu(f)
-
\eta^*(g)
=
\epsilon
H\bigl(
\pi^*
\mid
r_{f,h}(\mu\otimes\eta^*)
\bigr).
\]
\end{lemma}

\begin{proof}
By the definition of $g=\Gamma^{\epsilon}_{f,h}$, for $\eta^*$-a.e.\
$z\in\mathbb R^d$,
\[
e^{-g(z)/\epsilon}
=
\int_{\mathbb R^d}
\exp\biggl(
\frac{
f(x)+h(x)\cdot(z-x)-c(x,z)
}{\epsilon}
\biggr)\mu(dx).
\]
Consequently, $
\int_{\mathbb R^d}r_{f,h}(x,z)\,\mu(dx)
=1$, $
\text{for }\eta^*\text{-a.e. }z.$
It follows that $r_{f,h}(\mu\otimes\eta^*)$ is a probability measure
with second marginal $\eta^*$.

Since $\pi^*$ is the primal optimizer and
$\operatorname{EMOT}(\mu,\eta^*)<\infty$, the boundedness of $c$
implies that
$H(\pi^*\mid\mu\otimes\eta^*)<\infty$. In particular,
$\pi^*\ll\mu\otimes\eta^*$. Moreover,
$r_{f,h}>0$ $(\mu\otimes\eta^*)$-a.e., and hence
$\pi^*\ll r_{f,h}(\mu\otimes\eta^*)$.  The change-of-reference identity
for relative entropy gives
\begin{align*}
&\epsilon
H\bigl(
\pi^*
\mid
r_{f,h}(\mu\otimes\eta^*)
\bigr)
=
\epsilon H(\pi^*\mid\mu\otimes\eta^*)
-
\epsilon
\int_{\mathbb R^{2d}}
\log r_{f,h}(x,z)\,\pi^*(dx,dz)\\
=&
\epsilon H(\pi^*\mid\mu\otimes\eta^*)
+
\int_{\mathbb R^{2d}}c(x,z)\,\pi^*(dx,dz)
-\mu(f)
-\eta^*(g)
-
\int_{\mathbb R^{2d}}
h(x)\cdot(z-x)\,\pi^*(dx,dz).
\end{align*}
All the terms above are well defined by the definition of
$\mathcal D_{\operatorname{EMOT}}(\mu,\eta^*)$. Since
$\pi^*\in\Pi_M(\mu,\eta^*)$, the last integral vanishes. Finally,
the optimality of $\pi^*$ yields
\[
\epsilon
H\bigl(
\pi^*
\mid
r_{f,h}(\mu\otimes\eta^*)
\bigr)
=
\operatorname{EMOT}(\mu,\eta^*)
-
\mu(f)
-
\eta^*(g),
\]
which proves the claim.
\end{proof}

Fix a dual optimizer
$(f_0,h_0)\in
\mathcal D_{\operatorname{EMOT}}(\mu,\eta^*)$
and set 
\begin{equation*}
g_0:=\Gamma^\epsilon_{f_0,h_0}, \qquad r_0:=r_{f_0,h_0}=\exp\biggl(
\frac{
f_0(x)+h_0(x)\cdot(z-x)-c(x,z)+g_0(z)
}{\epsilon}
\biggr)
\end{equation*}
Dual attainment yields that $
\operatorname{EMOT}(\mu,\eta^*)
=
\mu(f_0)+\eta^*(g_0)$.
Applying Lemma~\ref{lem:fixed-marginal-entropy-gap} to
$(f_0,h_0)$ therefore implies that $
H\bigl(
\pi^*
\mid
r_0(\mu\otimes\eta^*)
\bigr)
=
0.$
Since relative entropy between probability measures vanishes only when
the two measures coincide, it follows that $
\pi^*
=
r_0(\mu\otimes\eta^*).$
Moreover, the exponential definition of $r_0$ implies that
$r_0>0$, $(\mu\otimes\eta^*)$-a.e.

We next prove a key ingredient for the dual attainment of WEMOT.  Let $(f_n,h_n)\in\mathcal D$ be a maximizing sequence for the WEMOT dual problem and set $g_n:=\Gamma^\epsilon_{f_n,h_n}$. The following lemma shows that, if the second-marginal potentials $g_n$ converge in $\eta^*$-probability after affine gauge transformations and $\int[\lambda|z-y|^2+g_n(z)-Q_\lambda g_n(y)]\,\gamma(dy,dz)\to0$ for a coupling $\gamma\in\Pi(\nu,\eta^*)$, then the slopes of the affine transformations are uniformly bounded. 

\begin{lemma}
\label{lem:affine-slope}
Let $\gamma\in\Pi(\nu,\eta^*)$. Suppose there exist $a_n\in\mathbb R^d$ such that
\[
\xi_n(z)
:=
g_n(z)+a_n\cdot z-g_0(z)
\to
0
\qquad
\text{in }\eta^*\text{-probability}.
\]
Define $
s_n(y,z)
:=
\lambda|y-z|^2
+
g_n(z)
-
Q_\lambda g_n(y).$ If $
\int_{\mathbb R^d\times\mathbb R^d}
s_n(y,z)\,\gamma(dy,dz)
\to
0,$
then the sequence $(a_n)_n$ is bounded, that is, $\sup_n|a_n|<\infty$.
\end{lemma}

\begin{proof}
Since $s_n\geq0$, $s_n(Y,Z)\to0$ in $\gamma$-probability, where $(Y,Z)\sim\gamma$. Let $Z'$ be independent of $(Y,Z)$ with law $\eta^*$. Then $\xi_n(Z)\to0$ and $\xi_n(Z')\to0$ in $(\gamma\otimes\eta^*)$-probability.

For every $y,z,z'$, $
Q_\lambda g_n(y)
\leq
\lambda|y-z'|^2+ g_n(z'),$
and therefore
\[
0
\leq
a_n\cdot(Z-Z')
+
R(Y,Z,Z')
+
\chi_n,
\]
$(\gamma\otimes\eta^*)$-almost surely, where $
R(Y,Z,Z')
:=
\lambda\bigl(
|Y-Z'|^2-|Y-Z|^2
\bigr)
+
g_0(Z')-g_0(Z)$
and $
\chi_n
:=
s_n(Y,Z)+\xi_n(Z')-\xi_n(Z)\to0$
in $(\gamma\otimes\eta^*)$-probability.

Suppose that $|a_n|\to\infty$ along a subsequence and, after passing
to a further subsequence, write $
e_n:=\frac{a_n}{|a_n|}\to e,$ $
|e|=1.$
By Proposition \ref{prop:input-condition}, $\eta^*$ has full affine dimension. Thus $e\cdot Z$ is not almost
surely constant, otherwise there exist a constant $c$ such that $\eta^*$ is concentrated on the affine hyperplane $\{z\in \mathbb{R}^d: e\cdot z =c\}$. Hence there exist bounded Borel sets
$B_-,B_+$ of positive $\eta^*$-measure and $\delta>0$ such that
\[
e\cdot(z-z')\leq-2\delta
\qquad
\text{for }z\in B_-,\ z'\in B_+.
\]
For all sufficiently large $n$, $
e_n\cdot(z-z')\leq-\delta$ on $B_-\times B_+.$
Since $R(Y, Z, Z')$ is finite almost surely, we can choose $M$ such that the event $\{Z\in B_-,\,Z'\in B_+,\,|R(Y,Z,Z')|\leq M\}$ has positive probability.
Since $\chi_n\to0$ in probability, for all sufficiently large $n$,
\[
(\gamma\otimes\eta^*)\bigl(\{Z\in B_-,\,Z'\in B_+,\,|R(Y,Z,Z')|\leq M,\,|\chi_n|\leq1\}\bigr)>0.
\]
On this event,
\[
a_n\cdot(Z-Z')+R+\chi_n
\leq
-|a_n|\delta+M+1<0
\]
for large $n$, and it leads to a contradiction. Hence $(a_n)_n$ is bounded.
\end{proof}

We next show that the $Q_\lambda$ transform defined in \eqref{eq:Q transform} is continuous along some sequence $(g_n)_{n}$.  

\begin{lemma}
\label{lem:wasserstein-complementarity}
Let
$\gamma\in\Pi(\nu,\eta^*)$ be an optimal coupling for
$W_2^2(\eta^*,\nu)$. $f_n:\mathbb R^d\to\mathbb R$ and $h_n:\mathbb R^d\to\mathbb R^d$ be measurable functions, and set $g_n:=\Gamma^\epsilon_{f_n,h_n}$. Suppose that there exist measurable functions $f, h$ and $g:=\Gamma^{\epsilon}_{f,h}\in L^1(\eta^*)$ such that 
\[
f_n\to f,\qquad h_n\to h
\quad \mu\text{-a.e.},
\qquad
g_n\to g
\quad \eta^*\text{-a.e.}
\]
Define $
s_n(y,z)
:=
\lambda|y-z|^2+g_n(z)-Q_\lambda g_n(y).$
If $
\int_{\mathbb R^{2d}}s_n(y,z)\,\gamma(dy,dz)\to0,$
then
\[
Q_\lambda g(y)
=
\lambda|y-z|^2+g(z)
\qquad
\gamma\text{-a.e.}
\]
Consequently, $
Q_\lambda g\in L^1(\nu)$
and $
\nu(Q_\lambda g)
=
\lambda W_2^2(\eta^*,\nu)+\eta^*(g).$
\end{lemma}

\begin{proof}
Since $s_n\geq0$, 
$s_n\to0$ in $\gamma$-probability. After passing to a further
subsequence, we may assume $
s_n(y,z)\to0$, $
\gamma\text{-a.e.}$
Since the second marginal of $\gamma$ is $\eta^*$, we may also assume that
$g_n(z)\to g(z)$, $\gamma$-a.e. Hence, for $\gamma$-a.e. $(y,z)$,
\begin{equation}\label{eq:limit-Qgn}
Q_\lambda g_n(y)
=
\lambda|y-z|^2+g_n(z)-s_n(y,z)
\to
\lambda|y-z|^2+g(z).
\end{equation}

We next establish a pointwise upper bound for $g_n$. Since
$g_n=\Gamma^{\epsilon}_{f_n,h_n}$,
\[
e^{-g_n(u)/\epsilon}
=
\int_{\mathbb R^d}
\exp\left(
\frac{
f_n(x)+h_n(x)\cdot(u-x)-c(x,u)
}{\epsilon}
\right)\mu(dx)
\]
for every $u\in\mathbb R^d$. By the almost-everywhere convergence of
$f_n$ and $h_n$ and Fatou's lemma,
\[
\begin{aligned}
\liminf_{n\to\infty}e^{-g_n(u)/\epsilon}
&\geq
\int_{\mathbb R^d}
\exp\left(
\frac{
f(x)+h(x)\cdot(u-x)-c(x,u)
}{\epsilon}
\right)\mu(dx) \\
&=
e^{-g(u)/\epsilon}.
\end{aligned}
\]
Since $r\mapsto-\epsilon\log r$ is decreasing,
\begin{equation}\label{eq:limsup-gn}
\limsup_{n\to\infty}g_n(u)\leq g(u)
\qquad
\text{for every }u\in\mathbb R^d.
\end{equation}
Fix now $(y,z)$ in the $\gamma$-full-measure set on which \eqref{eq:limit-Qgn} holds. For every $u\in\mathbb R^d$, the definition of
$Q_\lambda g_n$ gives
\[
Q_\lambda g_n(y)
\leq
\lambda|y-u|^2+g_n(u).
\]
Taking the limit on the left and the limsup on the right, and using
\eqref{eq:limit-Qgn}-\eqref{eq:limsup-gn}, yields
\[
\lambda|y-z|^2+g(z)
\leq
\lambda|y-u|^2+g(u)
\qquad
\text{for every }u\in\mathbb R^d.
\]
Taking the infimum over $u$ gives $
\lambda|y-z|^2+g(z)
\leq
Q_\lambda g(y).$
The reverse inequality follows directly from the definition of
$Q_\lambda g$ by choosing $u=z$. Therefore
\begin{align}
Q_\lambda g(y)
=
\lambda|y-z|^2+g(z)
\qquad
\gamma\text{-a.e.}
\label{eq:equality-Q-lambda g}
\end{align}

Since $\gamma$ is an optimal coupling for $W_2^2(\eta^*,\nu)$, we have
$
\lambda|y-z|^2\in L^1(\gamma).$
Moreover, since $g\in L^1(\eta^*)$ and the second marginal of $\gamma$
is $\eta^*$, $g\in L^1(\gamma)$.
Therefore, $
\lambda|y-z|^2+g(z)\in L^1(\gamma).$
By the equality \eqref{eq:equality-Q-lambda g}, it follows that $Q_\lambda g(y)\in L^1(\gamma)$. Since the first
marginal of $\gamma$ is $\nu$, this is equivalent to
$Q_\lambda g\in L^1(\nu)$. Integrating the equality with respect to
$\gamma$ gives $
\nu(Q_\lambda g)
=
\lambda W_2^2(\eta^*,\nu)+\eta^*(g).$
\end{proof}

With the above ingredients in place, we are now ready to prove Theorem~\ref{thm:dual-attainment}.

\begin{proof}[Proof of Theorem~\ref{thm:dual-attainment}]
Set
$M:=M^{\lambda,\epsilon}(\mu,\nu)$ and
$D:=D^{\lambda,\epsilon}(\mu,\nu)$. By Theorem~\ref{thm:strong-duality},
$M=D$. We first consider the case where $\mu$ is non-Dirac. Recall that we denote $(\pi^*,\eta^*)$ to be an optimal WEMOT optimizer and $\pi^*$ is also the optimizer for $\operatorname{EMOT}(\mu, \eta^*)$. Set $T^*:=\lambda W_2^2(\eta^*,\nu)$ and $V^*:=\operatorname{EMOT}(\mu,\eta^*)$, so that $M=T^*+V^*$. 

By Proposition~\ref{prop:input-condition} and
Proposition~\ref{prop:emot-dual-attainment}, the fixed-marginal problem
$\operatorname{EMOT}(\mu,\eta^*)$ admits a dual optimizer
$(f_0,h_0)\in\mathcal D_{\operatorname{EMOT}}(\mu,\eta^*)$. Set
$g_0:=\Gamma^\epsilon_{f_0,h_0}$. Then
$V^*=\mu(f_0)+\eta^*(g_0)$. Define 
\begin{equation*}
r_0(x,z)
:=
\exp\biggl(
\frac{
f_0(x)+g_0(z)+h_0(x)\cdot(z-x)-c(x,z)
}{\epsilon}
\biggr).
\end{equation*} Therefore, by
Lemma~\ref{lem:fixed-marginal-entropy-gap}, $
\frac{d\pi^*}{d(\mu\otimes\eta^*)}(x,z)
=
r_0(x,z)$
and $r_0>0$, $\mu\otimes\eta^*$-a.e.
\medskip 

\noindent\emph{Step 1:}
Choose a maximizing sequence
$(f_n,h_n)\subset\mathcal D$ for the WEMOT dual problem and set
$g_n:=\Gamma^\epsilon_{f_n,h_n}$ and
$J_n:=\mu(f_n)+\nu(Q_\lambda g_n)$. Then $J_n\to D=M$. Since
$\mathcal D\subset\mathcal D_{\operatorname{EMOT}}(\mu,\eta^*)$, $
\mu(f_n)+\eta^*(g_n)\leq V^*.$
Let $\gamma^*\in\Pi(\nu,\eta^*)$ be an optimal coupling for
$W_2^2(\eta^*,\nu)$. Since
$Q_\lambda g_n(y)\leq\lambda|y-z|^2+g_n(z)$, we also have
$\nu(Q_\lambda g_n)\leq T^*+\eta^*(g_n)$. Define
\[
\Delta_n
:=
V^*-\mu(f_n)-\eta^*(g_n),
\qquad
\Sigma_n
:=
T^*+\eta^*(g_n)-\nu(Q_\lambda g_n).
\]
Then $\Delta_n,\Sigma_n\geq0$ and
$M-J_n=\Delta_n+\Sigma_n$. Hence $
\Delta_n\to0$ and $
\Sigma_n\to0.$

For each $n$, define $
r_n(x,z)
:=
\exp(
\frac{
f_n(x)+g_n(z)+h_n(x)\cdot(z-x)-c(x,z)
}{\epsilon}).$
Lemma~\ref{lem:fixed-marginal-entropy-gap} gives
$\epsilon H(\pi^*| r_n(\mu\otimes\eta^*))=\Delta_n\to0$.
Since $\pi^*=r_0(\mu\otimes\eta^*)$, Pinsker's inequality yields $
r_n\to r_0$ in $L^1(\mu\otimes\eta^*).$
Passing to a subsequence, we may assume that
$r_n\to r_0$ $(\mu\otimes\eta^*)$-a.e. Since $r_0>0$, setting
$v_n:=\epsilon\log r_n+c$ and $v_0:=\epsilon\log r_0+c$ gives
$v_n\to v_0$ $(\mu\otimes\eta^*)$-a.e., where $
v_n(x,z)
=
f_n(x)+g_n(z)+h_n(x)\cdot(z-x),$
and similarly
$v_0(x,z)=f_0(x)+g_0(z)+h_0(x)\cdot(z-x)$. 
\medskip

\noindent\emph{Step 2:}
By Lemma~\ref{lem:additive_compactness}, there exist
$\beta_n\in\mathbb R$ and $a_n\in\mathbb R^d$ such that
\[
f_n-f_0-\beta_n-a_n\cdot\operatorname{id}\to0,
\qquad
h_n-h_0-a_n\to0
\quad
\mu\text{-a.e.},
\]
and
\[
g_n-g_0+\beta_n+a_n\cdot\operatorname{id}\to0
\quad
\eta^*\text{-a.e.}
\]
Normalize the additive constants by setting
$\widetilde f_n:=f_n-\beta_n$,
$\widetilde g_n:=g_n+\beta_n$, and
$\widetilde h_n:=h_n$. Then
$\widetilde g_n=\Gamma^\epsilon_{\widetilde f_n,\widetilde h_n}$ and $
\mu(\widetilde f_n)
+
\nu(Q_\lambda\widetilde g_n)
=
J_n,$
while $\Delta_n$ and $\Sigma_n$ are unchanged. Moreover,
\[
\widetilde f_n-f_0-a_n\cdot\operatorname{id}\to0,
\qquad
\widetilde h_n-h_0-a_n\to0
\quad
\mu\text{-a.e.},
\]
and
\[
\widetilde g_n-g_0+a_n\cdot\operatorname{id}\to0
\quad
\eta^*\text{-a.e.}
\]

Define the Wasserstein slack
$s_n(y,z):=\lambda|y-z|^2+\widetilde g_n(z)
-Q_\lambda\widetilde g_n(y)$. Then $s_n\geq0$ and
\[
\int_{\mathbb R^d\times\mathbb R^d}s_n(y,z)\,\gamma^*(dy,dz)
=
\Sigma_n
\to0.
\]
Lemma~\ref{lem:affine-slope} therefore implies that $(a_n)_n$ is
bounded. Passing to a further subsequence, let $a_n\to a\in\mathbb R^d$.
It follows that
\[
\widetilde f_n\to f^*
:=
f_0+a\cdot\operatorname{id},
\qquad
\widetilde h_n\to h^*
:=
h_0+a
\quad
\mu\text{-a.e.},
\]
and
\[
\widetilde g_n\to g^*
:=
g_0-a\cdot\operatorname{id}
\quad
\eta^*\text{-a.e.}
\]
The affine transformation is compatible with $\Gamma$, and hence $
\Gamma^\epsilon_{f^*,h^*}
=
g^*.$

\medskip

\noindent\emph{Step 3:}
Since $g_0\in L^1(\eta^*)$ and $\eta^*\in\mathcal P_2(\mathbb R^d)$,
we have $g^*\in L^1(\eta^*)$. We may therefore apply
Lemma~\ref{lem:wasserstein-complementarity}, which gives $
Q_\lambda g^*(y)
=
\lambda|y-z|^2+g^*(z)$ $
\gamma^*\text{-a.e.}$
and $
\nu(Q_\lambda g^*)
=
T^*+\eta^*(g^*).$
Since $\pi^*\in\Pi_M(\mu,\eta^*)$, the measures $\mu$ and $\eta^*$
have the same barycenter. Therefore
\[
\mu(f^*)+\eta^*(g^*)
=
\mu(f_0)+\eta^*(g_0)
=
V^*.
\]
Consequently,
\[
\mu(f^*)
+
\nu(Q_\lambda\Gamma^\epsilon_{f^*,h^*})
=
V^*+T^*
=
D.
\]

\medskip

\noindent\emph{Step 4:}
We claim that $(f^*, h^*)\in \mathcal D$. First, since
$f^*=f_0+a\cdot\operatorname{id}$ and $\mu\in\mathcal P_2(\mathbb R^d)$,
$f_0\in L^1(\mu)$ implies that $f^*\in L^1(\mu)$. 

Next, we show that $g^*\in L^1(\eta)$ for every $\eta\in\mathcal P_2(\mathbb R^d)$, and that $Q_\lambda g^*\in L^1(\nu)$. By Lemma~\ref{lem:wasserstein-complementarity}, we have $Q_\lambda g^*\in L^1(\nu)$. It remains to show that $g^*\in L^1(\eta)$ for every $\eta\in\mathcal P_2(\mathbb R^d)$. Choose $y_0\in\mathbb R^d$ such that $Q_\lambda g^*(y_0)\in\mathbb R$.  The definition of $Q_\lambda$ then implies
\[
g^*(z)\geq Q_\lambda g^*(y_0)-\lambda|y_0-z|^2,\qquad z\in\mathbb R^d,
\]
which provides a global quadratic lower bound. For the upper bound, since $f^*$ and $h^*$ are finite $\mu$-a.e., there exist a measurable set $E\subset\operatorname{supp}(\mu)$ with $\mu(E)>0$ and a constant $C>0$ such that $|f^*(x)|+|h^*(x)|\leq C$ on $E$. By compactness of $\operatorname{supp}(\mu)$ and boundedness of $c$, enlarging $C$ if necessary gives
\[
f^*(x)+h^*(x)\cdot(z-x)-c(x,z)\geq-C(1+|z|),\qquad x\in E,\ z\in\mathbb R^d.
\]
Restricting the integral defining $g^*=\Gamma^\epsilon_{f^*,h^*}$ to $E$, we obtain $e^{-g^*(z)/\epsilon}\geq\mu(E)e^{-C(1+|z|)/\epsilon}$, and hence
\[
g^*(z)\leq C(1+|z|)-\epsilon\log\mu(E),\qquad z\in\mathbb R^d.
\]
Combining these bounds yields, for some constant $C'>0$, $
|g^*(z)|\leq C'(1+|z|^2),$ for all $z\in\mathbb R^d.$
Thus, for every $\eta\in\mathcal P_2(\mathbb R^d)$, $g^*\in L^1(\eta)$, as claimed. 

Finally, we verify the admissibility condition for $h^*$. 
The definition of $\Gamma^\epsilon_{f^*,h^*}$ and the
boundedness of $c$ imply
\begin{equation}\label{eq: gaussian-bound}
\int_{\mathbb R^d}
\exp\left(
\frac{f^*(x)+h^*(x)\cdot(z-x)}{\epsilon}
\right)\mu(dx)
\leq
\exp\left(
\frac{\|c\|_\infty-g^*(z)}{\epsilon}
\right).
\end{equation}
As $g^*(z)\geq-C'(1+|z|^2)$, we fix some constant $A>C'$ and set $
T(x):=
f^*(x)-h^*(x)\cdot x+\frac{|h^*(x)|^2}{4A}.$ Multiplying \eqref{eq: gaussian-bound} by
$e^{-A|z|^2/\epsilon}$ and integrating over $z$ yields
\[
\int_{\mathbb R^d} e^{T(x)/\epsilon}\,\mu(dx)
\leq
\exp\left(\frac{\|c\|_\infty+C'}{\epsilon}\right)
\left(\frac{A}{A-C'}\right)^{d/2}
<\infty.
\]
Since $t^+\leq \epsilon e^{t/\epsilon}$ for every
$t\in\mathbb R$, we obtain $T^+\in L^1(\mu)$. By Young's inequality, $
|h^*(x)\cdot x|
\leq
\frac{|h^*(x)|^2}{8A}+2A|x|^2,$
and hence $
\frac{|h^*(x)|^2}{8A}
\leq
T^+(x)+(f^*(x))^-+2A|x|^2.$
Because $f^*\in L^1(\mu)$ and
$\mu\in\mathcal P_2(\mathbb R^d)$, this proves that
$h^*\in L^2(\mu)$.
Consequently, for every $\eta\in\mathcal P_2(\mathbb R^d)$ and
$\pi\in\Pi_M(\mu,\eta)$,
\[
\int |h^*(x)\cdot(z-x)|\,d\pi
\leq
\|h^*\|_{L^2(\mu)}
\left(\int |z-x|^2\,d\pi\right)^{1/2}<\infty.
\]
Thus $h^*$ satisfies the required admissibility condition
for every such martingale coupling, in particular for
those satisfying $\pi\ll\mu\otimes\eta$. Hence
$(f^*,h^*)\in\mathcal D$.

\medskip
The case where $\mu$ is a Dirac measure is treated separately below. Consider the case $\mu=\delta_m$ for some
$m\in\mathbb R^d$. For $a\in\mathbb R^d$, define
$f_a\equiv0$ and $h_a\equiv a$. Then 
\[
\Gamma^{\epsilon}_{f_a,h_a}(z)
=
c(m,z)-a\cdot(z-m)
=:g_a(z).
\]
Define $F(a):=\nu(Q_\lambda g_a)$, where $
Q_\lambda g_a(y)
=
\inf_{z\in\mathbb R^d}\{
\lambda|y-z|^2+c(m,z)-a\cdot(z-m)\}.$
Then we have
\[
-\|c\|_{\infty}-a\cdot(y-m)-\frac{|a|^2}{4\lambda}
\leq
Q_\lambda g_a(y)
\leq
\|c\|_{\infty}-a\cdot(y-m)-\frac{|a|^2}{4\lambda}.\]
Hence $Q_\lambda g_a\in L^1(\nu)$ and $
F(a)
\leq
\|c\|_{\infty}
-
a\cdot(m_\nu-m)
-
\frac{|a|^2}{4\lambda},$
where $m_\nu:=\int y\,\nu(dy)$. Therefore
$F(a)\to-\infty$ as $|a|\to\infty$. Moreover, $F$ is concave, since
for each $y$ the map $a\mapsto Q_\lambda g_a(y)$ is the pointwise
infimum of affine functions of $a$. Since $F$ is finite on
$\mathbb R^d$, it is continuous. Consequently, there exists
$a^*\in\mathbb R^d$ such that $
F(a^*)=\sup_{a\in\mathbb R^d}F(a).$

For $\mu=\delta_m$, every admissible pair $(f,h)$ enters the dual
objective only through $f(m)$ and $h(m)$. Indeed, setting
$a:=h(m)\in\mathbb R^d$, one has
\[
\Gamma^{\epsilon}_{f,h}(z)
=
c(m,z)-f(m)-a\cdot(z-m)
=
g_a(z)-f(m).
\]
Since $Q_\lambda\Gamma^{\epsilon}_{f,h}=Q_\lambda g_a-f(m)$, the constant
$f(m)$ cancels between $\mu(f)$ and
$\nu(Q_\lambda\Gamma^{\epsilon}_{f,h})$. Hence $
D=\sup_{a\in\mathbb R^d}F(a).$
By strong duality, $F(a^*)=M$. Taking
$f^*\equiv0$ and $h^*\equiv a^*$, we obtain
\[
\mu(f^*)
+
\nu(Q_\lambda\Gamma^\epsilon_{f^*,h^*})
=D.
\]
Thus $(f^*,h^*)$ attains the dual problem also in the Dirac case.
\end{proof}

\section{Stability and Convergence}\label{sec:quantitative_stability}

\subsection{Quantitative stability under marginal perturbations}

In this section, we prove the quantitative stability estimate stated in
Theorem~\ref{thm:quantitative-stability}. For fixed $\lambda,\epsilon > 0$
and arbitrary marginals
$\mu,\nu,\tilde{\mu},\tilde{\nu} \in \mathcal{P}_2(\mathbb{R}^d)$,
we bound $
|M^{\lambda,\epsilon}(\mu,\nu)
- M^{\lambda,\epsilon}(\tilde{\mu},\tilde{\nu})|$
in terms of $W_2(\mu,\tilde{\mu})$ and $W_2(\nu,\tilde{\nu})$,
under the Lipschitz condition on $c$ in
Assumption~\ref{ass:c-lipschitz}, for any dimension $d \geq 1$.
We begin by recalling the objective function and this regularity assumption.

For $\eta\in\mathcal P_2(\mathbb R^d)$ and
$\pi\in\Pi_M(\mu,\eta)$, recall that the objective function in \eqref{eq:objective-function} is defined as
\[
J_{\mu, \nu}(\pi,\eta)
:=
\int_{\mathbb R^d\times\mathbb R^d} c(x,y)\,\pi(dx,dy)
+
\epsilon H(\pi\mid\mu\otimes\eta)
+
\lambda W_2^2(\eta,\nu),
\]
where $c:\mathbb R^d\times\mathbb R^d\to\mathbb R$ is the cost
function and $\Pi_M(\mu,\eta)$ denotes the set of martingale
couplings of $\mu$ and $\eta$. Thus,
\[
M^{\lambda,\epsilon}(\mu,\nu)
=
\inf_{\substack{\eta\in\mathcal P_2(\mathbb R^d)\\
\pi\in\Pi_M(\mu,\eta)}}
J_{\mu, \nu}(\pi,\eta).
\]
Recall also that Assumption~\ref{ass:c-lipschitz} requires
the existence of a constant $L_c>0$ such that
\[
|c(x,y)-c(x',y')|
\leq
L_c\bigl(|x-x'|+|y-y'|\bigr)
\]
for all $(x,y),(x',y')\in\mathbb R^d\times\mathbb R^d$.
This condition will be used to control the variation of the cost
term. 

To start with, We recall the following Gaussian mutual information bound \cite[Section~6.2]{polyanskiy2025information}.
This bound will be used to control the entropy term when constructing
a martingale coupling with  perturbed marginals.
For completeness, we state the lemma as follows. 

For random variables $U$ and $V$, we write the mutual information between $U$ and $V$ as 
\[
I(U;V)
:=
H\bigl(\mathcal L(U,V)
\mid \mathcal L(U)\otimes\mathcal L(V)\bigr),
\]
where $H$ denotes the relative entropy defined in \eqref{eq:relative-entropy} and $\mathcal L(\cdot)$ denotes
the law of a random variable.

\begin{lemma}[{\citealp[Section~6.2]{polyanskiy2025information}}]
\label{lem:gaussian_information}
Let $X$ be an $\mathbb R^d$-valued random variable satisfying
$\mathbb E|X|^2<\infty$, and let $G\sim N(0,I_d)$ be independent of $X$.
Then, for every $\sigma>0$,
\begin{equation}
I(X;X+\sigma G)\leq
\frac d2
\log\Big(
1+\frac{\mathbb E|X|^2}{d\sigma^2}
\Big)
\leq
\frac{\mathbb E|X|^2}{2\sigma^2}.
\label{eq:gaussian_information}
\end{equation}
\end{lemma}

We next consider how to modify a martingale coupling when its first
marginal is perturbed. Starting from $\pi\in\Pi_M(\mu,\eta)$ for some $\eta\in \mathcal{P}_2(\mathbb{R}^d)$ and a given perturbed marginal $\tilde{\mu}\in \mathcal{P}_2(\mathbb{R}^d)$, the following lemma constructs a coupling
$\tilde{\pi}\in\Pi_M(\tilde{\mu},\tilde{\eta})$.
The Gaussian mutual information bound controls the increase in
relative entropy, while the construction also yields bounds on
$W_2(\eta,\tilde{\eta})$.

\begin{lemma}
\label{lem:martingale_repair}
Let $\mu,\widetilde\mu,\eta\in\mathcal P_2(\mathbb R^d)$ and let
$\pi\in\Pi_M(\mu,\eta)$ satisfy $ H(\pi\mid\mu\otimes\eta)<\infty.$
Write $
\delta:=W_2(\mu,\widetilde\mu).$
Then, for every $\sigma>0$, there exist
$\widetilde\eta\in\mathcal P_2(\mathbb R^d)$ and
$\widetilde\pi\in\Pi_M(\widetilde\mu,\widetilde\eta)$ such that
\begin{align}
W_2(\eta,\widetilde\eta)
&\leq
\sqrt{\delta^2+d\sigma^2},
\label{eq:repair_eta}
\\
H(\widetilde\pi\mid
\widetilde\mu\otimes\widetilde\eta)
&\leq
H(\pi\mid\mu\otimes\eta)
+
\frac{\delta^2}{2\sigma^2}.
\label{eq:repair_entropy}
\end{align}
Moreover, if $c$ satisfies Assumption \ref{ass:c-lipschitz}, 
then
\begin{equation}
\left|
\int c\,d\widetilde\pi-\int c\,d\pi
\right|
\leq
L_c
\left(
\delta+\sqrt{\delta^2+d\sigma^2}
\right).
\label{eq:repair_cost}
\end{equation}
\end{lemma}

\begin{proof}
Let $(X,\widetilde X)$ be an optimal quadratic Wasserstein coupling of
$\mu$ and $\widetilde\mu$, so that $
\mathbb E|X-\widetilde X|^2
=
\delta^2.$
On the other hand, let $(X,Y)$ have law $\pi$. By the gluing lemma, these two couplings may be realized on a common probability space in such a way that $\widetilde X$ and $Y$ are conditionally independent given $X$. Equivalently, $
\widetilde X-X-Y$
forms a Markov chain. Since $\pi$ is a martingale coupling,
$\mathbb E[Y\mid X]=X$, and therefore
\begin{align}
\mathbb E[Y\mid\widetilde X]
&=
\mathbb E\left[
\mathbb E[Y\mid X,\widetilde X]
\,\middle|\,
\widetilde X
\right]=
\mathbb E\left[
\mathbb E[Y\mid X]
\,\middle|\,
\widetilde X
\right]
=
\mathbb E[X\mid\widetilde X].
\label{eq:conditional_Y}
\end{align}

Introduce the martingale defect $
A
:=
\widetilde X-\mathbb E[X\mid\widetilde X].$
Since $
A
=
\mathbb E[\widetilde X-X\mid\widetilde X],$
conditional Jensen's inequality yields
\begin{equation}
\mathbb E|A|^2
\leq
\mathbb E|\widetilde X-X|^2
=
\delta^2.
\label{eq:A_bound}
\end{equation}
Let $G\sim N(0,I_d)$ be independent of
$(X,\widetilde X,Y)$ and define $
\Xi:=A+\sigma G,$ $
\widetilde Y:=Y+\Xi.$
Using \eqref{eq:conditional_Y} and the fact that $G$ is centered and independent of $\widetilde X$, we obtain
\[
\mathbb E[\widetilde Y\mid\widetilde X]
=
\mathbb E[X\mid\widetilde X]
+
\widetilde X-\mathbb E[X\mid\widetilde X]
=
\widetilde X.
\]
Hence, with $
\widetilde\pi
:=
\operatorname{Law}(\widetilde X,\widetilde Y),
$, $
\widetilde\eta
:=
\operatorname{Law}(\widetilde Y),$
we have
$\widetilde\pi\in\Pi_M(\widetilde\mu,\widetilde\eta)$.

We next estimate $W_2^2(\eta, \widetilde\eta)$. Since $G$ is independent of $A$ and has mean zero,
\[
\mathbb E|\widetilde Y-Y|^2
=
\mathbb E|A+\sigma G|^2
=
\mathbb E|A|^2+d\sigma^2
\leq
\delta^2+d\sigma^2.
\]
The joint law of $(Y,\widetilde Y)$ is therefore an admissible coupling of
$\eta$ and $\widetilde\eta$, and \eqref{eq:repair_eta} follows.

We now prove the entropy estimate \eqref{eq:repair_entropy}. For random variables $(U,V)$, write $
I(U;V)
:=
H(
\mathbb{P}_{U,V}\mid \mathbb{P}_U\otimes \mathbb{P}_V
)$
for their mutual information. Consequently, $H(\pi\mid\mu\otimes\eta)
=
I(X;Y)$ and  $
 H(
\widetilde\pi
\mid
\widetilde\mu\otimes\widetilde\eta)
=
I(\widetilde X;\widetilde Y).$
Since $\widetilde Y=Y+\Xi$ is a measurable function of $(Y,\Xi)$, the data processing inequality of mutual information \cite[Theorem~2.8.1]{CoverThomas2006} gives $
I(\widetilde X;\widetilde Y)
\leq
I(\widetilde X;(Y,\Xi)).$
By the chain rule for mutual information, $
I(\widetilde X;(Y,\Xi))
=
I(\widetilde X;Y)
+
I(\widetilde X;\Xi\mid Y).$
The Markov property $\widetilde X-X-Y$ implies, again by the data processing inequality,
\begin{equation}
I(\widetilde X;Y)
\leq
I(X;Y).
\label{eq:DP_first}
\end{equation}
Conditional on $\widetilde X$, the variable $\Xi=A+\sigma G$ depends only on the independent Gaussian variable $G$. Hence $
Y\perp\!\!\!\perp\Xi\mid\widetilde X$, 
and $
I(Y;\Xi\mid\widetilde X)=0.$
Therefore applying the chain rule \cite[Theorem~2.5.2]{CoverThomas2006}, to $I((\widetilde X,Y);\Xi)$  gives $I((\widetilde X,Y);\Xi)=
I(\widetilde X;\Xi)+I(Y;\Xi\mid\widetilde X)=I(\widetilde X;\Xi)$. Then applying the chain rule in the other direction yields
\[I(\widetilde X;\Xi)=I((\widetilde X,Y);\Xi)
=
I(Y;\Xi)
+
I(\widetilde X;\Xi\mid Y)
\geq
I(\widetilde X;\Xi\mid Y),
\]
so that
\begin{equation}
I(\widetilde X;\Xi\mid Y)
\leq
I(\widetilde X;\Xi).
\label{eq:conditional_information_bound}
\end{equation}
Combining \eqref{eq:DP_first} and \eqref{eq:conditional_information_bound}, we have 
\begin{align}
I(\widetilde{X}; \widetilde{Y})\leq I(X; Y)+I(\widetilde{X}; \Xi)\label{eq:entropy_bound_1}
\end{align}

We next provide a bound for $I(\widetilde{X}; \Xi)$. By the chain rule for mutual information,
\[
I((\widetilde X,A);\Xi)
=
I(\widetilde X;\Xi)
+
I(A;\Xi\mid\widetilde X).
\]
Since $A$ is completely determined by $\widetilde X$, we have $
I(A;\Xi\mid\widetilde X)=0,$
and hence $
I((\widetilde X,A);\Xi)
=
I(\widetilde X;\Xi).$
On the other hand, applying the chain rule in the other direction gives
\[
I((\widetilde X,A);\Xi)
=
I(A;\Xi)
+
I(\widetilde X;\Xi\mid A).
\]
Recall that $
\Xi=A+\sigma G,$
where $G$ is independent of $\widetilde X$. Since $A$ is a measurable function of $\widetilde X$, $G$ is also independent of the pair $(\widetilde X,A)$. Consequently, conditional on $A$, the law of $\Xi$ is $
\mathcal L(\Xi\mid A)
=
N(A,\sigma^2 I_d),$
and does not depend further on $\widetilde X$. Equivalently, $
\widetilde X\perp\!\!\!\perp\Xi\mid A,$
so that $
I(\widetilde X;\Xi\mid A)=0.$
Combining the above equations yields $
I(\widetilde X;\Xi)
=
I(A;\Xi).$
Since $\Xi=A+\sigma G$, Lemma~\ref{lem:gaussian_information} and Equation
\eqref{eq:A_bound} imply
\begin{equation}\label{eq:entropy_bound_2}
I(\widetilde X;\Xi)
=
I(A;A+\sigma G)
\leq
\frac{\mathbb E|A|^2}{2\sigma^2}
\leq
\frac{\delta^2}{2\sigma^2}.
\end{equation}
Combining
\eqref{eq:entropy_bound_1} and
\eqref{eq:entropy_bound_2} yields $
I(\widetilde X;\widetilde Y)
\leq
I(X;Y)
+
\frac{\delta^2}{2\sigma^2},$
which is precisely \eqref{eq:repair_entropy}.

Finally, under the Lipschitz assumption on $c$,
\begin{align*}
\left|
\mathbb E[c(\widetilde X,\widetilde Y)]
-
\mathbb E[c(X,Y)]
\right|\leq
L_c
\left(
\mathbb E|\widetilde X-X|
+
\mathbb E|\widetilde Y-Y|
\right)\leq
L_c
\left(
\delta+\sqrt{\delta^2+d\sigma^2}
\right),
\end{align*}
where the last inequality follows from Cauchy--Schwarz. This proves
\eqref{eq:repair_cost}.
\end{proof}

We can now derive a quantitative stability estimate for the primal
problem. Recall from \eqref{eq:K estimate of stability} that
\[
K(\mu,\nu)
:=
\Big(
\frac{2\|c\|_\infty}{\lambda}
+
\frac{\epsilon}{2\lambda}m_2(\mu)
+
\big(\sqrt d+W_2(\mu,\nu)\big)^2
\Big)^{1/2}.
\]
We now state a version of Theorem~\ref{thm:quantitative-stability}
with an explicit quantitative rate.

\begin{stabilityrestatement}
Suppose Assumption~\ref{ass:c-lipschitz} holds. Fix $\lambda,\epsilon>0$, let $\mu,\widetilde\mu,\nu,\widetilde\nu\in\mathcal P_2(\mathbb R^d)$, and define
\[
\delta:=W_2(\mu,\widetilde\mu),\qquad \rho:=W_2(\nu,\widetilde\nu),\qquad \overline K:=\max\{K(\mu,\nu),K(\widetilde\mu,\widetilde\nu)\}.
\]
Then
\begin{equation}\label{eq:quantitative_stability}
\begin{aligned}
\bigl|M^{\lambda,\epsilon}(\widetilde\mu,\widetilde\nu)-M^{\lambda,\epsilon}(\mu,\nu)\bigr|
&\leq L_c\bigl(\delta+\sqrt{\delta^2+d\delta^{4/3}}\bigr)+\frac{\epsilon}{2}\delta^{2/3}\\
&\quad+2\lambda\overline K\bigl(\sqrt{\delta^2+d\delta^{4/3}}+\rho\bigr)+\lambda\bigl(\sqrt{\delta^2+d\delta^{4/3}}+\rho\bigr)^2.
\end{aligned}
\end{equation}
Consequently, setting $s(\delta,\rho):=\delta^{2/3}+\delta+\rho$, we have
\[
\bigl|M^{\lambda,\epsilon}(\widetilde\mu,\widetilde\nu)-M^{\lambda,\epsilon}(\mu,\nu)\bigr|\leq C\bigl\{s(\delta,\rho)+s(\delta,\rho)^2\bigr\},
\]
where one may take
\[
C:=\max\bigl\{(1+\sqrt d)L_c+\epsilon/2+2\lambda\overline K\sqrt d,\ \lambda d\bigr\}.
\]
For fixed $(\mu,\nu)$, the definition of $K$ ensures that $\overline K$ remains bounded as $\delta,\rho\to0$. In particular,
\[
\bigl|M^{\lambda,\epsilon}(\widetilde\mu,\widetilde\nu)-M^{\lambda,\epsilon}(\mu,\nu)\bigr|=O\bigl(\delta^{2/3}+\rho\bigr).
\]
\end{stabilityrestatement}

\begin{proof}
By Theorem~\ref{thm:primal-attainment}, choose an optimal pair $(\pi^*,\eta^*)$ for $M^{\lambda,\epsilon}(\mu,\nu)$ and set $r:=W_2(\eta^*,\nu)$. Since $c\geq-\|c\|_\infty$ and relative entropy is nonnegative, \eqref{eq:primal-value-bound} and the definition of $K$ give
\[
\lambda r^2\leq M^{\lambda,\epsilon}(\mu,\nu)+\|c\|_\infty\leq\lambda K(\mu,\nu)^2.
\]
Hence $r\leq K(\mu,\nu)\leq\overline K$. Moreover, the optimal pair has finite relative entropy.

Fix $\sigma>0$ and set $s:=\sqrt{\delta^2+d\sigma^2}$. Applying Lemma~\ref{lem:martingale_repair} to $\pi^*$ and the perturbed first marginal $\widetilde\mu$ gives $\widetilde\eta\in\mathcal P_2(\mathbb R^d)$ and $\widetilde\pi\in\Pi_M(\widetilde\mu,\widetilde\eta)$ satisfying $W_2(\eta^*,\widetilde\eta)\leq s$. Consequently,
\[
W_2(\widetilde\eta,\widetilde\nu)\leq r+s+\rho.
\]
Combining the Wasserstein penalty estimate with the entropy and cost estimates in Lemma~\ref{lem:martingale_repair}, we obtain
\[
M^{\lambda,\epsilon}(\widetilde\mu,\widetilde\nu)-M^{\lambda,\epsilon}(\mu,\nu)\leq J_{\widetilde\mu,\widetilde\nu}(\widetilde\pi,\widetilde\eta)-J_{\mu,\nu}(\pi^*,\eta^*)\leq L_c(\delta+s)+\frac{\epsilon\delta^2}{2\sigma^2}+2\lambda\overline K(s+\rho)+\lambda(s+\rho)^2.
\]

Interchanging $(\mu,\nu)$ and
$(\widetilde\mu,\widetilde\nu)$ yields the same bound
for the opposite difference, since $\delta$, $\rho$,
and $\overline K$ are unchanged. Therefore, for every $\sigma>0$,
{\small
\begin{equation}\label{eq:quantitative_stability_sigma}
\bigl|M^{\lambda,\epsilon}(\widetilde\mu,\widetilde\nu)-M^{\lambda,\epsilon}(\mu,\nu)\bigr|
\leq L_c(\delta+\sqrt{\delta^2+d\sigma^2})+\frac{\epsilon\delta^2}{2\sigma^2}
+2\lambda\overline K(\sqrt{\delta^2+d\sigma^2}+\rho)
+\lambda(\sqrt{\delta^2+d\sigma^2}+\rho)^2.
\end{equation}
}
If $\delta>0$, choosing $\sigma=\delta^{2/3}$ in
\eqref{eq:quantitative_stability_sigma} gives
\eqref{eq:quantitative_stability}. If $\delta=0$, letting $\sigma\downarrow0$ instead yields
\[
\left|
M^{\lambda,\epsilon}(\widetilde\mu,\widetilde\nu)
-
M^{\lambda,\epsilon}(\mu,\nu)
\right|
\leq
2\lambda\overline K\rho+\lambda\rho^2,
\]
which is precisely \eqref{eq:quantitative_stability}
in this case.

By \eqref{eq:quantitative_stability} and the inequality $
\sqrt{\delta^2+d\delta^{4/3}}
\leq \delta+\sqrt d\,\delta^{2/3},$
we obtain
\begin{align}\label{eq:rought rate}
\big|
M^{\lambda,\epsilon}(\widetilde\mu,\widetilde\nu)
-
M^{\lambda,\epsilon}(\mu,\nu)\big|\leq &
L_c\big(2\delta+\sqrt d\,\delta^{2/3}\big)
+\frac{\epsilon}{2}\delta^{2/3}+
2\lambda\overline K
\big(\delta+\sqrt d\,\delta^{2/3}+\rho\big)\nonumber\\
&+\lambda
\big(\delta+\sqrt d\,\delta^{2/3}+\rho\big)^2.
\end{align}
Since $d\geq1$, $
\delta+\sqrt d\,\delta^{2/3}+\rho
\leq \sqrt d\,s.$
After collecting terms, \eqref{eq:rought rate} therefore becomes
\[
\begin{aligned}
&\big|
M^{\lambda,\epsilon}(\widetilde\mu,\widetilde\nu)
-
M^{\lambda,\epsilon}(\mu,\nu)\big|\leq
\big(2L_c+2\lambda\overline K\big)\delta
+
\Big(
\sqrt d\big(L_c+2\lambda\overline K\big)
+\frac{\epsilon}{2}
\Big)\delta^{2/3}
+
2\lambda\overline K\rho
+
\lambda d\,s^2.
\end{aligned}
\]
Thus, defining
\[
C:=
\max\Big\{
2L_c+2\lambda\overline K,\,
\sqrt d\big(L_c+2\lambda\overline K\big)
+\frac{\epsilon}{2},\,
\lambda d
\Big\},
\]
we conclude that
\[
\big|
M^{\lambda,\epsilon}(\widetilde\mu,\widetilde\nu)
-
M^{\lambda,\epsilon}(\mu,\nu)
\big|
\leq C(s+s^2).
\]

Finally, fix $(\mu,\nu)$ and let $\delta,\rho\to 0$.
By the triangle inequality,
$m_2(\widetilde\mu)^{1/2}\leq m_2(\mu)^{1/2}+\delta$
and
$W_2(\widetilde\mu,\widetilde\nu)\leq W_2(\mu,\nu)+\delta+\rho$.
Thus, both quantities remain bounded, so the definition of $K$ in \eqref{eq:K estimate of stability}
implies that $K(\widetilde\mu,\widetilde\nu)$ remains bounded.
Since $K(\mu,\nu)$ is fixed, $\overline K$ and $C$ also remain bounded. Moreover, for $\delta\leq1$,  $
s\leq2\delta^{2/3}+\rho.$
Since $s\to0$, it follows that
\[
\big|
M^{\lambda,\epsilon}(\widetilde\mu,\widetilde\nu)
-
M^{\lambda,\epsilon}(\mu,\nu)
\big|
=
O\,\big(\delta^{2/3}+\rho\big).
\]
This completes the proof.
\end{proof}

\subsection{Convergence to the martingale optimal transport problem}

In this subsection, we prove the convergence results stated in
Theorem~\ref{thm:convergence-to-mot}. Recall that $\lambda$ penalizes deviations
of the relaxed second marginal from $\nu$, while $\epsilon$
controls the strength of the entropic regularization.
This suggests that the classical martingale optimal transport
problem can be recovered by letting $\lambda\to\infty$ and
$\epsilon\to0$. The two limits play different roles and require
separate arguments. We first let $\epsilon\to0$ with $\lambda>0$
fixed to remove the entropic regularization, and then let
$\lambda\to\infty$ to recover the second marginal constraint.

Recall from \eqref{eq:wmot} that the Wasserstein-relaxed
martingale optimal transport value is defined by
\[
M^\lambda(\mu,\nu)
:=
\inf_{\eta\in\mathcal P_2(\mathbb R^d)}
\left\{
\lambda W_2^2(\eta,\nu)
+
\inf_{\pi\in\Pi_M(\mu,\eta)}
\int_{\mathbb R^d\times\mathbb R^d}
c(x,y)\,\pi(dx,dy)
\right\}.
\]
Our first step is to show that, for every fixed $\lambda>0$, $
M^{\lambda,\epsilon}(\mu,\nu)
\xrightarrow{\epsilon\to0}
M^\lambda(\mu,\nu).$
We then combine this vanishing-entropy limit with the convergence of the Wasserstein-relaxed problem, $
M^\lambda(\mu,\nu)\xrightarrow{\lambda\to\infty}M(\mu,\nu).$
Consequently, we prove that $
\lim_{\lambda\to\infty}\lim_{\epsilon\to0^+}
M^{\lambda,\epsilon}(\mu,\nu)
=
M(\mu,\nu).$
The convergence result holds in arbitrary dimension. Under additional regularity assumptions, we further obtain quantitative estimates for the approximation error and derive explicit convergence rates with respect to $\lambda$ and $\epsilon$.

\begin{proposition}\label{prop:epsilon_convergence}
Suppose Assumption~\ref{ass:c-lipschitz} holds. Let $\mu,\nu\in\mathcal P_2(\mathbb R^d)$ satisfy $\mu\preceq_{cx}\nu$. Then, for every fixed $\lambda>0$,
\[
\lim_{\epsilon\to0^+} M^{\lambda,\epsilon}(\mu,\nu)=M^\lambda(\mu,\nu).
\]
More precisely, set $C_\lambda:=(\sqrt{m_2(\nu)}+\sqrt{2\|c\|_\infty/\lambda})^2$ and $C_d:=\mathbb E|G|=\frac{\sqrt{2}\,\Gamma((d+1)/2)}{\Gamma(d/2)}$, where $G\sim\mathcal N(0,I_d)$. Then, for every $\epsilon>0$ and $\sigma>0$,
\begin{equation}\label{eq:entropy_approximation_bound}
0\leq M^{\lambda,\epsilon}(\mu,\nu)-M^\lambda(\mu,\nu)
\leq \lambda d\sigma^2+C_dL_c\sigma+\frac{\epsilon d}{2}\log\biggl(1+\frac{C_\lambda}{d\sigma^2}\biggr).
\end{equation}
In particular, for fixed $\lambda>0$, choosing $\sigma=\epsilon$ gives $M^{\lambda,\epsilon}(\mu,\nu)-M^\lambda(\mu,\nu)=O\bigl(\epsilon|\log\epsilon|\bigr)$ as $\epsilon\to0^+$.
\end{proposition}

\begin{proof}
Since relative entropy is nonnegative, $M^{\lambda,\epsilon}(\mu,\nu)\geq M^\lambda(\mu,\nu)$. It therefore suffices to prove the upper bound.

Fix $\lambda>0$. Since $\mu\preceq_{cx}\nu$, Strassen's theorem gives $\Pi_M(\mu,\nu)\neq\varnothing$. Taking $\eta=\nu$ and using the boundedness of $c$, we obtain $-\|c\|_\infty\leq M^\lambda(\mu,\nu)\leq\|c\|_\infty$. For every $\delta>0$, choose an almost optimal pair $(\eta, \pi)$ such that  $\eta\in\mathcal P_2(\mathbb R^d)$, $\pi\in\Pi_M(\mu,\eta)$ and
\begin{equation}\label{eq:delta_optimal_pair}
\lambda W_2^2(\eta,\nu)+\int c\,d\pi\leq M^\lambda(\mu,\nu)+\delta.
\end{equation}
Since $\int c\,d\pi\geq-\|c\|_\infty$, we have $W_2^2(\eta,\nu)\leq(2\|c\|_\infty+\delta)/\lambda$. Consequently,
\begin{equation}\label{eq:eta_moment_bound}
\sqrt{m_2(\eta)}\leq\sqrt{m_2(\nu)}+W_2(\eta,\nu)\leq\sqrt{m_2(\nu)}+\sqrt{\frac{2\|c\|_\infty+\delta}{\lambda}}.
\end{equation}

Let $(X,Y)$ have law $\pi$. By gluing an optimal coupling of $\eta$ and $\nu$ to $(X,Y)$, choose $V\sim\nu$ such that $\mathbb E|V-Y|^2=W_2^2(\nu,\eta)$. Let $G\sim\mathcal N(0,I_d)$ be independent of $(X,Y,V)$. For $\sigma>0$, define
\[
Y^\sigma:=Y+\sigma G,\qquad \pi^\sigma:=\mathcal L(X,Y^\sigma),\qquad \eta^\sigma:=\mathcal L(Y^\sigma).
\]
Since $\eta^\sigma\in\mathcal P_2(\mathbb R^d)$ and $\mathbb E[Y^\sigma\mid X]=\mathbb E[Y\mid X]=X$, we have $\pi^\sigma\in\Pi_M(\mu,\eta^\sigma)$. By independence and $\mathbb E[G]=0$,
\begin{equation}\label{eq:w2_gaussian_bound}
W_2^2(\nu,\eta^\sigma)\leq\mathbb E|V-Y-\sigma G|^2=\mathbb E|V-Y|^2+\sigma^2\mathbb E|G|^2=W_2^2(\nu,\eta)+d\sigma^2.
\end{equation}
The Lipschitz condition in Assumption~\ref{ass:c-lipschitz} gives
\begin{equation}\label{eq:cost_gaussian_bound}
\int c\,d\pi^\sigma-\int c\,d\pi=\mathbb E\bigl[c(X,Y+\sigma G)-c(X,Y)\bigr]\leq L_c\sigma\mathbb E|G|=C_dL_c\sigma.
\end{equation}
Moreover, $X\to Y\to Y^\sigma$ is a Markov chain. By the data processing inequality and Lemma~\ref{lem:gaussian_information},
\begin{equation}\label{eq:entropy_gaussian_bound}
H(\pi^\sigma\mid\mu\otimes\eta^\sigma)=I(X;Y^\sigma)\leq I(Y;Y^\sigma)\leq\frac d2\log\biggl(1+\frac{m_2(\eta)}{d\sigma^2}\biggr).
\end{equation}

Evaluating the regularized objective at $(\pi^\sigma,\eta^\sigma)$ and combining \eqref{eq:delta_optimal_pair}, \eqref{eq:w2_gaussian_bound}, \eqref{eq:cost_gaussian_bound}, and \eqref{eq:entropy_gaussian_bound}, we obtain
\begin{equation}\label{eq:before_delta_limit}
M^{\lambda,\epsilon}(\mu,\nu)-M^\lambda(\mu,\nu)\leq\delta+\lambda d\sigma^2+C_dL_c\sigma+\frac{\epsilon d}{2}\log\biggl(1+\frac{m_2(\eta)}{d\sigma^2}\biggr).
\end{equation}
Substituting $m_2(\eta)$ with the bound from \eqref{eq:eta_moment_bound} and letting $\delta\downarrow0$ yields \eqref{eq:entropy_approximation_bound}.

Finally, choosing $\sigma=\epsilon$ gives
\[
0\leq M^{\lambda,\epsilon}(\mu,\nu)-M^\lambda(\mu,\nu)\leq\lambda d\epsilon^2+C_dL_c\epsilon+\frac{\epsilon d}{2}\log\biggl(1+\frac{C_\lambda}{d\epsilon^2}\biggr).
\]
For fixed $\lambda>0$, the right-hand side is $O\bigl(\epsilon|\log\epsilon|\bigr)$ as $\epsilon\to0^+$ and therefore tends to zero. This proves both the convergence and the stated rate.
\end{proof}

We next show that letting $\lambda\to\infty$ recovers the second marginal constraint and yields convergence to the classical martingale optimal transport problem.

\begin{proposition}\label{prop:lambda_convergence}
Suppose Assumption~\ref{ass:c-lipschitz} holds. Let $\mu,\nu\in\mathcal P_2(\mathbb R^d)$ satisfy $\mu\preceq_{cx}\nu$. Then
\[
\lim_{\lambda\to\infty}M^\lambda(\mu,\nu)=M(\mu,\nu).
\]
If, in addition, Assumption~\ref{ass:dual-attainment} and ~\ref{ass:nu-compact} hold, then there exists a constant $C_{\mu,\nu,L_c}>0$, independent of $\lambda$, such that
\[
0\leq M(\mu,\nu)-M^\lambda(\mu,\nu)\leq\frac{C_{\mu,\nu,L_c}}{\lambda},\qquad \lambda>0.
\]
\end{proposition}

\begin{proof}
Taking $\eta=\nu$ in $M^{\lambda}$ gives $M^\lambda(\mu,\nu)\leq M(\mu,\nu)$. Moreover, $M^\lambda(\mu,\nu)$ is nondecreasing in $\lambda$ and bounded below by $-\|c\|_\infty$. Thus, its limit $\ell$ exists and satisfies $\ell\leq M(\mu,\nu)$.

Let $\lambda_n\to\infty$. For each $n$, choose $\eta_n\in\mathcal P_2(\mathbb R^d)$ and $\pi_n\in\Pi_M(\mu,\eta_n)$ such that
\begin{equation}\label{eq:lambda_approximate_minimizer}
\lambda_nW_2^2(\eta_n,\nu)+\int c\,d\pi_n\leq M^{\lambda_n}(\mu,\nu)+\frac1n.
\end{equation}
By Strassen's theorem, $\mu\preceq_{cx}\nu$ implies $\Pi_M(\mu,\nu)\neq\varnothing$, so $M(\mu,\nu)\leq\|c\|_\infty$. Using also $\int c\,d\pi_n\geq-\|c\|_\infty$, we obtain
\[
W_2^2(\eta_n,\nu)\leq\frac{2\|c\|_\infty+1/n}{\lambda_n}\to0.
\]
Hence $\eta_n\to\nu$ in $W_2$. Since the first marginal of $\pi_n$ is $\mu$ and its second marginal converges weakly to $\nu$, the sequence $(\pi_n)_n$ is tight. By Prokhorov's theorem, a subsequence converges weakly to some $\pi\in\Pi(\mu,\nu)$.

The convergence $\eta_n\to\nu$ in $W_2$ gives $\sup_n m_2(\eta_n)<\infty$. Since the first marginal of each $\pi_n$ is $\mu$, the function $q(x)\cdot(y-x)$ is uniformly integrable with respect to $(\pi_n)_n$ for every $q\in C_b(\mathbb R^d;\mathbb R^d)$. Therefore, uniform integrability and weak convergence of $(\pi_n)_n$ yield
\[
\int q(x)\cdot(y-x)\,\pi(dx,dy)=\lim_{n\to\infty}\int q(x)\cdot(y-x)\,\pi_n(dx,dy)=0.
\]
Thus $\pi\in\Pi_M(\mu,\nu)$. Since $c\in C_b(\mathbb{R}^{2d})$, \eqref{eq:lambda_approximate_minimizer} implies, along this subsequence,
\[
\ell=\lim_{n\to\infty}M^{\lambda_n}(\mu,\nu)\geq\lim_{n\to\infty}\int c\,d\pi_n=\int c\,d\pi\geq M(\mu,\nu).
\]
Together with $\ell\leq M(\mu,\nu)$, this proves the convergence.

We now prove the quantitative estimate under Assumptions~\ref{ass:dual-attainment} and~\ref{ass:nu-compact}. Set $K_\mu:=\operatorname{supp}\mu$, $K_\nu:=\operatorname{supp}\nu$, and  $m:=\int x\,\mu(dx)=\int y\,\nu(dy)$. We use \cite[Theorems~4.1 and~4.6(i)]{carlier2026weak}, applied to the MOT problem with their $\theta=0$, $f=0$, and $g(x,y)=y-x$. For compactly supported marginals and a continuous cost, these results give a bounded multiplier $h\in L^\infty(\mu;\mathbb R^d)$ satisfying
\begin{equation}\label{eq:mot_multiplier_representation}
M(\mu,\nu)=\inf_{\gamma\in\Pi(\mu,\nu)}\int\bigl[c(x,y)+h(x)\cdot(y-x)\bigr]\,\gamma(dx,dy),
\end{equation}
provided that $(\mu,\nu)$ satisfies the conditions in \cite[Definition~2.1]{carlier2026weak}. Since $\mu$ and $\nu$ have the same mean, the first condition holds automatically. The second requires the existence of $a>0$ such that, for every $u\in L^\infty(\mu;\mathbb R^d)$ with $\int u\,d\mu=0$ and $\|u\|_{L^\infty(\mu)}\leq a$, there is a coupling $\pi^u\in\Pi(\mu,\nu)$ whose disintegration satisfies $
\int(y-x)\,\pi^{u}_{x}(dy)=u(x)$ $\text{for $\mu$-a.e. }x.$

We first verify this condition. Recall that the proof of Proposition~\ref{prop:input-condition} constructs a reference martingale coupling between $\mu$ and the mean-corrected measure $\widehat\nu$. Here, $\mu\preceq_{cx}\nu$ implies $\widehat\nu=\nu$.  That construction therefore yields $\bar\pi=\bar r\,(\mu\otimes\nu)\in\Pi_M(\mu,\nu)$ with $\bar r\geq\kappa>0$ for some constant $\kappa$.
Define
\[
\Sigma_\nu:=\int(y-m)(y-m)^\top\,\nu(dy),\qquad R_\nu:=\sup_{y\in K_\nu}|y-m|.
\]
The matrix $\Sigma_\nu$ is positive definite because $\operatorname{aff}(K_\mu)=\mathbb R^d$ and $\mu\preceq_{cx}\nu$. Set $a:=\kappa/(2\|\Sigma_\nu^{-1}\|_{\mathrm{op}}R_\nu)$. Given $u\in L^\infty(\mu;\mathbb R^d)$ with $\int u\,d\mu=0$ and $\|u\|_{L^\infty(\mu)}\leq a$, define
\[
r_u(x,y):=\bar r(x,y)+u(x)\cdot\Sigma_\nu^{-1}(y-m),\qquad \pi_u:=r_u\,(\mu\otimes\nu).
\]
The choice of $a$ ensures that $r_u\geq\kappa/2$. Since $\int(y-m)\,\nu(dy)=0$ and $\int u\,d\mu=0$, the perturbation preserves both marginals, so $\pi_u\in\Pi(\mu,\nu)$. Moreover, the martingale property of $\bar\pi$ gives
\[
\int(y-x)r_u(x,y)\,\nu(dy)=\Sigma_\nu\Sigma_\nu^{-1}u(x)=u(x).
\]
This verifies the required condition. Therefore \eqref{eq:mot_multiplier_representation} holds for some bounded multiplier $h$. Furthermore, the a priori estimate \cite[equation~(51), in the proof of Theorem~4.1]{carlier2026weak} allows us to choose $h$ with $\|h\|_{L^\infty(\mu)}\leq H_{\mu,\nu,L_c}$, where $H_{\mu,\nu,L_c}<\infty$ is a constant depending only on $\mu$, $\nu$, and $L_c$, and is independent of $\lambda$.

Set $B_{\mu,\nu,L_c}:=L_c+H_{\mu,\nu,L_c}$. Let $\eta\in\mathcal P_2(\mathbb R^d)$ and $\pi\in\Pi_M(\mu,\eta)$. By the gluing lemma, choose $(X,Z,Y)$ such that $\mathcal L(X,Z)=\pi$ and $(Z,Y)$ is an optimal coupling of $\eta$ and $\nu$. Then $\mathcal L(X,Y)\in\Pi(\mu,\nu)$, while the martingale property gives $\mathbb E[h(X)\cdot(Z-X)]=0$. Using \eqref{eq:mot_multiplier_representation} and Assumption~\ref{ass:c-lipschitz}, we obtain
\[
M(\mu,\nu)\leq\mathbb E\bigl[c(X,Y)+h(X)\cdot(Y-Z)\bigr]\leq\int c\,d\pi+B_{\mu,\nu,L_c}W_2(\eta,\nu).
\]
Taking the infimum over $\pi\in\Pi_M(\mu,\eta)$ yields
\[
M(\mu,\nu)-M(\mu,\eta)\leq B_{\mu,\nu,L_c}W_2(\eta,\nu).
\]
Consequently, writing $r:=W_2(\eta,\nu)$, we have
\[
M(\mu,\eta)+\lambda W_2^2(\eta,\nu)\geq M(\mu,\nu)+\lambda r^2-B_{\mu,\nu,L_c}r\geq M(\mu,\nu)-\frac{B_{\mu,\nu,L_c}^2}{4\lambda}.
\]
Taking the infimum over $\eta$ and recalling that $M^\lambda(\mu,\nu)\leq M(\mu,\nu)$ gives
\[
0\leq M(\mu,\nu)-M^\lambda(\mu,\nu)\leq\frac{C_{\mu,\nu,L_c}}{\lambda},\qquad C_{\mu,\nu,L_c}:=\frac{(L_c+H_{\mu,\nu,L_c})^2}{4}.
\]
This completes the proof.
\end{proof}

Combining Propositions~\ref{prop:epsilon_convergence} and~\ref{prop:lambda_convergence}, we obtain the following explicit estimate for the convergence stated in Theorem~\ref{thm:convergence-to-mot}.

\begin{convergencerestatement}\label{thm:convergence-to-mot-sub}
Suppose Assumption~\ref{ass:c-lipschitz} holds. Let $\mu,\nu\in\mathcal P_2(\mathbb R^d)$ satisfy $\mu\preceq_{cx}\nu$. Then
\[
\lim_{\lambda\to\infty}\lim_{\epsilon\to0^+}M^{\lambda,\epsilon}(\mu,\nu)=M(\mu,\nu).
\]
Suppose additionally that Assumptions~\ref{ass:dual-attainment} and~\ref{ass:nu-compact} hold. Let $C_{\mu,\nu,L_c}$ be the constant in Proposition~\ref{prop:lambda_convergence}, and let $C_d:=\mathbb E|G|$ for $G\sim\mathcal N(0,I_d)$. Then, for every $\lambda,\epsilon>0$,
\begin{equation}\label{eq:two_parameter_rate}
\bigl|M^{\lambda,\epsilon}(\mu,\nu)-M(\mu,\nu)\bigr|
\leq\frac{C_{\mu,\nu,L_c}+d+C_dL_c}{\lambda}+\frac{d\epsilon}{2}\log\biggl(1+\frac{\lambda^2}{d}\bigl(\sqrt{m_2(\nu)}+\sqrt{2\|c\|_\infty/\lambda}\bigr)^2\biggr).
\end{equation}
\end{convergencerestatement}

\begin{proof}
By Proposition~\ref{prop:epsilon_convergence}, for every fixed $\lambda>0$, we have $M^{\lambda,\epsilon}(\mu,\nu)\to M^\lambda(\mu,\nu)$ as $\epsilon\to0^+$. Proposition~\ref{prop:lambda_convergence} then gives $M^\lambda(\mu,\nu)\to M(\mu,\nu)$ as $\lambda\to\infty$, proving the iterated limit.

Under the additional assumptions, combining the estimates in Proposition \ref{prop:epsilon_convergence} and \ref{prop:lambda_convergence} yields, for every $\sigma>0$,
\[
\bigl|M^{\lambda,\epsilon}(\mu,\nu)-M(\mu,\nu)\bigr|\leq\frac{C_{\mu,\nu,L_c}}{\lambda}+\lambda d\sigma^2+C_dL_c\sigma+\frac{d\epsilon}{2}\log\biggl(1+\frac{\bigl(\sqrt{m_2(\nu)}+\sqrt{2\|c\|_\infty/\lambda}\bigr)^2}{d\sigma^2}\biggr).
\]
We choose $\sigma=1/\lambda$ so that $\lambda d\sigma^2+C_dL_c\sigma=(d+C_dL_c)/\lambda$, which has the same order as the first term $C_{\mu,\nu,L_c}/\lambda$. Substituting this choice proves \eqref{eq:two_parameter_rate}. 

To recover the rate stated in Theorem~\ref{thm:convergence-to-mot}, set $D:=1+\bigl(m_2(\nu)+2\|c\|_\infty\bigr)/d$. Using $(a+b)^2\leq2a^2+2b^2$, we obtain
\[
1+\frac{\lambda^2}{d}\bigl(\sqrt{m_2(\nu)}+\sqrt{2\|c\|_\infty/\lambda}\bigr)^2\leq1+\frac{2m_2(\nu)}{d}\lambda^2+\frac{4\|c\|_\infty}{d}\lambda\leq(1+D\lambda)^2.
\]
Since $D\geq1$, Bernoulli's inequality gives $1+D\lambda\leq(1+\lambda)^D$. Thus the logarithmic term in \eqref{eq:two_parameter_rate} is bounded by $dD\epsilon\log(1+\lambda)$, and hence
\[
\bigl|M^{\lambda,\epsilon}(\mu,\nu)-M(\mu,\nu)\bigr|\leq C\biggl(\frac1\lambda+\epsilon\log(1+\lambda)\biggr),
\]
where
\[
C:=\max\bigl\{C_{\mu,\nu,L_c}+d+C_dL_c,\ d+m_2(\nu)+2\|c\|_\infty\bigr\}.
\]
Since $C$ is independent of $\lambda$ and $\epsilon$, this gives the claimed $O\bigl(\lambda^{-1}+\epsilon\log(1+\lambda)\bigr)$ bound. In particular, $M^{\lambda,\epsilon}(\mu,\nu)\to M(\mu,\nu)$ whenever $\lambda\to\infty$ and $\epsilon\log(1+\lambda)\to0$. This also completes the proof of Theorem~\ref{thm:convergence-to-mot}.
\end{proof}

\begin{remark}
\begin{enumerate}
\item[(i)]
The order of the two limits cannot in general be reversed. Consider $\mu=\nu=\mathcal N(0,1)$ and $c\equiv0$, so that $M(\mu,\mu)=0$. Fix $\epsilon>0$. For any $\eta\in\mathcal P_2(\mathbb R)$ and $\pi\in\Pi_M(\mu,\eta)$, let $(X,Y)\sim\pi$ and set $D:=\mathbb E|Y-X|^2$. The martingale property gives $\mathbb E[Y^2]=\mathbb E[X^2]+D=1+D$. Hence,
\[
W_2^2(\eta,\mu)\geq\bigl(\sqrt{m_2(\eta)}-\sqrt{m_2(\mu)}\bigr)^2=\bigl(\sqrt{1+D}-1\bigr)^2.
\]
Moreover, $H(\pi\mid\mu\otimes\eta)=I(X;Y)$. For $D>0$, the conditional representation of mutual information and the Gaussian density formula give
\[
\mathbb E\bigl[H\bigl(\mathcal L(X\mid Y)\mid\mathcal N(Y,D)\bigr)\bigr]=I(X;Y)+\frac12\log D-\frac12\mathbb E[X^2]+\frac{\mathbb E|X-Y|^2}{2D}=I(X;Y)+\frac12\log D.
\]
Since relative entropy is nonnegative, $I(X;Y)\geq\frac12\log(1/D)$.
Together with $I(X;Y)\geq0$, this implies $I(X;Y)\geq\frac12\log^+(1/D)$, where $\log^+t:=\max\{\log t,0\}$. If $D=0$, then $Y=X$ almost surely and $I(X;Y)=+\infty$, since $\mu$ is nonatomic. Consequently,
\[
M^{\lambda,\epsilon}(\mu,\mu)\geq\inf_{D>0}\biggl\{\lambda\bigl(\sqrt{1+D}-1\bigr)^2+\frac{\epsilon}{2}\log^+(1/D)\biggr\}.
\]
For any $a\in(0,1)$, splitting the infimum into $D\leq a$ and $D>a$ gives
\[
M^{\lambda,\epsilon}(\mu,\mu)\geq\min\biggl\{\frac{\epsilon}{2}\log(1/a),\,\lambda\bigl(\sqrt{1+a}-1\bigr)^2\biggr\}.
\]
Letting $\lambda\to\infty$ and then $a\downarrow0$ shows that $M^{\lambda,\epsilon}(\mu,\mu)\to+\infty$ for every fixed $\epsilon>0$. Therefore,
\[
\lim_{\epsilon\to0^+}\lim_{\lambda\to\infty}M^{\lambda,\epsilon}(\mu,\mu)=+\infty\neq M(\mu,\mu)=0,
\]
so the two limits cannot in general be exchanged.

\item[(ii)]
The sufficient condition in Theorem~\ref{thm:convergence-to-mot} allows a broad range of choices of $\lambda(\epsilon)\to\infty$. For example, if $\lambda(\epsilon)=\exp(\epsilon^{-\alpha})$ with $0<\alpha<1$, then
\[
\epsilon\log\bigl(1+\lambda(\epsilon)\bigr)=\epsilon^{1-\alpha}+\epsilon\log\bigl(1+\exp(-\epsilon^{-\alpha})\bigr)\to0.
\]
Thus, under Assumption \ref{ass:dual-attainment}-\ref{ass:nu-compact}, $M^{\lambda(\epsilon),\epsilon}(\mu,\nu)\to M(\mu,\nu)$ as $\epsilon\to0^+$. 
\end{enumerate}
\end{remark}


\section{Numerical Experiments}
\label{sec:numerics}

In this section, we conduct numerical experiments to illustrate three aspects of the WEMOT formulation: its stability under marginal perturbations, its convergence to the underlying MOT problem, and the computational benefits of entropy regularization. Unless otherwise stated, all experiments are carried out in dimension $d=2$ with the continuous bounded cost function
\[
c_\sigma(x,y)
=
1-\exp\left(-\frac{|x-y|^2}{2\sigma^2}\right),
\qquad \sigma=0.45.
\]

\subsection{A Sinkhorn-type solver for WEMOT}

Motivated by the dual formulation derived in
Section~\ref{sec:proof-strong-duality}, we solve the WEMOT problem using a
Sinkhorn-type alternating scheme. We describe the numerical procedure
below and refer to it as the \emph{WEMOT solver}. 

Let $
\mu=\sum_{i=1}^N a_i\delta_{x_i}$ and $
\nu=\sum_{j=1}^M b_j\delta_{y_j}$. 
For discrete marginals, the dual problem takes the finite-dimensional form
\[
\sup_{f,h}\inf_{z_1,\ldots,z_M}\mathcal D(f,h,z),
\]
where
\[
\mathcal D(f,h,z)
=
\sum_{i=1}^N a_i f_i
+
\sum_{j=1}^M b_j
\Bigg[
\lambda |z_j-y_j|^2
-
\epsilon
\log\bigg(
\sum_{i=1}^N
a_i
\exp\Big(
\frac{
f_i+h_i\cdot(z_j-x_i)-c(x_i,z_j)
}{\epsilon}
\Big)
\bigg)
\Bigg].
\]
Here, $f=(f_i)_{i=1}^N$ and $h=(h_i)_{i=1}^N$ are the dual variables
associated with the first-marginal and martingale constraints,
respectively, while $z_j\in\mathbb R^d$ denotes the contact point
associated with $y_j$ through the Wasserstein relaxation.

Starting from $(f^k,h^k,z^k)$, we first keep the dual variables
$(f^k,h^k)$ fixed and perform a fixed number $L$ of inner fixed-point iterations to update the contact points $(z_j)_{j=1}^{M}$. At the
$\ell$-th inner iteration, $\ell=0,\ldots,L-1$,
the entropy term induces the Gibbs weights
\[
p_{ij}^{\ell}
=
\frac{
a_i
\exp\Big(
\frac{
f_i+h_i\cdot(z_j-x_i)-c(x_i,z_j)
}{\epsilon}
\Big)
}{
\sum_{r=1}^N
a_r
\exp\Big(
\frac{
f_r+h_r\cdot(z_j-x_r)-c(x_r,z_j)
}{\epsilon}
\Big)
}, \qquad
i=1,\ldots,N,\quad j=1,\ldots,M.
\]
After updating the contact points, we recompute the Gibbs weights and
use the resulting first-marginal and martingale constraint to update
$f$ and $h$, respectively. The procedure is repeated until both
residuals fall below a prescribed numerical tolerance. The complete
iteration is summarized in Algorithm~\ref{alg:WEMOT}.

\par\medskip
\begingroup
\small

\refstepcounter{algorithm}
\label{alg:WEMOT}
\addcontentsline{loa}{algorithm}
  {\protect\numberline{\thealgorithm}Sinkhorn-Type Solver for WEMOT}

\hrule height 0.8pt
\nobreak\vspace{2pt}

\noindent
\textbf{Algorithm \thealgorithm}
Sinkhorn-Type Solver for WEMOT
\par\nobreak\vspace{2pt}

\hrule height 0.4pt
\nobreak\vspace{2pt}

\begin{algorithmic}[1]

\Require
$\mu=\sum_{i=1}^N a_i\delta_{x_i}$,
$\nu=\sum_{j=1}^M b_j\delta_{y_j}$,
$\lambda,\epsilon,\tau>0$,
$\rho\in(0,1)$,
and $L\in\mathbb N$

\State Initialize $f_i^0=0$, $h_i^0=0$, and $z_j^0=y_j$

\For{$k=0,1,2,\ldots$}

    \State Set $z_j^{k,0}=z_j^k$ for $j=1,\ldots,M$

    \For{$\ell=0,\ldots,L-1$}

    \Statex \hspace{\algorithmicindent}
    $\displaystyle
    p_{ij}^{k,\ell}
    =
    \frac{
    a_i
    \exp\!\left(
    \frac{
    f_i^k+h_i^k\cdot(z_j^{k,\ell}-x_i)
    -c(x_i,z_j^{k,\ell})
    }{\epsilon}
    \right)
    }{
    \sum_{r=1}^N
    a_r
    \exp\!\left(
    \frac{
    f_r^k+h_r^k\cdot(z_j^{k,\ell}-x_r)
    -c(x_r,z_j^{k,\ell})
    }{\epsilon}
    \right)
    }.
    $

    \Statex \hspace{\algorithmicindent}
    $\displaystyle
    z_j^{k,\ell+1}
    =
    (1-\rho)z_j^{k,\ell}
    +
    \rho\left[
    y_j
    +
    \frac{1}{2\lambda}
    \sum_{i=1}^N
    p_{ij}^{k,\ell}
    \bigl(
    h_i^k-\nabla_z c(x_i,z_j^{k,\ell})
    \bigr)
    \right].
    $

    \EndFor

    \State Set $z_j^{k+1}=z_j^{k,L}$

    \Statex \hspace{\algorithmicindent}
    $\displaystyle
    p_{ij}^{k,+}
    =
    \frac{
    a_i
    \exp\!\left(
    \frac{
    f_i^k+h_i^k\cdot(z_j^{k+1}-x_i)
    -c(x_i,z_j^{k+1})
    }{\epsilon}
    \right)
    }{
    \sum_{r=1}^N
    a_r
    \exp\!\left(
    \frac{
    f_r^k+h_r^k\cdot(z_j^{k+1}-x_r)
    -c(x_r,z_j^{k+1})
    }{\epsilon}
    \right)
    }.
    $

    \State Compute
    $\displaystyle
    r_i^{k,+}
    =
    \sum_{j=1}^M b_jp_{ij}^{k,+}
    $
    and
    $\displaystyle
    m_i^{k,+}
    =
    \sum_{j=1}^M
    b_jp_{ij}^{k,+}(z_j^{k+1}-x_i)
    $

    \If{$\max_i |r_i^{k,+}-a_i|<\mathrm{tol}$ \textbf{and}
         $\max_i \|m_i^{k,+}\|<\mathrm{tol}$}
        \State \textbf{break}
    \EndIf

    \Statex \hspace{\algorithmicindent}
    $\displaystyle
    f_i^{k+1}
    =
    f_i^k
    +
    \epsilon\log\!\left(\frac{a_i}{r_i^{k,+}}\right),
    \qquad
    h_i^{k+1}
    =
    h_i^k
    -
    \tau\frac{m_i^{k,+}}{a_i}.
    $

\EndFor

\State \Return
$\displaystyle
\pi_{ij}=b_jp_{ij}^{k,+},
\quad
\eta=\sum_{j=1}^M b_j\delta_{z_j^{k+1}}
$

\end{algorithmic}
\par\nobreak
\hrule height 0.8pt

\endgroup
\medskip
This solver is used
throughout the numerical experiments in this section. We use two inner contact-point
updates per outer iteration, i.e., $L=2$, with damping parameter
$\rho=1/2$. The martingale dual variable is updated with step size
$\tau=0.03$. We terminate the iteration when both the first-marginal
residual and the martingale residual satisfy
\[
\max_i |r^{k,+}_i-a_i| < 3\times 10^{-6},
\qquad
\max_i \|m_i^{k,+}\| < 3\times 10^{-6},
\]
that is, we set $\mathrm{tol}=3\times 10^{-6}$.
These parameters are kept fixed across all experiments unless stated
otherwise.

\subsection{Stability under marginal perturbations}
\label{subsec:stability-comparison}

\begin{figure}[H]
    \centering
    \includegraphics[width=0.95\textwidth]{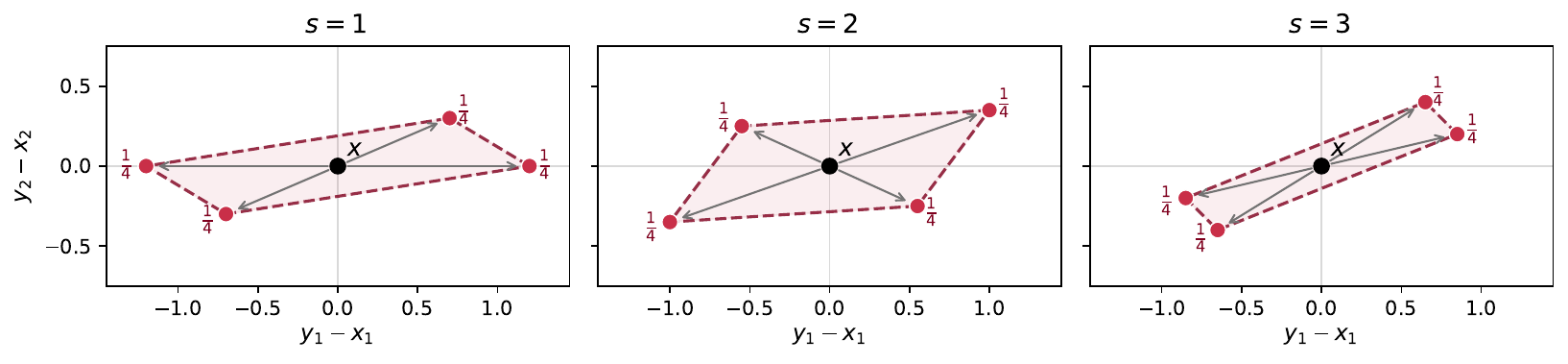}
    \caption{Geometry of $Y^{(s)}$ given $X=x$. The four red
vertices are $x\pm A_s e_1$ and $x\pm A_s e_2$, each with conditional
probability $1/4$.}
    \label{fig:conditional-geometry}
\end{figure}

We first compare the stability of classical MOT and WEMOT to
perturbations of marginals. We work in dimension $d=2$ and
consider the following family of models. Let $
X\sim\mathcal N(0,I_2)$
and let $Z$ be independent of $X$ and uniformly distributed on $
\{\pm e_1,\pm e_2\}$, where $e_i$ are coordinate vectors. 
For each case $s\in\{1,2,3\}$, define
\[
Y^{(s)}=X+A_s Z,
\]
where
\[
A_1=
\begin{pmatrix}
1.20 & 0.70\\
0 & 0.30
\end{pmatrix},
\qquad
A_2=
\begin{pmatrix}
1.00 & -0.55\\
0.35 & 0.25
\end{pmatrix},
\qquad
A_3=
\begin{pmatrix}
0.85 & 0.65\\
0.20 & 0.40
\end{pmatrix}.
\]
Here, $s$ indexes three different choices of the transition matrix
$A_s$, corresponding to different martingale transition geometries; see
Figure~\ref{fig:conditional-geometry} for an illustration .
We consider several $A_s$ to verify that the observed
stability behavior is not specific to a single configuration. Since $\mathbb E[Z]=0$, we have $
\mathbb E[Y^{(s)}\mid X]=X,$
so that, denoting
\[
\mu=\operatorname{Law}(X),
\qquad
\nu^{(s)}=\operatorname{Law}(Y^{(s)}),
\]
we have $\mu\preceq_{cx}\nu^{(s)}$ for each $s$.

To generate empirical approximation of the marginals while preserving the convex order, for each sample size $N$ we draw
independent samples $
X_1,\ldots,X_N\sim\mathcal N(0,I_2)$
and define 
\begin{align}\label{eq:stability_experiment_setting}
\mu_N
=
\frac1N\sum_{i=1}^N\delta_{X_i},
\qquad
\nu_N^{(s)}
=
\frac1{4N}
\sum_{i=1}^N
\sum_{z\in\{\pm e_1,\pm e_2\}}
\delta_{X_i+A_s z}.
\end{align}
For every realization and every $N$, the coupling that sends each
$X_i$ uniformly to the four points
$X_i+A_s z$ is a martingale coupling. Hence $\mu_N\preceq_{cx}\nu_N^{(s)}$
holds exactly. Thus the classical MOT problem remains feasible for every realization. 

For each test case $s\in\{1,2,3\}$ and each sample size $
N\in\{10,20,40,80,160\},$
we generate $R=20$ independent realizations of the empirical pair,
denoted by $
(\mu_N^{(r)},\nu_N^{(s,r)}),$ $ r=1,\ldots,R.$
For each realization, we compute both the classical MOT value $
M\bigl(\mu_N^{(r)},\nu_N^{(s,r)}\bigr)$
and the WEMOT value $
M^{\lambda,\epsilon}
\bigl(\mu_N^{(r)},\nu_N^{(s,r)}\bigr),$
where we fix $
\lambda=2,$ $
\epsilon=0.1.$ For each pair $(s,N)$, we compute the sample mean and sample standard
deviation of the values over the $R$ realizations. In
particular, for WEMOT we consider
\begin{align}\label{eq:stability_std}
\widehat{\operatorname{Std}}_{\mathrm{WEMOT}}(s,N)
=
\left[
\frac{1}{R-1}
\sum_{r=1}^R
\left(
M^{\lambda, \epsilon}(\mu_N^{(r)},\nu_N^{(s,r)})
-
\overline M_{\mathrm{WEMOT}}(s,N)
\right)^2
\right]^{1/2},
\end{align}
and analogously define
$\widehat{\operatorname{Std}}_{\mathrm{MOT}}(s,N)$ for classical MOT.
A smaller standard deviation indicates that the corresponding optimal
value is more stable to the random perturbations of the empirical
marginals.

\begin{figure}[H]
    \centering

    \begin{subfigure}[t]{0.32\textwidth}
        \centering
        \includegraphics[width=\textwidth]{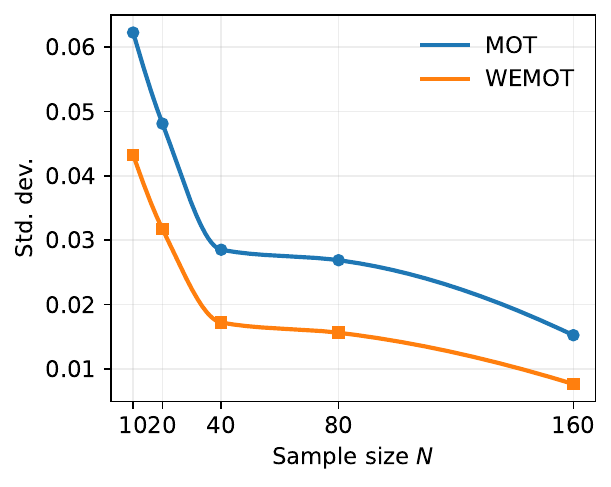}
        \caption{$s=1$}
        \label{fig:stability-s1}
    \end{subfigure}
    \hfill
    \begin{subfigure}[t]{0.32\textwidth}
        \centering
        \includegraphics[width=\textwidth]{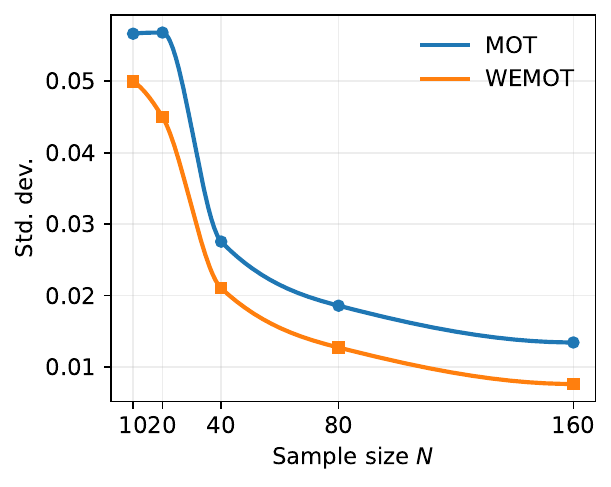}
        \caption{$s=2$}
        \label{fig:stability-s2}
    \end{subfigure}
    \hfill
    \begin{subfigure}[t]{0.32\textwidth}
        \centering
        \includegraphics[width=\textwidth]{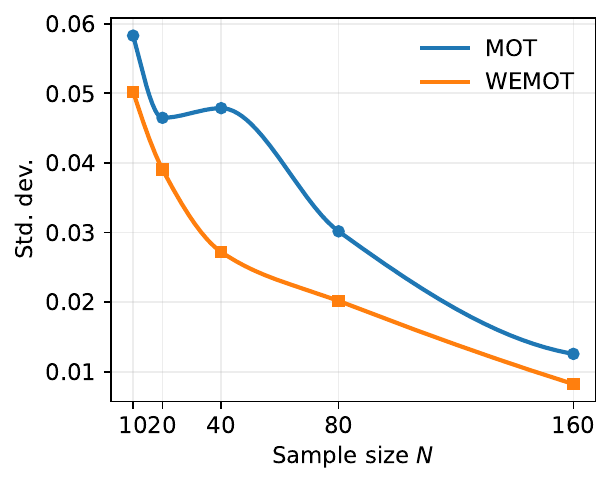}
        \caption{$s=3$}
        \label{fig:stability-s3}
    \end{subfigure}

    \caption{Standard deviation of the MOT and WEMOT values across
    independent realizations of the empirical marginals for the three
    cases.}
    \label{fig:mot-WEMOT-stability}
\end{figure}

Figure~\ref{fig:mot-WEMOT-stability} reports the resulting standard
deviations as a function of the sample size $N$, separately for the
three test cases. Across all three test cases and all sample sizes considered, the
standard deviation of the WEMOT value is smaller than that of the
classical MOT value. Moreover, the variability of both formulations
decreases as $N$ increases, while WEMOT consistently exhibits a lower
level of sensitivity to the empirical marginals. 

\subsection{Quantitative stability rate}\label{subsec:stability_rate}
Our quantitative stability result, Theorem~\ref{thm:quantitative-stability},
suggests that perturbations of the WEMOT value, $|
M^{\lambda,\epsilon}(\mu_N,\nu_N)
-
M^{\lambda,\epsilon}(\mu,\nu)|,$
are controlled by the Wasserstein scale
\[
\Delta_N
:=
W_2(\mu_N,\mu)^{2/3}
+
W_2(\nu_N,\nu).
\]
To illustrate this dependence numerically, we use the same Gaussian
setting as in Section~\ref{subsec:stability-comparison} with $s=1$,
and fix $\lambda=2$ and $\epsilon=0.1$. We consider
$N\in\{10,20,40,80,160,320\}$. For each $N$, we generate
$R=100$ independent realizations, denoted by
$(\mu_N^{(r)},\nu_N^{(r)})_{r=1}^R$, and compute the sample standard deviation
of the corresponding WEMOT values as in
Equation~\eqref{eq:stability_std}.

Since the reference marginals $\mu$ and $\nu$ are continuous, the
Wasserstein distances appearing in $\Delta_N$ are approximated
numerically. Specifically, we construct Monte Carlo reference measures
$\widetilde{\mu}$ and $\widetilde{\nu}$ using $8000$ samples from the
respective distributions, and approximate
$W_2(\mu_N,\mu)$ and $W_2(\nu_N,\nu)$ by
$W_2(\mu_N,\widetilde{\mu})$ and
$W_2(\nu_N,\widetilde{\nu})$, respectively. These discrete
$W_2$ distances are computed by solving the quadratic optimal
transport problem with squared Euclidean cost. For each $N$, we average
the resulting stability scale over the same $R=100$ repetitions and
denote it by $\overline{\Delta}_N$.

\begin{figure}[htbp]
    \centering
    \includegraphics[width=0.68\textwidth]{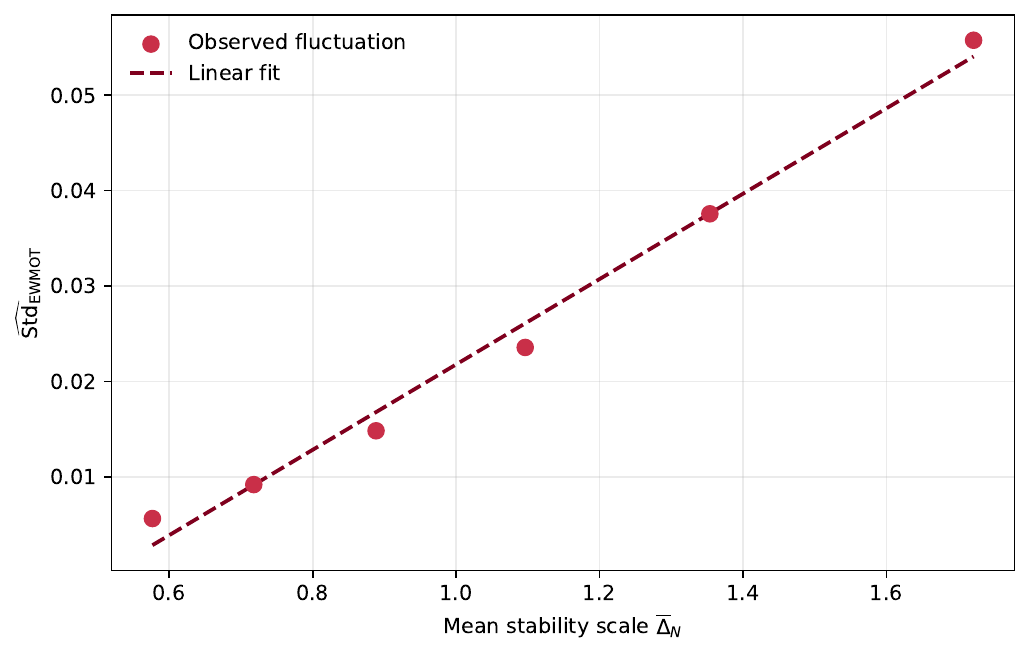}
    \caption{Empirical fluctuation of the WEMOT value against the mean
    stability scale $\overline{\Delta}_N$. The points represent the sample standard deviation
    $\widehat{\mathrm{Std}}_{\mathrm{WEMOT}}$ computed from $R=100$
    independent realizations for each sample size $N$, while the dashed line
    shows a linear least-squares fit.}
    \label{fig:stability-rate}
\end{figure}

Figure~\ref{fig:stability-rate} plots the empirical standard deviation
$\widehat{\operatorname{Std}}_{\mathrm{WEMOT}}$ against $\overline{\Delta}_N$. The observed fluctuations exhibit an approximately linear dependence on
$\Delta_N$, as indicated by the dashed linear fit. This is consistent
with the quantitative stability estimate, which predicts that the
variation of the WEMOT value is controlled at the scale $\Delta_N$.

\subsection{Stability of WEMOT on a MOT counterexample}
\label{subsec:bruckerhoff-instability}
Classical MOT is known to be unstable under marginal perturbations in
dimension $d\geq 2$. To illustrate how WEMOT can mitigate this
instability, we consider a modification of the counterexample in
\cite{bruckerhoff2022instability}.

Fix $m=5$ and let $
x_i=(i,0)$ for $ i=1,\ldots,m.$ Define $
\mu=\frac1m\sum_{i=1}^m\delta_{x_i}.$
For $n\geq 2$, let
\[
\theta_n=\frac{\pi}{2n},
\qquad
u_n=(\cos\theta_n,\sin\theta_n),
\]
and define the perturbed second marginal by $
\nu_n
=
\frac1{2m}
\sum_{i=1}^m
\left(
\delta_{x_i+u_n}
+
\delta_{x_i-u_n}
\right).$
Equivalently, $\nu_n$ is obtained by sending each atom $x_i$ with
equal probability to $x_i+u_n$ and $x_i-u_n$. This  construction defines a martingale coupling between $\mu$ and
$\nu_n$, and hence $
\mu\preceq_{cx}\nu_n$
for every $n$.

As $n\to\infty$, we have $\theta_n\to0$ and $u_n\to e_1$, and hence $
\nu_n $ converges to $\nu
:=
\frac{1}{2m}\sum_{i=1}^m
(
\delta_{x_i+e_1}+\delta_{x_i-e_1})$
in $W_2$. Following the instability analysis in
\cite{bruckerhoff2022instability}, classical MOT is not stable along
this sequence of marginals. Indeed, for the bounded cost $c_\sigma$ considered
here,
\[
M(\mu,\nu_n)
=
1-\exp\left(-\frac{1}{2\sigma^2}\right)
\approx 0.915342
\qquad\text{for every } n,
\]
whereas $
M(\mu,\nu)\approx 0.185890.$
Consequently, $
\lim_{n\to\infty} M(\mu,\nu_n)
\neq
M(\mu,\nu),$
despite $\nu_n\to\nu$ in $W_2$.

We next apply WEMOT to the same sequence of marginals. We consider
\[
n\in\{2,3,5,10,20,50,100,200\}
\]
and choose
\[
\epsilon_n=\frac{0.15}{\sqrt n},
\qquad
\lambda_n^{(\alpha)}
=
\exp\left(
\bigl(\log(1/\epsilon_n)\bigr)^\alpha
\right),
\qquad
\alpha\in\{0.25,0.5,0.75\}.
\]
For each fixed $\alpha>0$, we have
\[
\epsilon_n\to0,
\qquad
\lambda_n^{(\alpha)}\to\infty,
\qquad
\epsilon_n\log\lambda_n^{(\alpha)}
=
\epsilon_n
\bigl(\log(1/\epsilon_n)\bigr)^\alpha
\to0.
\]
Hence, all three parameter schedules satisfy the sufficient condition
in our joint convergence result. For each $\alpha$, we compute $
M^{\lambda_n^{(\alpha)},\epsilon_n}(\mu,\nu_n)$
and compare the resulting values with both the classical MOT value
$M(\mu,\nu_n)$ and the limiting MOT value $M(\mu,\nu)$.

\begin{figure}[H]
    \centering
    \includegraphics[width=0.7\textwidth]{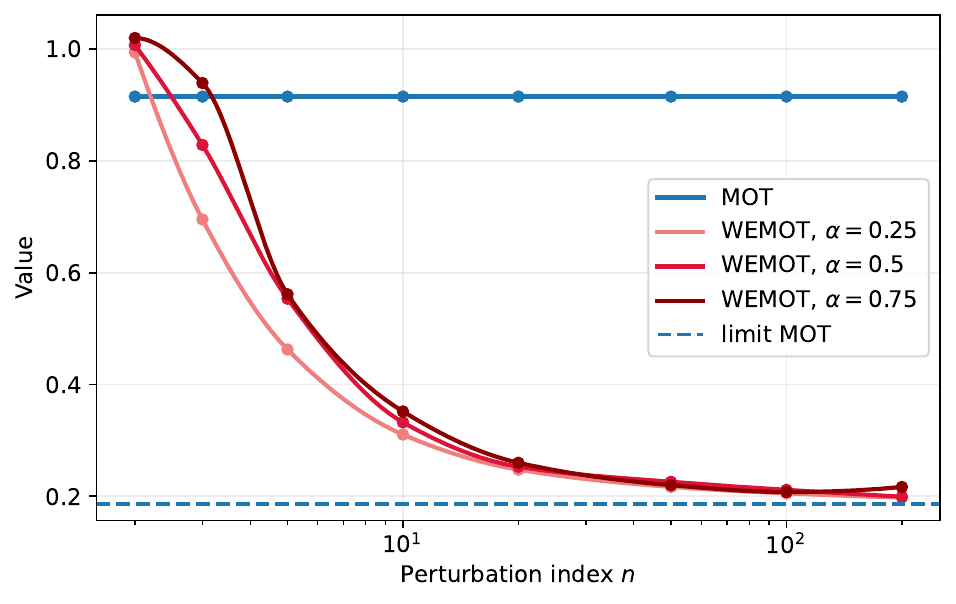}
    \caption{MOT and WEMOT values for the Br\"uckerhoff--Juillet-type perturbation.}
    \label{fig:Bruckerhoff_joint_stability_WEMOT}
\end{figure}

Figure~\ref{fig:Bruckerhoff_joint_stability_WEMOT} shows the resulting values.
While the classical MOT $M(\mu, \nu_n)$ value remains separated from
$M(\mu,\nu)$ along the perturbed marginal sequence, the WEMOT values associated
with all $\alpha$ move toward the limiting MOT value as $n$
increases. This illustrates both the stabilizing effect of the
Wasserstein relaxation on this classical MOT counterexample and the
joint convergence behavior predicted by our theoretical results.

\subsection{Computational effect of entropy regularization}
\label{subsec:computational-convergence}

We investigate the computational effect of entropy
regularization by comparing the numerical convergence of the
Wasserstein-relaxed MOT problem
$M^\lambda(\mu,\nu)$
with WEMOT
$M^{\lambda,\epsilon}(\mu,\nu)$.
The purpose of this subsection is to investigate whether the
Sinkhorn-type scheme enabled by entropy regularization leads to
faster numerical convergence compared with the unregularized
Wasserstein-relaxed formulation.

We fix the limiting marginal pair $(\mu,\nu)$ from
Section~\ref{subsec:bruckerhoff-instability}. More precisely, with
$m=5$ and $x_i=(i,0)$, $i=1,\ldots,5$, we consider
\[
\mu=\frac{1}{5}\sum_{i=1}^5\delta_{x_i},
\qquad
\nu=\frac{1}{10}\sum_{i=1}^5
\left(
\delta_{x_i+e_1}+\delta_{x_i-e_1}
\right).
\]
Since both marginals are finitely supported, no additional
discretization is required, and the numerical algorithms are applied
directly to these discrete measures.
We fix the Wasserstein penalty parameter at $\lambda=2$ and compare
the numerical convergence of
$M^\lambda(\mu,\nu)$
with its entropy-regularized counterpart
$M^{\lambda,\epsilon}(\mu,\nu)$.
For the latter, we consider
$\epsilon\in\{0.2,0.15,0.1\}$.
Thus, the marginals and $\lambda$ are kept fixed throughout the
experiment, while the entropy regularization parameter $\epsilon$ is
varied.

\begin{figure}[htbp]
    \centering
    \begin{subfigure}[t]{0.49\textwidth}
        \centering
        \includegraphics[width=\linewidth]{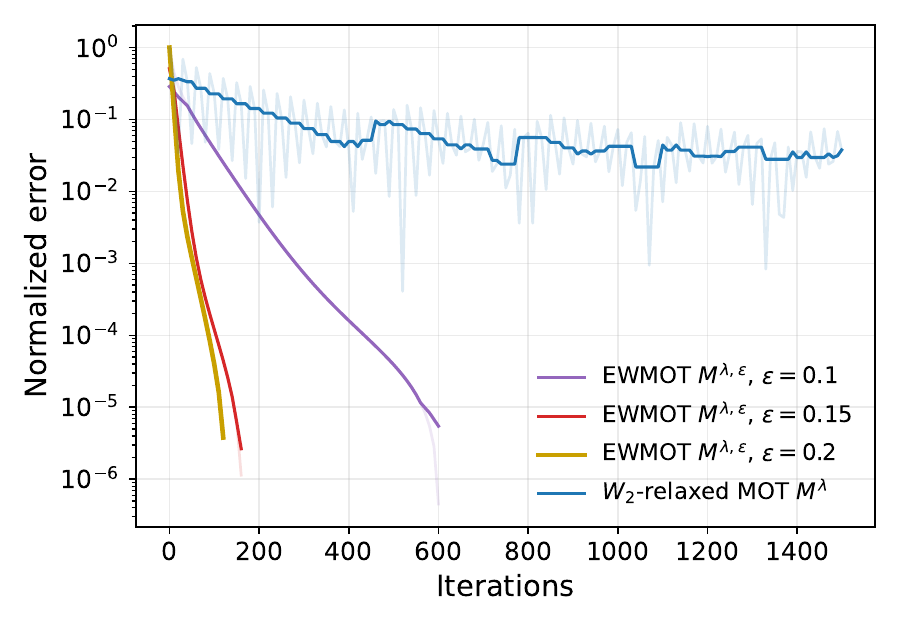}
        \caption{Convergence versus iteration.}
        \label{fig:exp5-iteration}
    \end{subfigure}
    \hfill
    \begin{subfigure}[t]{0.49\textwidth}
        \centering
        \includegraphics[width=\linewidth]{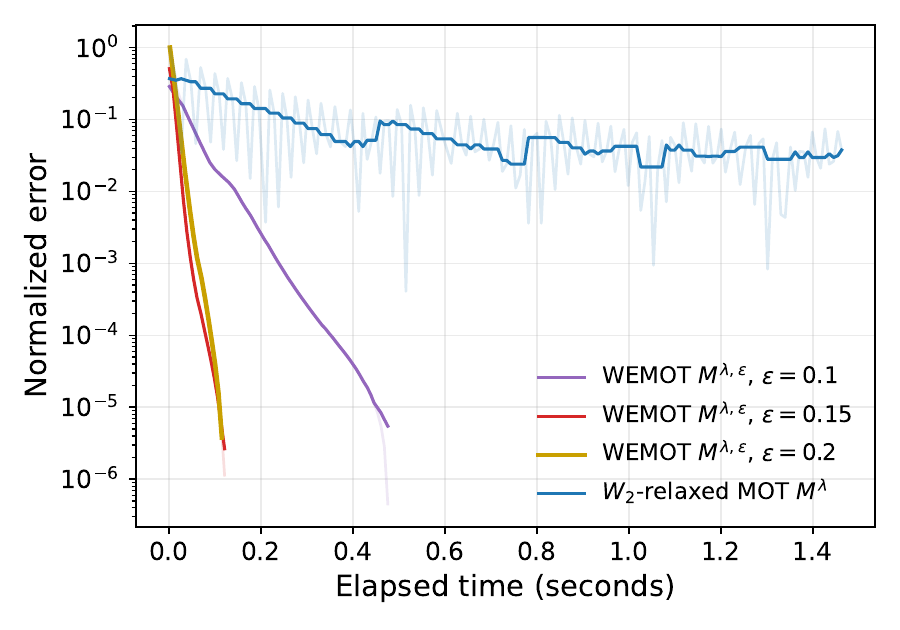}
        \caption{Convergence versus runtime.}
        \label{fig:exp5-time}
    \end{subfigure}
    \caption{Convergence speed of WEMOT $M^{\lambda,\epsilon}$ and
    $W_2$-relaxed MOT $M^\lambda$.
    Faint curves show the raw trajectories, while solid curves
    show their smoothed versions.}
    \label{fig:computational-convergence}
\end{figure}

For each formulation, we record the objective values
$M_k^\lambda$ and $M_k^{\lambda,\epsilon}$ at iteration $k$ and
compare them with the corresponding reference values
$M_{\mathrm{ref}}^\lambda$ and
$M_{\mathrm{ref}}^{\lambda,\epsilon}$, obtained from substantially
longer runs of the same formulations. We report the normalized numerical errors
\[
\operatorname{Err}^\lambda(k)
:=
\frac{|M_k^\lambda-M_{\mathrm{ref}}^\lambda|}{|M_0^\lambda-M_{\mathrm{ref}}^\lambda|},
\qquad
\operatorname{Err}^{\lambda,\epsilon}(k)
:=
\frac{|M_k^{\lambda,\epsilon}
-M_{\mathrm{ref}}^{\lambda,\epsilon}|}{|M_0^{\lambda,\epsilon}
-M_{\mathrm{ref}}^{\lambda,\epsilon}|}.
\]
Figure~\ref{fig:computational-convergence} compares these numerical
errors against both the iteration number and the elapsed
computational time.

Figure~\ref{fig:computational-convergence} shows a clear computational
advantage of WEMOT. In both panels,
WEMOT reduces the normalized error substantially faster than the
$W_2$-relaxed MOT problem $M^\lambda$. In particular, for
$\epsilon=0.15$ and $\epsilon=0.2$, the error decreases by several
orders of magnitude within relatively few iterations and a small amount
of computational time, whereas the numerical error for $M^\lambda$
decreases much more slowly and remains at a visibly higher level over
the same computational budget. The case $\epsilon=0.1$ exhibits an
intermediate behavior, converging more slowly than the strongly
regularized cases but still considerably faster than $M^\lambda$.
These results illustrate the computational benefit of the
Sinkhorn-type structure induced by entropy regularization, as well as
the trade-off that smaller values of $\epsilon$ generally lead to
slower numerical convergence.

\bibliographystyle{plainnat}
\bibliography{refs}

\appendix
\section*{Appendix}
\section{Dual attainment for EMOT}
\label{app:emot-dual-attainment}

Our proof for Proposition \ref{prop:emot-dual-attainment} is based on a Galerkin approximation of the martingale constraint. We choose an increasing sequence
$(\mathcal H_n)_{n\geq1}$ of finite-dimensional subspaces of
$C(K;\mathbb R^d)$ and replace the martingale constraint, tested
against all functions in $C(K;\mathbb R^d)$, by its Galerkin
approximation, tested only against functions in $\mathcal H_n$. We first prove that these truncated constraints are locally nondegenerate, that is,  every sufficiently small  violation of the martingale constraint can be realized by a coupling with uniformly bounded relative entropy. This local surjectivity helps us to obtain a Lagrange multiplier $h_n$ for each Galerkin problem, which serves as an approximation of the martingale potential. Together with the marginal potentials $f_n$ and $g_n$, these multipliers yield Gibbs representations of the Galerkin optimizers. We then show that the associated Gibbs densities converge strongly to the optimal density, while an additive compactness argument yields convergence of the potentials, after suitable affine gauge transformations. Passing to the limit produces an admissible dual optimizer for the EMOT problem. To the best of our knowledge, this Galerkin approach to dual attainment for EMOT has not previously been used in the literature.

Given $\mu,\eta\in\mathcal{P}_2(\mathbb{R}^d)$ with
$\mu\preceq_{cx}\eta$, write $K:=\operatorname{supp}(\mu)$
and $\rho:=\mu\otimes\eta$. For $\pi\in\Pi(\mu,\eta)$, define
\begin{equation}
\mathcal{E}(\pi):=\int_{\mathbb{R}^d\times\mathbb{R}^d}
c(x,z)\,\pi(dx,dz)+\epsilon H(\pi\mid\rho).
\end{equation}
We write $V:=\operatorname{EMOT}(\mu,\eta)
=\inf_{\pi\in\Pi_M(\mu,\eta)}\mathcal{E}(\pi)$.
We begin by establishing several lemmas used in the Galerkin
approximation argument. The first shows that full-dimensionality
of the support is preserved under convex order.

\begin{lemma}
\label{lem:support-convex-order}
Let $\mu,\eta\in\mathcal P_2(\mathbb R^d)$ satisfy
$\mu\preceq_{\mathrm{cx}}\eta$. If
$\operatorname{aff}(\operatorname{supp}(\mu))=\mathbb R^d$, then
\[
\operatorname{aff}(\operatorname{supp}(\eta))
=
\mathbb R^d.
\]
\end{lemma}

\begin{proof}
Let
$C:=\overline{\operatorname{co}}(\operatorname{supp}(\eta))$.
Since $C$ is closed and convex, the function
$x\mapsto\operatorname{dist}(x,C)$ is convex, nonnegative, and
$1$-Lipschitz, and is therefore integrable under both $\mu$ and
$\eta$. By $\mu\preceq_{\mathrm{cx}}\eta$,
\[
0
\leq
\int_{\mathbb R^d}\operatorname{dist}(x,C)\,\mu(dx)
\leq
\int_{\mathbb R^d}\operatorname{dist}(z,C)\,\eta(dz)
=
0.
\]
Hence $\mu(C)=1$. Since $C$ is closed, this implies
$\operatorname{supp}(\mu)\subseteq C$.

Moreover,
$C\subseteq\operatorname{aff}(\operatorname{supp}(\eta))$, because
the latter is a closed affine subspace containing
$\operatorname{supp}(\eta)$, while
$\operatorname{supp}(\eta)\subseteq C$. Therefore
\[
\operatorname{aff}(C)
=
\operatorname{aff}(\operatorname{supp}(\eta)).
\]
It follows that
\[
\mathbb R^d
=
\operatorname{aff}(\operatorname{supp}(\mu))
\subseteq
\operatorname{aff}(C)
=
\operatorname{aff}(\operatorname{supp}(\eta)),
\]
and hence
$\operatorname{aff}(\operatorname{supp}(\eta))=\mathbb R^d$.
\end{proof}

The next lemma shows that the truncated martingale constraints can be
perturbed locally in every direction while retaining a uniformly
bounded relative entropy. This local surjectivity will ensure the
existence of Lagrange multipliers for the finite-dimensional Galerkin
problems.

\begin{lemma}
\label{lem:local_surjectivity}
Let $\mathcal H\subset C(K;\mathbb R^d)$ be finite-dimensional and
contain the constant vector fields, and set $\mathcal H_0
:=\{
q\in\mathcal H:
\int_K q\,d\mu=0
\}.$ Let $q_1,\ldots,q_N$ be a basis of $\mathcal H_0$, and define
\[
T_{\mathcal H}(\pi)
:=
\left(
\int q_i(x)\cdot(z-x)\,\pi(dx,dz)
\right)_{i=1}^N.
\]
Under the assumptions of
Proposition~\ref{prop:emot-dual-attainment}, there exist an open
neighborhood $U$ of $0$ in $\mathbb R^N$ and a constant
$C_0<\infty$ such that, for every $\zeta\in U$, there exists
$\pi^\zeta\in\Pi(\mu,\eta)$ satisfying
\[
T_{\mathcal H}(\pi^\zeta)=\zeta,
\qquad H(\pi^\zeta\,|\,\mu\otimes\eta)\leq C_0.
\]
\end{lemma}

\begin{proof}
Let $m_\eta:=\int_{\mathbb R^d}z\,\eta(dz)$ and define
\[
b_0(z):=\frac{z-m_\eta}{1+|z-m_\eta|},
\qquad
b(z):=b_0(z)-\int_{\mathbb R^d}b_0\,d\eta.
\]
Then $b\in L^\infty(\eta;\mathbb R^d)$ and $\int b\,d\eta=0$. Consider the matrix $B:=\int_{\mathbb R^d}(z-m_\eta)b(z)^\top\,\eta(dz)$.
 $\int(z-m_\eta)\,\eta(dz)=0$ yields that
\[
B
=
\int_{\mathbb R^d}
\frac{(z-m_\eta)(z-m_\eta)^\top}{1+|z-m_\eta|}
\,\eta(dz).
\]

We now show that $B$ is positive definite. 
For every $a\in\mathbb R^d$,
\[
a^\top Ba
=
\int_{\mathbb R^d}
\frac{\bigl[a\cdot(z-m_\eta)\bigr]^2}
{1+|z-m_\eta|}
\,\eta(dz)\geq 0.
\]
If $a^\top Ba=0$, then
$a\cdot(z-m_\eta)=0$ $\eta$-a.e. Since the map
$z\mapsto a\cdot(z-m_\eta)$ is continuous, it vanishes on
$\operatorname{supp}(\eta)$. Hence, if $a\neq0$,
$\operatorname{supp}(\eta)$ would be contained in the proper affine
hyperplane
$\{z:a\cdot(z-m_\eta)=0\}$, and therefore
$\operatorname{aff}(\operatorname{supp}(\eta))\neq\mathbb R^d$.
On the other hand, Lemma~\ref{lem:support-convex-order}  implies that
$\operatorname{aff}(\operatorname{supp}(\eta))=\mathbb R^d$, leading a
contradiction. Therefore $a^\top Ba>0$ for every $a\neq0$, so $B$ is
positive definite and hence invertible.

We next construct $\pi^\zeta$ for arbitrary sufficiently
small values $\zeta$ such that $T_{\mathcal{H}}(\pi^\zeta)=\zeta$. Let
\[
G_{ij}:=\int_K q_i(x)\cdot q_j(x)\,\mu(dx),
\qquad
1\leq i,j\leq N.
\]
The Gram matrix $G$ is positive definite.  For $j=1,\ldots,N$, define
\[
a_j(x):=
\sum_{k=1}^N(G^{-1})_{jk}q_k(x),
\qquad
r_\zeta(x):=
\sum_{j=1}^N\zeta_j a_j(x),
\quad
\zeta\in\mathbb R^N.
\]
Since $q_k\in\mathcal H_0$, one has $\int a_j\,d\mu=0$, and by the
definition of $G$,
$\int q_i(x)\cdot a_j(x)\,\mu(dx)=\delta_{ij}$.
Consequently,
$\int r_\zeta\,d\mu=0$ and
$\int q_i(x)\cdot r_\zeta(x)\,\mu(dx)=\zeta_i$.
Moreover, each $a_j$ is continuous on the compact set $K$, so
$\|r_\zeta\|_\infty\to0$ as $\zeta\to0$. For $\zeta$ sufficiently close to zero, define
\[
r^\zeta(x,z)
:=
\bar r(x,z)
+
b(z)^\top B^{-1}r_\zeta(x).
\]
Since $b$ and $B^{-1}$ are bounded, there exists an open
neighborhood $U$ of $0$ such that
$|b(z)^\top B^{-1}r_\zeta(x)|\leq\kappa/2$ for every
$\zeta\in U$, $(\mu\otimes\eta)$-a.e. Hence
\[
\frac{\kappa}{2}
\leq
r^\zeta(x,z)
\leq
\kappa^{-1}+\frac{\kappa}{2}.
\]
In particular, $r^\zeta$ is nonnegative. Moreover, using
$\int b\,d\eta=0$ and $\int r_\zeta\,d\mu=0$, we obtain
$\eta(r^\zeta)=1$ for $\mu$-a.e. $x$ and
$\int r^\zeta(x,z)\,\mu(dx)=1$ for $\eta$-a.e. $z$. Therefore
$\pi^\zeta:=r^\zeta(\mu\otimes\eta)$ belongs to $\Pi(\mu,\eta)$.

It remains to compute its conditional drift. Since $\bar\pi$ is a
martingale,
$\int(z-x)\bar r(x,z)\,\eta(dz)=0$ for $\mu$-a.e. $x$. Moreover,
using $\int b\,d\eta=0$,
\[
\int_{\mathbb R^d}(z-x)b(z)^\top\,\eta(dz)
=
\int_{\mathbb R^d}(z-m_\eta)b(z)^\top\,\eta(dz)
=
B.
\]
Hence
\[
\int_{\mathbb R^d}(z-x)r^\zeta(x,z)\,\eta(dz)
=
BB^{-1}r_\zeta(x)
=
r_\zeta(x).
\]
It follows that, for every $i=1,\ldots,N$,
\[
[T_{\mathcal H}(\pi^\zeta)]_i
=
\int_K q_i(x)\cdot r_\zeta(x)\,\mu(dx)
=
\zeta_i,
\]
and therefore $T_{\mathcal H}(\pi^\zeta)=\zeta$.

Finally, the preceding uniform upper bound gives
$r^\zeta\leq C$ for some constant $C<\infty$ independent of
$\zeta\in U$. Since $\int r^\zeta\,d(\mu\otimes\eta)=1$,
\[ H(\pi^\zeta\,|\,\mu\otimes\eta)
=
\int r^\zeta\log r^\zeta\,d(\mu\otimes\eta)
\leq
\log C.
\]
Thus the entropy is uniformly bounded for $\zeta\in U$, which
completes the proof.
\end{proof}

The next lemma shows that convergence of the combined potentials
implies almost-everywhere convergence of the individual potentials,
after suitable affine normalization.

\begin{lemma}
\label{lem:additive_compactness}
Let $f_n,g_n:\mathbb{R}^d\to\mathbb{R}$ and
$h_n:\mathbb{R}^d\to\mathbb{R}^d$ be measurable, and set
$v_n(x,z):=f_n(x)+g_n(z)+h_n(x)\cdot(z-x)$.
Suppose that $v_n\to v$ in $\rho$-probability, where $v$ is
finite $\rho$-a.e., and that
$\operatorname{aff}(\operatorname{supp}(\eta))=\mathbb{R}^d$.

After passing to a subsequence and applying affine gauge
transformations
\[
(f_n(x),g_n(z),h_n(x))
\mapsto
(f_n(x)-\beta_n-a_n\cdot x,\,
 g_n(z)+\beta_n+a_n\cdot z,\,
 h_n(x)-a_n),
\]
with $\beta_n\in\mathbb{R}$ and $a_n\in\mathbb{R}^d$,
the relabeled potentials satisfy
$f_n\to f$ and $h_n\to h$ $\mu$-a.e., and
$g_n\to g$ $\eta$-a.e., for finite measurable functions $f,g,h$.
Moreover, $v(x,z)=f(x)+g(z)+h(x)\cdot(z-x)$ $\rho$-a.e.,
and this representation is unique up to an affine gauge
transformation.
\end{lemma}

\begin{proof}
Passing to a subsequence, assume that $v_n\to v$ $(\mu\otimes\eta)$-a.e. Let $E$ be a $(\mu\otimes\eta)$-full-measure set such that, for every $(x,z)\in E$ and every $n$, one has $v_n(x,z)=f_n(x)+g_n(z)+h_n(x)\cdot(z-x)$ and $v_n(x,z)\to v(x,z)$. 

For $z\in\mathbb R^d$, write $E^z:=\{x:(x,z)\in E\}$. By Fubini's theorem, the set
$Z:=\{z:\mu(E^z)=1\}$ satisfies $\eta(Z)=1$. Since $Z$ has full
$\eta$-measure, it is dense in $\operatorname{supp}(\eta)$. Hence
$\operatorname{aff}(Z)=\operatorname{aff}(\operatorname{supp}(\eta))
=\mathbb R^d$, and we may choose affinely independent points
$z_0,\ldots,z_d\in Z$. Similarly, for $x\in\mathbb R^d$, write
$E_x:=\{z:(x,z)\in E\}$. Again by Fubini's theorem, the set
$X:=\{x:\eta(E_x)=1\}$ has full $\mu$-measure. Since
$z_i\in Z$, one has $\mu(E^{z_i})=1$ for every
$i=0,\ldots,d$. Therefore
$X\cap E^{z_0}\cap\cdots\cap E^{z_d}$ has full $\mu$-measure.
Choose $x_0$ in this intersection. Then
$\eta(E_{x_0})=1$ and $(x_0,z_i)\in E$ for every
$i=0,\ldots,d$.

For each $n$, we use the affine gauge transform to normalize
$\widetilde h_n(x_0)=0$ and $\widetilde g_n(z_0)=0$. More precisely,
we may take $a_n:=h_n(x_0)$ and
$\beta_n:=-g_n(z_0)-a_n\cdot z_0$. Since the affine transformation
does not change $v_n$, we then have
$\widetilde f_n(x_0)=v_n(x_0,z_0)$. In particular,
$\widetilde f_n(x_0)$ converges. Moreover, for
$i=1,\ldots,d$,
\[
\widetilde g_n(z_i)
=
v_n(x_0,z_i)-\widetilde f_n(x_0),
\]
so each sequence $\widetilde g_n(z_i)$ also converges.

Let $A$ be the $d\times d$ matrix whose $i$-th row is
$(z_i-z_0)^\top$. Since $z_0,\ldots,z_d$ are affinely independent,
$A$ is invertible. Set
$X_1:=E^{z_0}\cap\cdots\cap E^{z_d}$, which has full $\mu$-measure
because each $z_i\in Z$. For every $x\in X_1$, evaluating the representation of $v_n$ at
$(x,z_i)$ and $(x,z_0)$ and subtracting the two identities gives
\[
A\widetilde h_n(x)
=
\begin{pmatrix}
v_n(x,z_1)-v_n(x,z_0)-\widetilde g_n(z_1)\\
\vdots\\
v_n(x,z_d)-v_n(x,z_0)-\widetilde g_n(z_d)
\end{pmatrix},
\]
where we used $\widetilde g_n(z_0)=0$. The right-hand side converges
for every $x\in X_1$. Since $A$ is invertible,
$\widetilde h_n(x)$ converges for every $x\in X_1$; denote its limit
by $h(x)$.

For $x\in X_1$, evaluating the representation of $v_n$ at $(x,z_0)$
and using $\widetilde g_n(z_0)=0$ yields
$\widetilde f_n(x)
=
v_n(x,z_0)-\widetilde h_n(x)\cdot(z_0-x)$, and hence
$\widetilde f_n(x)$ converges; denote its limit by $f(x)$. On the
other hand, since $\eta(E_{x_0})=1$, for every $z\in E_{x_0}$ we
have
$\widetilde g_n(z)
=
v_n(x_0,z)-\widetilde f_n(x_0)$, so
$\widetilde g_n(z)$ converges for $\eta$-a.e. $z$; denote its limit
by $g(z)$. Extending $f,g,h$ arbitrarily on the corresponding null
sets, we obtain finite measurable functions satisfying the claimed
pointwise convergence.

Finally, on the $(\mu\otimes\eta)$-full-measure set on which the
representation of $v_n$ and all three convergences hold, we may pass
to the limit in
$v_n(x,z)
=
\widetilde f_n(x)+\widetilde g_n(z)
+\widetilde h_n(x)\cdot(z-x)$.
Therefore
\[
v(x,z)
=
f(x)+g(z)+h(x)\cdot(z-x)
\qquad
(\mu\otimes\eta)\text{-a.e.}
\]

It remains to prove uniqueness. Suppose that
$(\bar f,\bar g,\bar h)$ is another representation of $v$, and set
$A(x):=\bar f(x)-f(x)-[\bar h(x)-h(x)]\cdot x$,
$H(x):=\bar h(x)-h(x)$, and $G(z):=\bar g(z)-g(z)$. Then
\begin{align}\label{eq:uniquess_v}
A(x)+H(x)\cdot z+G(z)=0
\qquad
(\mu\otimes\eta)\text{-a.e.}
\end{align}
By Fubini's theorem, there exists a set $X_2$ of full $\mu$-measure
such that, for every $x\in X_2$, \eqref{eq:uniquess_v} holds for
$\eta$-a.e. $z$. Fix $x'\in X_2$. Then for any $x\in X_2$, evaluating \eqref{eq:uniquess_v} at $(x,z)$ and $(x', z)$, and subtracting the corresponding
equalities give
\[
A(x)-A(x')
+
\bigl[H(x)-H(x')\bigr]\cdot z
=
0
\qquad
\eta\text{-a.e.}
\]
The left-hand side is an affine function of $z$. By continuity, it
therefore vanishes on $\operatorname{supp}(\eta)$. Since
$\operatorname{aff}(\operatorname{supp}(\eta))=\mathbb R^d$, an
affine function vanishing on $\operatorname{supp}(\eta)$ must be
identically zero. Hence $A(x)=A(x'):=\beta$ and $H(x)=H(x'):=a$ for $\mu-$a.e. $x$.
 Substituting this into
$A(x)+H(x)\cdot z+G(z)=0$ gives
$G(z)=-\beta-a\cdot z$ for $\eta$-a.e. $z$. By the definitions of
$A,H,G$, this is equivalent to
\[
\bar f(x)=f(x)+\beta+a\cdot x,\qquad
\bar h(x)=h(x)+a
\quad \mu\text{-a.e.}, \qquad
\bar g(z)=g(z)-\beta-a\cdot z
\quad \eta\text{-a.e.},
\]
which proves uniqueness up to the affine gauge transform.
\end{proof}

The next result shows that the limiting Gibbs representation
yields an admissible dual optimizer.

\begin{lemma}
\label{lem:limiting-potential-admissibility}
Assume the hypotheses of
Proposition~\ref{prop:emot-dual-attainment}, and let
$\pi_0=r_0\rho$ be a primal optimizer with $r_0>0$ $\rho$-a.e.
For each $n$, let $f_n\in L^1(\mu)$ and
$h_n:\mathbb{R}^d\to\mathbb{R}^d$ be measurable, and suppose
$g_n:=\Gamma^\epsilon_{f_n,h_n}\in L^1(\eta)$. Define
\[
r_n(x,z):=\exp\biggl(
\frac{f_n(x)+g_n(z)+h_n(x)\cdot(z-x)-c(x,z)}{\epsilon}
\biggr).
\]
Suppose that $\pi^n:=r_n\rho\in\Pi(\mu,\eta)$,
$r_n\to r_0$ in $L^1(\rho)$, and
$\mu(f_n)+\eta(g_n)\to V$.

After passing to a subsequence and applying the affine
normalizations of Lemma~\ref{lem:additive_compactness},
the relabeled potentials satisfy
$f_n\to f$ and $h_n\to h$ $\mu$-a.e., and
$g_n\to g$ $\eta$-a.e., where
$f\in L^1(\mu)$, $g\in L^1(\eta)$, and
$h:\mathbb{R}^d\to\mathbb{R}^d$ is finite and measurable.
Moreover, $g=\Gamma^{\epsilon}_{f,h}$ $\eta$-a.e.,
$(f,h)\in\mathcal D_{\operatorname{EMOT}}(\mu,\eta)$,
and $V=\mu(f)+\eta(g)$.
\end{lemma}

\begin{proof}
\noindent\emph{Step 1.}
Since $r_n\to r_0$ in $L^1(\rho)$, after passing to a subsequence
we may assume that $r_n\to r_0$ $\rho$-a.e. As $r_0>0$
$\rho$-a.e., the functions $v_n:=\epsilon\log r_n+c$ converge
$\rho$-a.e. to the finite function $v:=\epsilon\log r_0+c$.
By Lemma~\ref{lem:additive_compactness}, after applying affine
normalizations and relabeling, we have $f_n\to f$ and $h_n\to h$
$\mu$-a.e., and $g_n\to g$ $\eta$-a.e., for finite measurable
functions $f,g,h$, with
\[
\epsilon\log r_0(x,z)
=f(x)+g(z)+h(x)\cdot(z-x)-c(x,z)
\qquad \rho\text{-a.e.}
\]
These normalizations preserve $r_n$ and
$g_n=\Gamma^\epsilon_{f_n,h_n}$. They also preserve
$\mu(f_n)+\eta(g_n)$, since $\mu$ and $\eta$ have the same
barycenter.

Choose $x_0$ in a common full $\mu$-measure set on which
$f_n(x_0)$ and $h_n(x_0)$ converge and that
$\int r_n(x_0,z)\,\eta(dz)=1$ for every $n$.
Applying a further affine normalization with
$a_n=h_n(x_0)$ and $\beta_n=0$, whose parameters converge,
we may assume that $h_n(x_0)=0$ for every $n$.
The pointwise convergence of the potentials is preserved,
and $f_n(x_0)$ still converges.

\medskip
\noindent\emph{Step 2.}
Let $C_\eta:=\overline{\operatorname{co}}(\operatorname{supp}(\eta))$
and let $m$ denote the common barycenter of $\mu$ and $\eta$.
By Lemma~\ref{lem:support-convex-order},
$\operatorname{aff}(\operatorname{supp}(\eta))=\mathbb R^d$,
and hence $m\in\operatorname{int}C_\eta$.
Indeed, otherwise a supporting hyperplane of $C_\eta$ at $m$
would contain $\operatorname{supp}(\eta)$, contradicting its
full affine dimension. Choose $r>0$ such that
$B(m,r)\subset C_\eta$.

Let $\bar\pi=\bar r\rho$ and $\kappa\in(0,1)$ be as in
Proposition~\ref{prop:emot-dual-attainment}.
For $\mu$-a.e. $x$, set
$s_x(z):=(\bar r(x,z)-\kappa)/(1-\kappa)$.
Then $s_x\geq0$, $\int s_x\,d\eta=1$, and the martingale
property gives $x=\kappa m+(1-\kappa)y_x$, where
$y_x:=\int z s_x(z)\,\eta(dz)\in C_\eta$.
For $|u|<\kappa r$, convexity therefore yields
\[
x+u=\kappa\bigl(m+u/\kappa\bigr)+(1-\kappa)y_x\in C_\eta.
\]
Thus $B(x,\kappa r)\subset C_\eta$ for $\mu$-a.e. $x$.
By density of this full-measure set in $K$ and closedness
of $C_\eta$, the inclusion holds for every $x\in K$.

Choose a full $\eta$-measure set
$Z\subset\operatorname{supp}(\eta)$ on which $g_n\to g$ and
\[
e^{-g_n(z)/\epsilon}
=\int_{\mathbb R^d}
\exp\biggl(
\frac{f_n(x)+h_n(x)\cdot(z-x)-c(x,z)}{\epsilon}
\biggr)\mu(dx)
\]
for every $n$. Then
$\overline{\operatorname{co}}(Z)=C_\eta$, and
$\operatorname{int}C_\eta=\operatorname{int}\operatorname{co}(Z)$
by \cite[Theorem~6.3]{rockafellar1970convex}.
Fix $\delta\in(0,\kappa r)$.
Since $K+\delta\overline B(0,1)$ is a compact subset of
$\operatorname{int}C_\eta$, there exist $z_1,\ldots,z_N\in Z$
such that
$K+\delta\overline B(0,1)\subset\operatorname{co}\{z_1,\ldots,z_N\}$.
To see this, each point of the compact set lies in the
interior of the convex hull of finitely many points of $Z$;
a finite covering argument then gives the claim.

Consequently,
$\delta|p|\leq\max_{1\leq j\leq N}p\cdot(z_j-x)$
for every $x\in K$ and $p\in\mathbb R^d$.
Taking $p=h_n(x)$, we obtain some $j$ for which
$h_n(x)\cdot(z_j-x)\geq0$. Hence
\[
e^{f_n^+(x)/\epsilon}
\leq 1+e^{\|c\|_\infty/\epsilon}
\sum_{j=1}^N
\exp\biggl(
\frac{f_n(x)+h_n(x)\cdot(z_j-x)-c(x,z_j)}{\epsilon}
\biggr).
\]
Integrating and using the identity defining $Z$ gives
\[
\int e^{f_n^+/\epsilon}\,d\mu
\leq 1+e^{\|c\|_\infty/\epsilon}
\sum_{j=1}^N e^{-g_n(z_j)/\epsilon}.
\]
The right-hand side is uniformly bounded because
$g_n(z_j)\to g(z_j)$ for every $j$.

The first marginal condition at $x_0$, together with
$h_n(x_0)=0$, gives
\[
e^{-f_n(x_0)/\epsilon}
=\int_{\mathbb R^d}
\exp\biggl(\frac{g_n(z)-c(x_0,z)}{\epsilon}\biggr)\eta(dz).
\]
It follows that
$\int e^{g_n^+/\epsilon}\,d\eta
\leq 1+e^{\|c\|_\infty/\epsilon}e^{-f_n(x_0)/\epsilon}$,
which is also uniformly bounded.
The de la Vallée--Poussin criterion and Vitali's theorem
therefore imply that
$\mu(f_n^+)\to\mu(f^+)<\infty$ and
$\eta(g_n^+)\to\eta(g^+)<\infty$.

Write $P_n:=\mu(f_n^+)+\eta(g_n^+)$ and
$N_n:=\mu(f_n^-)+\eta(g_n^-)$.
Since $P_n\to P:=\mu(f^+)+\eta(g^+)<\infty$ and
$P_n-N_n\to V$, we have $N_n\to P-V$.
Fatou's lemma gives
$\mu(f^-)+\eta(g^-)\leq\liminf_n N_n=P-V<\infty$.
Thus $f\in L^1(\mu)$ and $g\in L^1(\eta)$.

Finally, the second marginal condition of $\pi_0$ and
the limiting Gibbs identity yield
\[
e^{-g(z)/\epsilon}
=\int_{\mathbb R^d}
\exp\biggl(
\frac{f(x)+h(x)\cdot(z-x)-c(x,z)}{\epsilon}
\biggr)\mu(dx)
\qquad \eta\text{-a.e.}
\]
Therefore $g=\Gamma^{\epsilon}_{f,h}$ $\eta$-a.e.

\medskip
\noindent\emph{Step 3.}
Let $\pi\in\Pi_M(\mu,\eta)$ satisfy
$H(\pi\mid\rho)<\infty$, and set $\ell:=\log r_0$.
Since $e^{\ell^+}\leq1+r_0$, we have
$\int e^{\ell^+}\,d\rho\leq2$.
Applying the entropy variational inequality to
$\ell^+\wedge M$ and then letting $M\to\infty$ gives
$\int\ell^+\,d\pi\leq H(\pi\mid\rho)+\log2<\infty$.
As $\pi\ll\rho$, the Gibbs identity holds $\pi$-a.e. and implies
\[
[h(x)\cdot(z-x)]^+
\leq\epsilon\ell^+(x,z)+|f(x)|+|g(z)|+\|c\|_\infty.
\]
The right-hand side is $\pi$-integrable because $\pi$ has
marginals $\mu$ and $\eta$. Thus
$[h(x)\cdot(z-x)]^+\in L^1(\pi)$.

Disintegrate $\pi(dx,dz)=\mu(dx)\pi_x(dz)$.
For $\mu$-a.e. $x$, the conditional first moment is finite,
and the martingale property gives
\[
\int h(x)\cdot(z-x)\,\pi_x(dz)
=h(x)\cdot\biggl(\int z\,\pi_x(dz)-x\biggr)=0.
\]
The positive and negative conditional integrals therefore
coincide. By Tonelli's theorem,
$\int[h(x)\cdot(z-x)]^-\,d\pi
=\int[h(x)\cdot(z-x)]^+\,d\pi<\infty$.
Consequently, $h(x)\cdot(z-x)\in L^1(\pi)$ and
$\int h(x)\cdot(z-x)\,d\pi=0$.
Since $\pi$ was arbitrary,
$(f,h)\in\mathcal D_{\operatorname{EMOT}}(\mu,\eta)$.

The bounded-density competitor $\bar\pi$ has finite entropy,
so the primal optimizer $\pi_0$ also satisfies
$H(\pi_0\mid\rho)<\infty$.
Applying the preceding conclusion to $\pi_0$ and integrating
the Gibbs identity gives
\[
V=\int c\,d\pi_0+\epsilon H(\pi_0\mid\rho)
=\mu(f)+\eta(g).
\]
Together with $g=\Gamma^{\epsilon}_{f,h}$ $\eta$-a.e., this proves
the claim.
\end{proof}

\begin{proof}[Proof of Proposition~\ref{prop:emot-dual-attainment}]
\noindent\emph{Step 1.}
Since $K$ is compact, $C(K;\mathbb R^d)$ is separable under
the uniform norm. Choose increasing finite-dimensional subspaces
$\mathcal H_n\subset C(K;\mathbb R^d)$ containing the constant
vector fields such that
$\overline{\bigcup_{n\geq1}\mathcal H_n}^{\|\cdot\|_\infty}
=C(K;\mathbb R^d)$.
Define
\[
\mathcal M_n:=
\biggl\{\pi\in\Pi(\mu,\eta):
\int q(x)\cdot(z-x)\,d\pi=0
\text{ for every }q\in\mathcal H_n\biggr\},
\]
and set $P_n:=\inf_{\pi\in\mathcal M_n}\mathcal E(\pi)$.

The set $\Pi(\mu,\eta)$ is weakly compact. For every
$q\in C(K;\mathbb R^d)$, the map
$\pi\mapsto\int q(x)\cdot(z-x)\,d\pi$ is continuous on
$\Pi(\mu,\eta)$ under weak convergence: $q$ is bounded,
$x$ lies in the compact set $K$, and the tails in $z$ are
uniformly controlled by the fixed marginal $\eta$.
Thus $\mathcal M_n$ is weakly compact.
The coupling $\bar\pi=\bar r\rho$ from the hypotheses belongs
to every $\mathcal M_n$ and has finite entropy, since
$\bar r\leq\kappa^{-1}$.
Lower semicontinuity of relative entropy and bounded continuity
of $c$ therefore yield a minimizer $\pi^n$ with $P_n<\infty$.
Strict convexity of relative entropy gives uniqueness.
The same argument yields the unique primal optimizer $\pi_0$,
with $\mathcal E(\pi_0)=V<\infty$.
In particular, $\pi^n=r_n\rho$ and $\pi_0=r_0\rho$ for
densities $r_n$ and $r_0$ of finite entropy.

Since $\mathcal M_{n+1}\subseteq\mathcal M_n$ and
$\Pi_M(\mu,\eta)\subseteq\mathcal M_n$, we have
$P_n\uparrow L\leq V$.
Let $\pi^{n_k}\Rightarrow\widetilde\pi$ along a subsequence.
For every $q\in\bigcup_n\mathcal H_n$, the corresponding
constraint holds for $\pi^{n_k}$ for all sufficiently large
$k$, and hence also for $\widetilde\pi$.
By uniform density, it holds for every
$q\in C(K;\mathbb R^d)$, so
$\widetilde\pi\in\Pi_M(\mu,\eta)$.
Lower semicontinuity gives
$V\leq\mathcal E(\widetilde\pi)\leq L$.
Consequently, $L=V$ and $\widetilde\pi=\pi_0$.
Since every subsequential limit is $\pi_0$, we obtain
$P_n\to V$ and $\pi^n\Rightarrow\pi_0$.

\medskip
\noindent\emph{Step 2.}
Fix $n$ and set
$\mathcal H_n^0:=\{q\in\mathcal H_n:\int q\,d\mu=0\}$.
Let $q_1,\ldots,q_{N_n}$ be a basis of $\mathcal H_n^0$,
and define
\[
T_n(\pi):=
\biggl(\int q_i(x)\cdot(z-x)\,d\pi\biggr)_{i=1}^{N_n},
\qquad
\Psi_n(\zeta):=
\inf_{\substack{\pi\in\Pi(\mu,\eta)\\T_n(\pi)=\zeta}}
\mathcal E(\pi).
\]
The constant vector fields impose no additional constraints,
because $\mu$ and $\eta$ have the same barycenter.
Thus $P_n=\Psi_n(0)$.

By Lemma~\ref{lem:local_surjectivity}, $\Psi_n$ is finite
in a neighborhood of $0$.
Since $\Psi_n$ is convex and bounded below by
$-\|c\|_\infty$, it has a subgradient at $0$.
Choose $\ell_n\in\partial\Psi_n(0)$ and set
$h_n:=\sum_{i=1}^{N_n}(\ell_n)_i q_i$.
The subgradient inequality yields
$\mathcal E(\pi)\geq P_n+\ell_n\cdot T_n(\pi)$
for every $\pi\in\Pi(\mu,\eta)$.
Since
$\ell_n\cdot T_n(\pi)=\int h_n(x)\cdot(z-x)\,d\pi$,
with equality in the preceding inequality at $\pi^n$,
we have
\begin{equation}
P_n=
\inf_{\pi\in\Pi(\mu,\eta)}
\biggl\{\mathcal E(\pi)
-\int h_n(x)\cdot(z-x)\,d\pi\biggr\}.
\label{eq:galerkin_lagrangian}
\end{equation}

Hence $\pi^n$ solves the entropic optimal transport problem
with effective cost $c_n(x,z):=c(x,z)-h_n(x)\cdot(z-x)$
on $K\times\mathbb R^d$.
Write $M_n:=\|h_n\|_\infty$ and
$R:=\max_{x\in K}|x|$, and set
$b_n(z):=\|c\|_\infty+M_n(R+|z|)$.
Then $|c_n(x,z)|\leq b_n(z)$ and $b_n\in L^1(\eta)$.
In particular, $c_n+b_n$ is nonnegative and belongs to
$L^1(\rho)$.
Applying the entropic optimal transport dual attainment result \cite[Propositions~A.1 and~A.2]{nutz2023stability} to the cost $(c_n+b_n)/\epsilon$, and then subtracting $b_n$ from
the rescaled second marginal potential, yields
$f_n\in L^1(\mu)$ and $g_n\in L^1(\eta)$ such that
\[
r_n(x,z)=
\exp\biggl(
\frac{f_n(x)+g_n(z)+h_n(x)\cdot(z-x)-c(x,z)}
{\epsilon}
\biggr)
\qquad \rho\text{-a.e.},
\]
and $P_n=\mu(f_n)+\eta(g_n)$.
The second marginal condition gives
$g_n=\Gamma^\epsilon_{f_n, h_n}$ $\eta-$a.e., and we choose this version of $g_n$.

\medskip
\noindent\emph{Step 3.}
Lemma~\ref{lem:fixed-marginal-entropy-gap} and its proof remain valid
with $\eta^*$ replaced by $\eta$.
Since $h_n$ is bounded, $f_n\in L^1(\mu)$, and
$g_n=\Gamma^\epsilon_{f_n,h_n}\in L^1(\eta)$, we have
$(f_n,h_n)\in\mathcal D_{\operatorname{EMOT}}(\mu,\eta)$.
The lemma therefore gives
\[
\epsilon H(\pi_0\mid\pi^n)
=V-\mu(f_n)-\eta(g_n)
=V-P_n\to0.
\]
By Pinsker's inequality,
$r_n\to r_0$ in $L^1(\rho)$.

To prove strict positivity, let $\bar\pi=\bar r\rho$
be the martingale coupling from the hypotheses.
Since $r_n>0$ $\rho$-a.e., we have $\bar\pi\ll\pi^n$ and
$d\bar\pi/d\pi^n=\bar r/r_n$.
Moreover, $f_n\in L^1(\mu)$, $g_n\in L^1(\eta)$,
$h_n$ is bounded, and $\bar\pi$ has finite first moments.
Thus $\log r_n\in L^1(\bar\pi)$, and
\[
H(\bar\pi\mid\pi^n)
=\int\log\frac{\bar r}{r_n}\,d\bar\pi
=H(\bar\pi\mid\rho)-\int\log r_n\,d\bar\pi.
\]
Since $\bar\pi\in\Pi_M(\mu,\eta)$,
$\int h_n(x)\cdot(z-x)\,d\bar\pi=0$.
The Gibbs identity therefore yields
\[
\epsilon\int\log r_n\,d\bar\pi
=\mu(f_n)+\eta(g_n)-\int c\,d\bar\pi
=P_n-\int c\,d\bar\pi.
\]
Consequently,
\[
\epsilon H(\bar\pi\mid\pi^n)
=\epsilon H(\bar\pi\mid\rho)+\int c\,d\bar\pi-P_n
=\mathcal E(\bar\pi)-P_n.
\]
As $\pi^n\Rightarrow\pi_0$ and $P_n\to V$,
$\epsilon H(\bar\pi\mid\pi_0)
\leq\liminf_n\epsilon H(\bar\pi\mid\pi^n)
=\mathcal E(\bar\pi)-V<\infty$.
Hence $\bar\pi\ll\pi_0$.
Together with $\bar r\geq\kappa>0$ $\rho$-a.e.,
this gives $\rho\ll\bar\pi\ll\pi_0$, and thus
$r_0>0$ $\rho$-a.e.

\medskip
\noindent\emph{Step 4.}
We have $\pi^n=r_n\rho\in\Pi(\mu,\eta)$,
$r_n\to r_0$ in $L^1(\rho)$, $r_0>0$ $\rho$-a.e.,
and $\mu(f_n)+\eta(g_n)=P_n\to V$.
Lemma~\ref{lem:limiting-potential-admissibility}
therefore applies. After passing to a subsequence and
applying affine normalizations, it yields
$f\in L^1(\mu)$, $g\in L^1(\eta)$, and a finite measurable
$h:\mathbb R^d\to\mathbb R^d$ such that
$g=\Gamma^{\epsilon}_{f,h}$ $\eta$-a.e.,
$(f,h)\in\mathcal D_{\operatorname{EMOT}}(\mu,\eta)$,
and $V=\mu(f)+\eta(g)$.

For completeness, weak duality follows from the entropy
variational inequality: for any
$(f,h)\in\mathcal D_{\operatorname{EMOT}}(\mu,\eta)$,
the definition of $\Gamma^{\epsilon}_{f,h}$ gives
\[
\int
\exp\biggl(
\frac{f(x)+\Gamma^{\epsilon}_{f,h}(z)+h(x)\cdot(z-x)-c(x,z)}
{\epsilon}
\biggr)\,d\rho=1.
\]
For every finite-entropy $\pi\in\Pi_M(\mu,\eta)$,
admissibility makes the exponent integrable under $\pi$
and gives $\int H(x)\cdot(z-x)\,d\pi=0$.
The entropy variational inequality thus yields
$\mu(f)+\eta(\Gamma^{\epsilon}_{f,h})\leq\mathcal E(\pi)$.
Taking the infimum over such $\pi$ proves weak duality.
Together with $V=\mu(f)+\eta(\Gamma^{\epsilon}_{f,h})$, this establishes
strong duality and dual attainment.
Setting $f_0:=f$ and $h_0:=h$ completes the proof.
\end{proof}

\end{document}